\documentclass[11pt]{amsart} 
\usepackage{amssymb,amsmath,amscd,amsfonts,a4,psfrag,graphicx,latexsym,epsfig}
\usepackage[active]{srcltx}

\usepackage[all]{xy}
\everymath{\displaystyle}

\newtheorem{teo}{Theorem}[section]
\newtheorem{lem}[teo]{Lemma}

\newtheorem{cor}[teo]{Corollary}

\newtheorem{prop}[teo]{Proposition}
\newtheorem{defi}[teo]{Definition}

\newtheorem{conj}[teo]{Conjecture}
\newtheorem{remark}[teo]{Remark}
\newtheorem{remarks}[teo]{Remarks}

\newcommand{\mc}{\mathbb{C}}
\newcommand{\mz}{\mathbb{Z}}

\newcommand{\Yy}{{\mathcal Y}}

\newcommand{\Bb}{{\mathcal B}}

\newcommand{\Dd}{{\mathcal D}}

\newcommand{\Ff}{{\mathcal F}}
\newcommand{\Gg}{{\mathcal G}}
\newcommand{\Hh}{{\mathcal H}}
\newcommand{\Ii}{{\mathcal I}}

\newcommand{\Ll}{{\mathcal L}}

\newcommand{\Nn}{{\mathcal N}}
\newcommand{\Oo}{{\mathcal O}}
\newcommand{\Pp}{{\mathcal P}}

\newcommand{\Rr}{{\mathcal R}}
\newcommand{\Ss}{{\mathcal S}}
\newcommand{\Tt}{{\mathcal T}}

\newcommand{\Ww}{{\mathcal W}}

\newcommand{\KG}{{\mathfrak K}}

\newcommand{\TG}{{\mathfrak T}}

\newcommand{\C}{{\mathbb C}}

\newcommand{\Z}{{\mathbb Z}}
\newcommand{\R}{{\mathbb R}}

\def\Dim{\emph{Proof. }}
\def\cvd{\hfill$\Box$}

\title[ ]{The Kashaev and quantum hyperbolic link invariants}

\author [ ]{St\'ephane Baseilhac$^1$ and Riccardo Benedetti$^2$}

\begin{document}
\today

\maketitle

\vspace{0.5cm}

\noindent $^1$ Universit\'e de Grenoble, Institut Joseph Fourier UMR
CNRS 5582, 100 rue des Maths, BP 74, F-38402 Saint-Martin-d'H\`eres
Cedex, FRANCE

\and 

\noindent Universit\'e Montpellier 2, Institut de Math\'ematiques et de Mod\'elisation, Case Courrier 51, Place Eug\`ene Bataillon, 34095 Montpellier Cedex 5, France (sbaseilh@math.univ-montp2.fr)

\smallskip

\noindent $^2$ Dipartimento di Matematica, Universit\`a di Pisa,
Largo Bruno Pontecorvo 5, 56127 Pisa, Italy (benedett@dm.unipi.it)\medskip

\begin{abstract}
We show that the link invariants derived from
$3$-dimensional quantum hyperbolic geometry can be defined by means of planar state sums based on link diagrams and a new family of enhanced Yang-Baxter operators (YBO) that we compute explicitly. By a local comparison of the respective YBO's we show that these invariants coincide with the Kashaev specializations of the colored Jones polynomials.  
As a further application we disprove a conjecture about the semi-classical limits of quantum hyperbolic partition functions, by showing that it conflicts with the existence of hyperbolic links that verify the volume conjecture.
\end{abstract}


\bigskip



\section{Introduction}\label{YBINTRO}
In this paper we describe the relationships between the following two sequences of complex valued invariants of links $L$ in the $3$-sphere:
\begin{enumerate}
\item the {\it Kashaev invariants} $<L>_n$, indexed by the integers
$n>1$ \cite{K2},
\item the {\it quantum hyperbolic invariants} $\Hh_N(L)$, indexed by
the {\it odd} integers $N>1$ and defined up to sign and multiplication by $N$th roots of unity \cite{Top}.
\end{enumerate}
We denote by $=_N$ the equality modulo such an ambiguity. We prove:
\begin{teo}\label{equivalenceHK} For every link $L$ and odd integer $N>1$ we have $<L>_N=_N\Hh_N(L)$.
\end{teo}
Due to some orientation conventions adopted in the present paper (see 
Remarks \ref{altrosegno} and \ref{sign-or}, and the remark after Theorem \ref{embedding}), we will actually prove that $<L>_N=_N \Hh_N(\overline L)$, where $\overline L$ denotes the mirror image
of the link $L$. This result puts on a solid ground the intersection of quantum hyperbolic geometry and colored Jones invariants, which are related respectively to $\Hh_N(L)$ and $<L>_N$, and based on different families of representations of the quantum group $U_{q}(sl_2)$. 

Following \cite{MM}, for each $n$ the Kashaev invariant $<*>_n$ can be defined by means of an enhanced Yang-Baxter operator including an R-matrix proposed by Kashaev in \cite{K2}. This R-matrix had been derived from the cyclic
representation theory of a Borel subalgebra $U_{\zeta_n}b$ of the quantum group
$U_{\zeta_n}(sl_2)$, where $\zeta_n = \exp(2\sqrt{-1}\pi/n)$ \cite[\S 6]{K3}. For every link $L$, $<L>_n$ is
computed by state sums based on planar link diagrams of $L$,
considered as the closure of a $(1,1)$-tangle. Surprisingly, Murakami-Murakami showed:
\begin{teo}\label{MM}\cite{MM}
For every link $L$ and integer $n>1$ we have $<L>_n = J'_n(L)$, the value at $q=\zeta_n$ of the colored Jones polynomial $J_n(L)\in \mz[q^{\pm 1}]$ normalized by $J_n(K_U)=1$ on the unknot $K_U$.
\end{teo}
The proof is by showing that the enhanced Yang-Baxter operator
of $<*>_n$ is {\it congruent} to the usual one of $J'_n(*)$, derived from the representation theory of
the {\it restricted} quantum group
$\overline{U}_{\zeta_n}(sl_2)$. Hence the corresponding
state sums take the same value on any given $(1,1)$-tangle
presentation of a link. Because the $n$-dimensional simple $\overline{U}_{\zeta_n}(sl_2)$-module $V_n$ has vanishing quantum dimension, $<*>_n$ vanishes on split links. Following Akutsu--Deguchi--Ohtsuki \cite{ADO}, we call {\it generalized
Alexander invariant} any link invariant constructed from an enhanced Yang-Baxter operator and having this property.

The quantum hyperbolic (QH) invariants
$\Hh_N(L)$ are
specializations to (see Section \ref{YBLINKINV} for details)
$$W=S^3, \ \ L=L^0, \ \ L_\Ff=\emptyset, \ \ \rho=\rho_{\rm triv}, \ \
\kappa=0$$ of invariants $\Hh_N(W,L_\Ff\cup L^0,
\rho,\kappa)$ defined in \cite{AGT} for compact closed oriented
$3$-manifolds $W$, where $L_\Ff\cup L^0$ is a link in $W$ made by a framed part
$L_\Ff$ and an unframed part $L^0$, $\rho$ is a $PSL(2,\C)$--valued character of $\pi_1(W\setminus L_\Ff)$, and
$\kappa$ is a collection, called cohomological weight, of
elements in the first cohomology groups of $W\setminus U(L_\Ff)$ and
$\partial U(L_\Ff)$, $U(L_\Ff)$ being a tubular neighbourhood of
$L_\Ff$ in $W$. For links in $S^3$ with $L_\Ff=\emptyset$, the character $\rho$ is necessarily the trivial one $\rho_{\rm triv}$ and $\kappa=0$. Each $\Hh_N(W,L_\Ff\cup L^0,
\rho,\kappa)$ is defined up to sign and multiplication by $N$th
roots of unity. It is computed by state sums $\Hh_N(\Tt)$ supported
by 3-dimensional pseudo--manifold triangulations $\Tt$ with additional structures
encoding $W$, $L_\Ff\cup L^0$, $\rho$ and $\kappa$, and made of tensors called {\it matrix dilogarithms},
associated to the tetrahedra of $\Tt$ and derived from the $6$j--symbols of the cyclic representations of the quantum group $U_{\zeta_N}(sl_2)$ \cite{Ba}. QH
invariants are defined also for cusped hyperbolic $3$-manifolds \cite{GT}.

By means of purely $3$-dimensional constructions we define in Section
\ref{YBO} a family of \emph{QH enhanced
Yang-Baxter operators} $({\rm R}_N,M_N,1,1)$ providing $\Hh_N(*)$ with planar state sums based on link diagrams. More precisely:
\begin{teo}\label{ATC} $\Hh_N(*)$ is the generalized Alexander invariant associated to $({\rm R}_N,M_N,1,1)$.
\end{teo}
The tensors ${\rm R}_N$ and $M_N$ are determined patterns of matrix dilogarithms where the dependence on the local parameters entering the triangulations $\Tt$ has been ruled out. We deduce Theorem \ref{equivalenceHK} from Theorem \ref{ATC} by a local comparison of enhanced Yang-Baxter operators.  Alltogether they give a $3$-dimensional existence proof and reconstruction of $<L>_N$, independent of the results of \cite{MM}. By the way, the Volume Conjecture \cite{K4,MM} is embedded in the general problem of the semi-classical asymptotics of QH invariants.

Theorem \ref{equivalenceHK} can be viewed as an unfolding in QH terms of \cite[Theorem 1]{K2}, which states that for every odd $N$ and every link $L$, $<L>_N$ can be computed up to multiplication by $N$th roots of unity by
certain 3-dimensional state sums $K_N(\TG)$ based on a
specific class $\TG$ of decorated triangulations of $S^3$ adapted to
$(1,1)$-tangle diagrams of $L$. By expanding
remarks of \cite{Top, GT, AGT}, we point out carefully in Section \ref{KvH} how the QH state sums both refine and generalize the 3-dimensional Kashaev's ones $K_N(\TG)$. In the case of links we find: 
\begin{prop} \label{eqKQHI} For every link $L$ and odd integer $N>1$ we have $\Hh_N(L)=_N K_N(\TG)$.
\end{prop}
As an application of Theorem \ref{equivalenceHK} and the existence of hyperbolic links verifying the Volume Conjecture, we disprove in Section \ref{AbS} a so called {\it asymptotics by
signatures} conjecture that would have predicted an attractive
general asymptotic behaviour of the QH state sums. All computations are collected in Section \ref{compute}.
\medskip

\noindent {\bf Notations.} In all the paper, for every integer $n>1$
we set $\zeta_n= \exp(2\sqrt{-1}\pi/n)$, or $\zeta$ when no confusion is possible, and we identify $\Ii_n = \{0,\dots,
n-1\}$ with $\Z/n\Z$ with its Abelian group structure. By
$\delta_n: \Ii_n \to \{0,1\}$ we mean the $n$-{\it periodic Kronecker
symbol}, satisfying $\delta_n(j)=1$ if $j=0$, and $\delta_n(j)=0$
otherwise. Odd integers bigger than $1$ will be denoted by $N$, and ``$=_N$'' means equality up to sign and multiplication by $N$th roots of unity.
\medskip

\noindent {\bf Acknowledgments.}  The first author's work was supported by the grant ANR-08-JCJC-0114-01 of the French Agence Nationale de la Recherche.

\section{Quantum hyperbolic link invariants}\label{YBLINKINV}
First we recall briefly some basic notions introduced in
\cite{Top, GT, AGT}. Then we will specialize them to the quantum
hyperbolic invariants of links in $S^3$.
\subsection{QH triangulated pseudo-manifolds and o-graphs}\label{o-gr}
A triangulated pseudo-manifold is a finite set of oriented, {\it
branched} tetrahedra $(\Delta,b)$, where the {\it branching} $b$
consists in edge orientations compatible with a total ordering of the
vertices of $\Delta$, together with orientation reversing face
pairings. We require that the quotient space $Z$ is a compact oriented
triangulated polyhedron with at most a finite set of non-manifold
points located at vertices, and that the local branchings match along
faces. Thus we have a {\it branched} (singular) {\it
triangulation} $(T,b)$ of $Z$ (an oriented {\it $\Delta$-complex} in the terminology of \cite{HATCHER}). By using the ambient orientation and the branching, every
tetrahedron can be given a $b$-{\it orientation}, whence a $b$-{\it
sign}, $*_b\in \{\pm 1 \}$.

For example, for a compact closed oriented $3$-manifold $W$ with a
link $L=L_\Ff\cup L^0$ as in the Introduction, the corresponding
pseudo-manifold $Z=Z(L_\Ff)$ is obtained by collapsing to one point
each component of $L_\Ff$; hence, if $L_\Ff=\emptyset$, then $Z=W$.
In the case of a cusped hyperbolic manifold $M$, $Z=Z(M)$ is obtained
by compactifying $M$ with a point at each cusp.
\smallskip

\begin{figure}[ht]
\begin{center}
 \includegraphics[width=9cm]{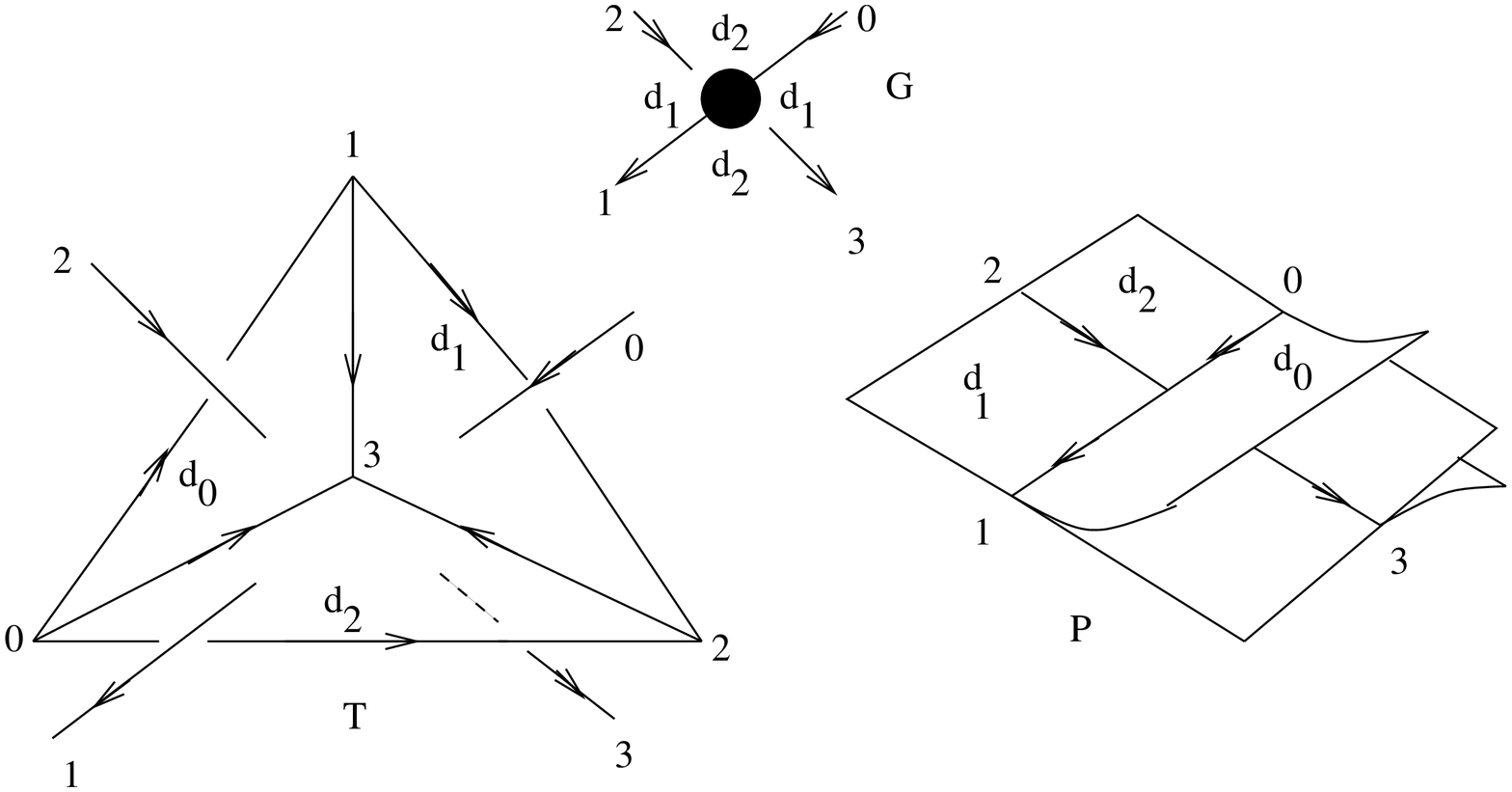}
\caption{\label{A-LoG+} $\Tt$, $\Pp$, and $\Gg$: $*_b=1$.} 
\end{center}
\end{figure}

\begin{figure}[ht]
\begin{center}
 \includegraphics[width=9cm]{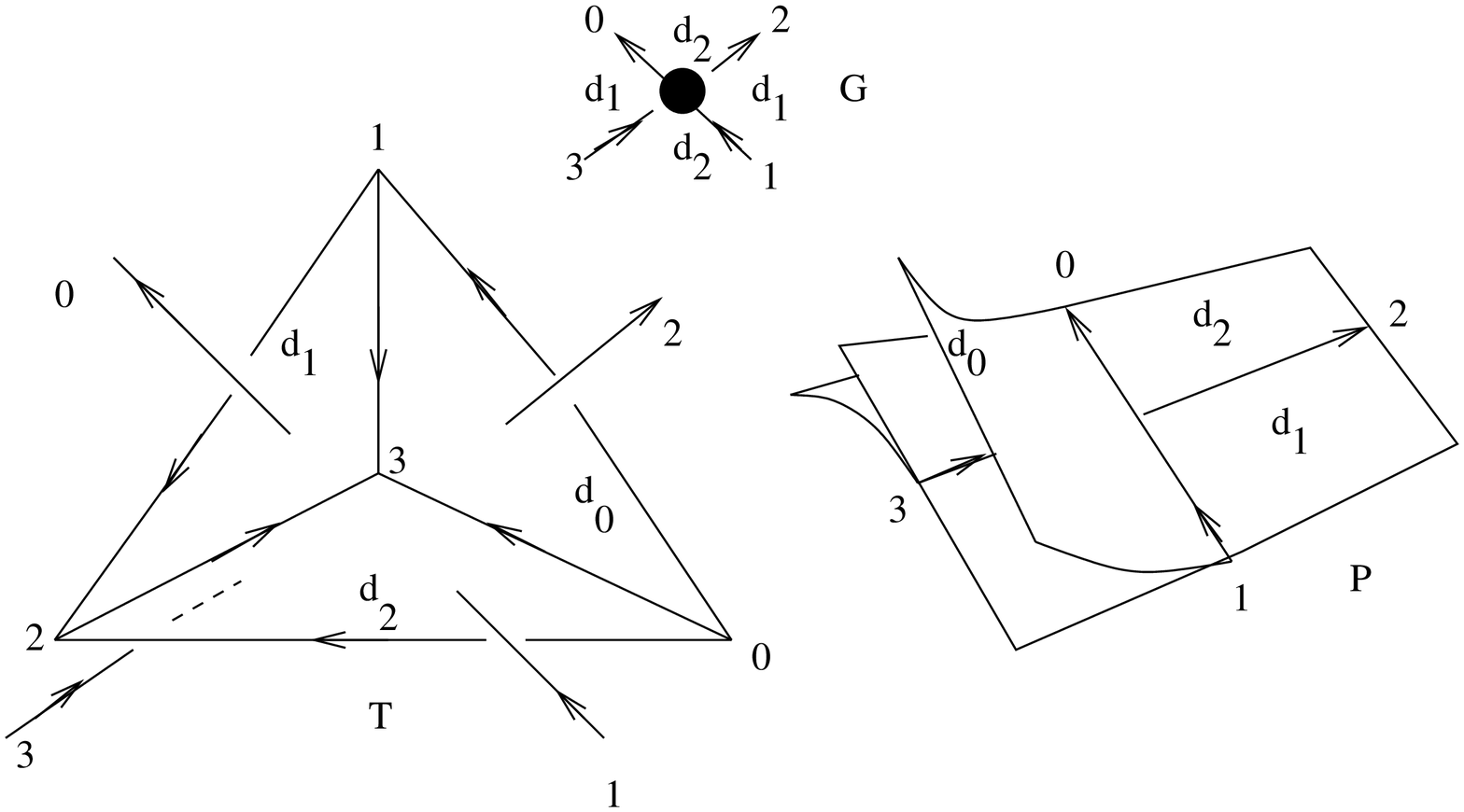}
\caption{\label{A-LoG-} $\Tt$, $\Pp$, and $\Gg$: $*_b=-1$.} 
\end{center}
\end{figure}

We have a QH triangulated pseudo-manifold $\Tt=(T,b,d)$ if
every tetrahedron $(\Delta,b)$ is equipped with
a {\it decoration} $d=(d_0,d_1,d_2)$ such that $d_j = (w_j,f_j,c_j)$ is attached to a pair of opposite edges of $(\Delta,b)$, the ordering of the $d_j$s is determined by the branching $b$ as in Figure
\ref{A-LoG+} and Figure \ref{A-LoG-}, and the following conditions (C1)-(C3) are satisfied:
\smallskip

(C1) $w_j \in \C\setminus \{0,1\}$, cyclically $w_{j+1}=(1-w_j)^{-1}$,
and $\textstyle \prod_j w_j = -1$; hence $w=(w_0,w_1,w_2)$
can be identified with the triple of cross ratio moduli of an ideal hyperbolic
tetrahedron;
\smallskip

(C2) $f_j \in \Z$ verify the {\it flattening} condition $\textstyle
\sum_j l_j =0$, where $ l_j := \log(w_j)+ f_j\sqrt{-1}\pi$ is called a
{\it log-branch} of $w_j$. Thus, if Im$(w_0)\geq 0$ (resp. $<0$), then
$(f_0,f_1,f_2)$ is a flattening iff $\textstyle \sum_j f_j = -1$
(resp. $\textstyle \sum_j f_j = 1$).
\smallskip

(C3) $c_j\in \Z$ verify the {\it charge} condition $\textstyle \sum_j
c_j = 1$.
\medskip

For every $N$, a decoration $d$ determines a system of {\it $N$th root
cross ratio moduli}
\begin{equation}\label{Nmoduli}
w'_j = (w'_N)_j =
\exp\left(\frac{\log(w_j)+\pi\sqrt{-1}(N+1)(f_j-*_jc_j)}{N}\right),\quad
j=1,2,3
\end{equation}
satisfying $\textstyle \prod_j w'_j = -\zeta_N^{-*_b(m+1)}$. 
\medskip

A QH triangulated pseudo-manifold $\Tt=(T,b,d)$ can be encoded by a
(normal) {\it QH o-graph} $\Gg=\Gg(\Tt)$, defined as follows \cite{BP}. The $2$-skeleton of the cell
decomposition dual to $T$ forms the {\it (standard) spine} $P=P(T)$ of
the complement of a regular neighbourhood of the vertices of
$T$. Every open $2$-cell of $P$, called a $2$-{\it region}, has the
orientation $\hat{b}$ dual to the $b$-orientation of the dual edge of
$T$. These region orientations define a {\it branching} of $P$. An o-graph
$G$ encoding $(T,b)$ is a suitable planar immersion with normal
crossings of the singular locus $S(P)$ of $P$. Every 2-face of $(T,b)$
has a prevailing {\it $b$-orientation} induced by the boundary edge orientations, thus $G$ is oriented in the dual way. $G$ has ``dotted'' crossings corresponding to the vertices
of $P$, and dual to the tetrahedra $\Delta$. The other ``virtual'' crossings of $G$ are immaterial. So $G$ encodes an
immersion in $\R^3$ of a branched regular neighbourhood $\Nn$ of
$S(P)$. It determines the whole branched spine
$(P,\hat{b})$ because every boundary component of $\Nn$ is filled uniquely by an oriented 2-disk. The QH o-graph $\Gg$ is defined by $G$ equipped with the decoration $d$ inherited from $(T,b,d)$.   

In {\rm Figure \ref{A-LoG+} and
Figure \ref{A-LoG-} we see a flat/charged branched tetrahedron
$\Delta(b,d)$, the neigbourhood $V(b,d)$ of the
corresponding vertex in $(P,\hat{b})$, the corresponding dotted crossing of the
o-graph $G$, and how the decoration $d$ transits to the $2$-regions of
$(P,\hat{b})$ and to the corners of the dotted crossing (in this case
$d_0$ is understood). Note that the pictures also indicates an ordering
$e_0,e_1,e_2$ of the edges opposite to the $3$-vertex. Sometimes we
denote the respective opposite edges by $e_j'$. Note that: 

(a) $\Delta$ is embedded in $\R^3$, with coordinates $(x,y,t)$ and the orientation determined by the standard basis, and inherits the induced orientation. The boundary is oriented by the
rule: {\it first the outgoing normal}. The $b$-orientation agrees with the boundary
orientation on two 2-faces of $\Delta$.

(b) $V(b,d)$ has 6 portions of oriented 2-regions of $P$.  Four
make the ``plate'' of $V(b,d)$, contained in the
$(x,y)$-plane with agreeing orientation. The two ``crests'' of $V(b,d)$
are over or down with respect to the $t$-coordinate. They are
oriented so that the o-graph is {\it left-turning}, that is,
its orientation coincides with the {\it prevailing} one among the
region boundary orientations, and the crests ``turn to the left'' with
respect to that orientation. $V$ is embedded in $\Delta$ as the
branched 2-skeleton of the dual cell decomposition, so that the plate goes onto a quadrilateral that cuts $\Delta$ by separating
the couples of vertices $(2,3)$ and $(0,1)$.  Note that the arc (resp. 2-region) orientation in $G$ (resp. $P$) is dual to the
$b$-orientation of the corresponding 2-face (resp. edge) of $\Delta$. The glueing rules of tetrahedra at common 2-faces (respecting all structures) can be read on $G$.

(c) Our convention for $*_b$-{\it signs} is such that it coincides at every dotted crossing of $\Gg$ with the
usual one at an oriented link diagram crossing.

\begin{remark}\label{altrosegno} {\rm From a ``simplicial'' point of view 
the opposite convention for $*_b$-signs is more
natural. We used it in \cite{Top,GT,AGT}, where we converted in a different way $(T,b)$ into a
planar graph that eventually supports the QH tensor networks. Here we follow the conventions of \cite{BP} so as to adopt an uniform sign rule for planar crossings of o-graphs and link diagrams. Both choices lead to QH theories isomorphic by reversing the orientation of QH pseudo-manifolds (ie. inverting the roles of $\Rr_N(\pm,d)$ below).}
\end{remark}  
The objects $\Tt$, $\Pp$ and $\Gg$ are equivalent
one to each other, and we use them indifferently. However, o-graphs are better suited when dealing with graphical encodings of tensor networks, as we will do all along the paper.

\subsection{QH state sums}
Given a QH triangulated pseudo-manifold $\Tt$, for every $N$ we
associate to every tetrahedron $(\Delta,b,d)$ of $\Tt$ (ie. to every
dotted crossing of the QH o-graph $\Gg$) the $N$-{\it matrix
dilogarithm} $\Rr_N(\Delta,b,d) = \Rr_N(*_b,d) \in {\rm Aut}(\C^N\otimes\C^N)$. More
precisely, the entries
$$\Rr_N(+,d)^{i,j}_{k,l}\ , \quad \Rr_N(-,d)_{i,j}^{k,l}$$ are
associated to the crossings of $\Gg$ with signs $*_b=+1$ and $*_b=-1$,
respectively, as on the left of Figure \ref{A-Rtens+} and Figure
\ref{A-Rtens-}, where $i,j,k,l \in \Ii_N$ label the edges of $\Gg$.

\begin{figure}[ht]
\begin{center}
 \includegraphics[width=7cm]{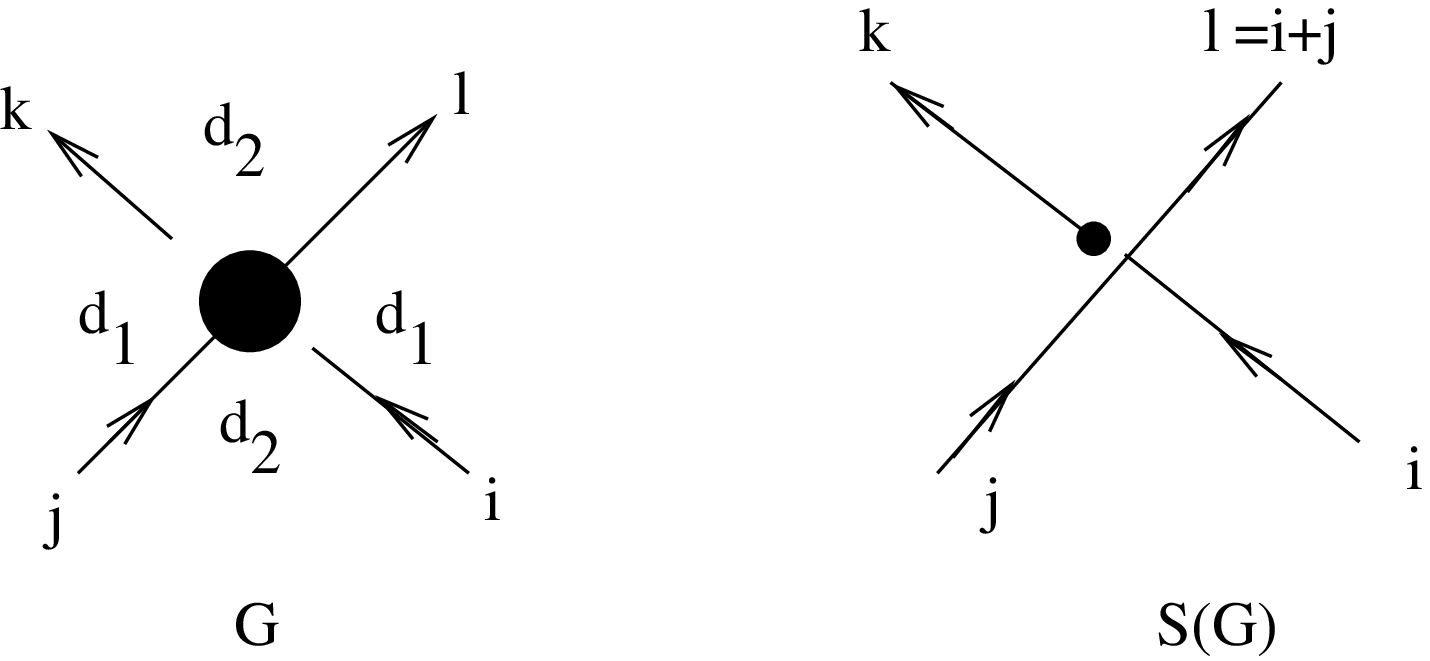}
\caption{\label{A-Rtens+}  Matrix dilogarithm and S-graph: $*_b=1$.} 
\end{center}
\end{figure}

\begin{figure}[ht]
\begin{center}
 \includegraphics[width=7cm]{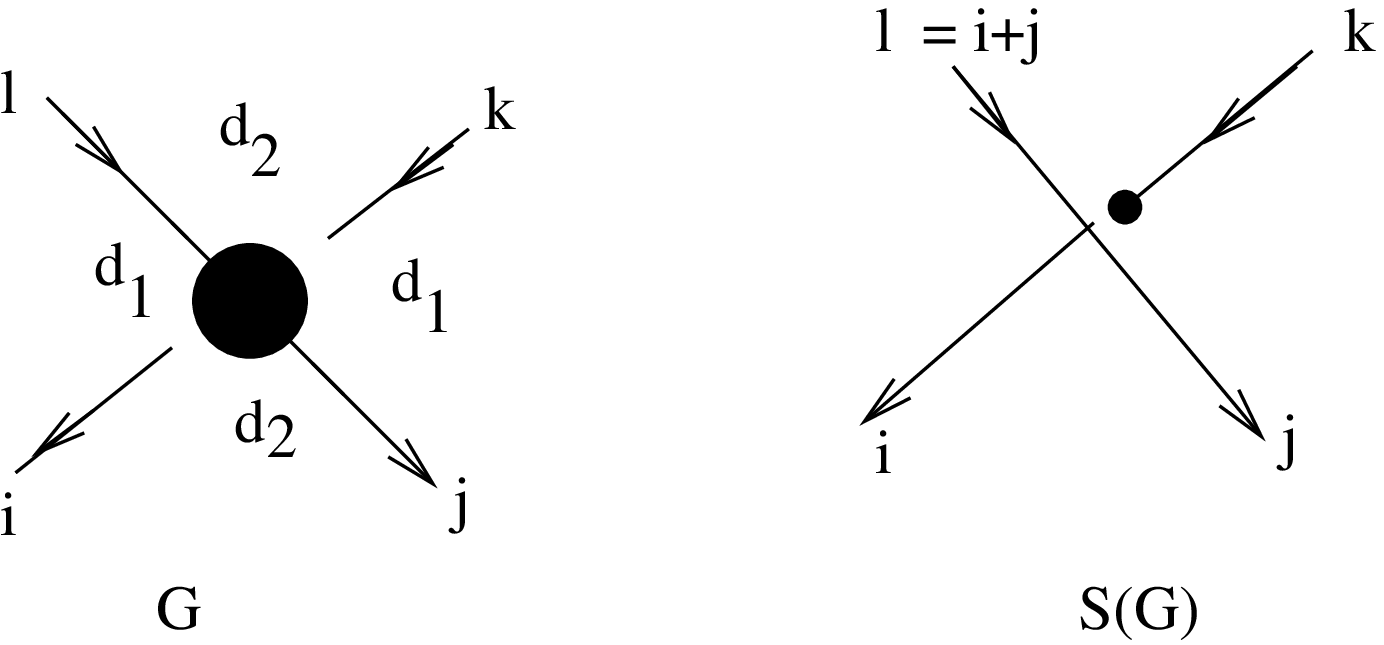}
\caption{\label{A-Rtens-}  Matrix dilogarithm and S-graph: $*_b=-1$.} 
\end{center}
\end{figure}

We call a {\it state} of $\Gg$ any labeling of its edges by indices in
$\Ii_N$. Every state $s$ selects an entry of $\Rr_N(\Delta,b,d)$, denoted by $\Rr_N(s,\Delta,b,d)$, at every crossing of $\Gg$. The QH
state sum $\Hh_N(\Tt)$ is defined by {\it tracing} (ie. contracting indices) the
resulting tensor network carried by $\Gg$:
\begin{equation}\label{ssum}
\Hh_N(\Tt)= N^{-(V-2)}\sum_s \prod_{(\Delta,b,d)}
\Rr_N(s,\Delta,b,d)
\end{equation} where $V$ is the number of vertices of $T$ that
are manifold points.

\begin{remarks}\label{diff-norm}{\rm 
(1) We will use {\it graphical} as well as {\it litteral} representations
    of tensors. Time by time, we
    must fix carefully the encoding/decoding rules in order to pass
    from one representation to the other. A first example is shown in Figure \ref{A-Rtens+} and Figure \ref{A-Rtens-}:
\medskip

\noindent {\bf Convention.} {\it The indices associated to ingoing (outgoing) arrows in a graphical representation correspond to top (bottom) indices in the litteral representation.}
\medskip

(2) In \cite{Top, GT, AGT} we used the normalization factor
$N^{-V}$. The present choice is more convenient to deal with the QH
link invariants, as it yields $\Hh_N(K_U)=1$ for every $N$ (see Lemma
\ref{no-crossing}).}
\end{remarks}

We refer to Section \ref{compute} for the explicit formulas of the
$N$-matrix dilogarithms. We just recall here that the non vanishing entries
$\Rr_N(s,\Delta,b,d)$ depend on the $N$th root
cross ratios $w'_0,w'_1$ given in (\ref{Nmoduli}), and correspond to
indices satisfying $i+j=l$ $(N)$.
\smallskip

Define $\Ss(\Gg,N)$ as the set of {\it efficient states}, such that
$\Rr_N(s,\Delta,b,d)\ne 0$ for all $\Delta$ in $T$. We have
$\Ss(\Gg,N) = H_1(S(\Gg),\partial S(\Gg);\ \Ii_N)$, where the {\it
S-graph} $S(\Gg)$ is the oriented (branched) graph with either
1-valent or 3-valent vertices, obtained from $\Gg$ by performing at
each vertex the modification shown on the right of Figure
\ref{A-Rtens+} and Figure \ref{A-Rtens-}; the 1-valent vertices form
the {\it boundary} $\partial S(\Gg)$. Hence $S(\Gg)$ determines the
actual range of summation in (\ref{ssum}), and governs the state sum $\Hh_N(\Tt)$. 

\subsection{ From link diagrams to $3$-dimensional triangulations}
\label{DIAG}
There is a very simple {\it tunnel construction}, introduced for
instance in Example 4.3 of \cite{Top}, that associates to a link
diagram a branched triangulation of $S^3$. It is very convenient to
describe this construction in terms of o-graphs.
\medskip

Let $\Dd$ be a link diagram on $\R^2\subset \R^2\cup \infty = S^2$,
representing a link $L$. We assume that $\Dd$ verifies the
following further condition (that can be always achieved for every
link $L$):
\smallskip

{\it Every connected component of $S^2\setminus |\Dd|$ is an open
$2$-disk, and $\Dd$ has at least one crossing.}
\smallskip
 
By $|\Dd|$ we mean the planar graph obtained by forgetting the
over/under crossings. The components of $S^2\setminus |\Dd|$ are called $\Dd$-{\it regions}.
\smallskip

Figure \ref{D-to-G} shows how to get from $\Dd$ two o-graphs $G'$ and
$G$ by replacing every crossing and every edge with an o-graph
portion. Both $G$ and $G'$ have no accidental virtual crossings.
\begin{figure}[ht]
\begin{center}
 \includegraphics[width=11cm]{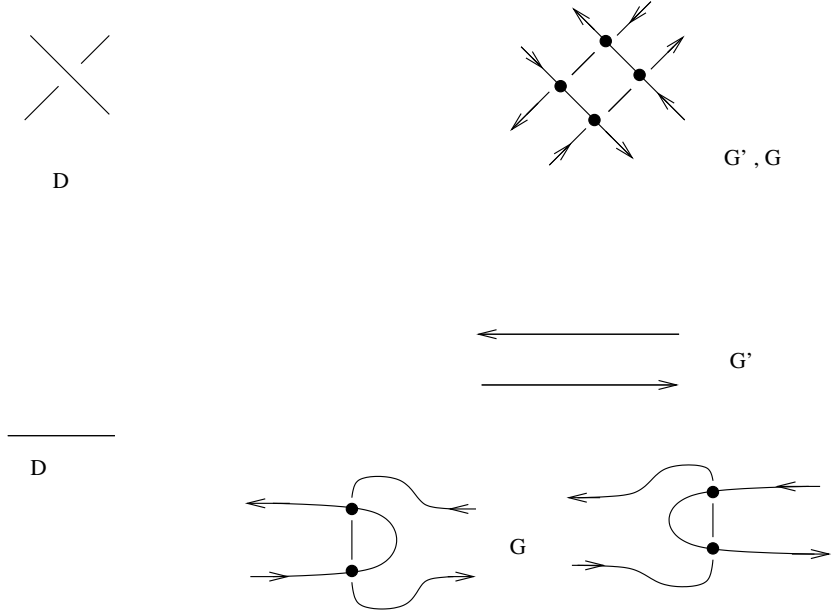}
\caption{\label{D-to-G} From link diagrams to o-graphs.} 
\end{center}
\end{figure} 

By forgetting the dots, the o-graph $G'$ appears as the
superposition of two oppositely oriented copies of the link diagram
$\Dd$.  It encodes a branched triangulation $(T',b')$ of the
pseudo-manifold $Z(L)$, with $4C$ tetrahedra ($C$ being the number of
crossings of $\Dd$), $E$ non-manifold vertices ($E$ being the number
of components of $L$), and $2$ further manifold vertices $V_\pm$. The
o-graph $G$ encodes a branched triangulation $(T,b)$ of $S^3$, with
$8C$ tetrahedra; the vertices $V_\pm$ persist in $T$, and there are
$2C$ further vertices.
\smallskip

Let us describe these triangulations. Start with
$$S^2 \times [-1,1] = S^3\setminus (B^3(-)\cup B^3(+))\subset S^3 $$
where $B^3(\pm)$ is an open $3$-ball with boundary $S^2\times \{\pm 1
\}$. Identify $\R^2 \subset \R^2\cup \infty = S^2$ with $S^2\times
\{0\}$, which is a non singular spine of $S^2 \times [-1,1]$. The
over/under crossings of the diagram $\Dd$ are thus specified with
respect to the coordinate $t$ on $[-1,1]$.

The o-graph $G'$ encodes a standard branched spine $(P',\hat{b}')$ of
$ S^2 \times [-1,1] \setminus U(L)$, where $U(L)$ is an open tubular
neighbourhood of $L$. The vertices $V_\pm$ correspond to the centers
of $B^3(\pm)$. In fact $P'$ is obtained by ``digging a tunnel'' in $
S^2 \times [-1,1]$ along $L$, by following the diagram $\Dd$. There is
a natural bijection between the $\Dd$-regions of the link diagram and
the $2$-regions of $P'$ contained in $S^2$, that we also call
$\Dd$-{\it regions}.

The o-graph $G$ encodes a standard branched spine $(P,\hat{b})$
obtained from $(P',\hat{b}')$ by inserting a {\it wall} in the tunnel
digged about every edge of $|\Dd|$; topologically, each wall is a
meridian $2$-disk of $U(L)$. In order to extend the branching we have
to fix the wall orientations. Both choices are admissible; at the
bottom of Figure \ref{D-to-G} we show the two possibilities for
$G$. Later we will fix a choice by using an auxiliary {\it diagram
orientation}. The vertices $V_\pm$ persist in the dual triangulation
$(T,b)$ and the $\Dd$-regions persist in
$(P,\hat{b})$. There is also a partition by pairs of the further $2C$
vertices of $T$ produced by the walls. Every pair, say $(v^-, v^+)$,
is associated to a crossing of $\Dd$. The two vertices of each pair
are separated by the spine $S^2\times\{0\}$ of $S^2 \times [-1,1]$,
and are the endpoints of an oriented edge $[v^-,v^+]$ of $(T,b)$.
\smallskip

The triangulation $(T,b)$ has the following further properties:
\smallskip

(P1) It is {\it quasi-regular}, that is, every edge of $T$ has
endpoints at distinct vertices.
\smallskip

(P2) The edges of $T$ dual to the walls realise $L$ as a subcomplex
$H'$ of the $1$-skeleton of $T$, containing all the further $2C$
vertices of $T$ but missing $V_\pm$.
\medskip

The definition of $G'$ works as well if $\Dd$ is the unknot diagram
without crossing; in such a case we stipulate that $G$ is obtained
from $G'$ by inserting two walls.

\begin{remark}\label{more-wall} {\rm 
If we insert several parallel oriented walls (more than one) around
every edge of $|\Dd|$, we get branched triangulations of $S^3$, with
more vertices, satisfying similar properties.}
\end{remark}

To simplify the figures sometimes we will indicate the o-graph $G$ by
means of {\it fat} diagrams, with black disks corresponding to the
walls of $G$ (the wall orientations will be specified time by
time). See Figure \ref{CARRIED}.

\subsection{Links carried by a link diagram}\label{A-LCARRIED}
Let $\Dd$, $(P,\hat b)$ and $(T,b)$ be as in Section \ref{DIAG}. We
indicate now two ways of selecting a {\it Hamiltonian subcomplex}
$H^0$ of the $1$-skeleton of $T$, that is, containing all the vertices
of $T$ (recall that in (P2) above, the subcomplex $H'$ realizing $L$
is not Hamiltonian since $V_\pm \notin H'$):
\medskip

(i) There is one $\Dd$-region of $P$, say $\Omega_0$, that contains
$\infty\in S^2$. Select a wall $B_0$ adjacent to $\Omega _0$.  Select two edges of $T$ dual to regions adjacent to $B_0$ and located at opposite sides of it, such that one has $V_+$ and the other $V_-$ as an endpoint. Remove from $H'$ the interior of the edge dual to $B_0$,
and take the union of the resulting triangulated arc with the two selected edges, and the edge dual to
$\Omega_0$. We get a complex $H^0$ that is Hamiltonian and
provides (up to isotopy) another realization of $L$. We denote it by
$L^0$.
\smallskip

(ii) Select two $\Dd$-regions of $G$. The dual edges have $V_+$ and
$V_-$ as endpoints, and their union realizes an unknot $K$ in $S^3$,
possibly linked with $L$. Take $H^0 = H' \cup K$. It realizes a link
$L^0=L\cup K$.
\begin{defi}\label{carried}{\rm For every link diagram $\Dd$, any link
$L^0$ obtained either as in (i) or (ii) is said to be {\it carried by $\Dd$}.
}
\end{defi}

In Figure \ref{CARRIED} we show some example of links carried by
diagrams. As indicated after Remark \ref{more-wall}, fat diagrams
correspond to o-graphs. The two distinguished regions
involved in the implementation of (ii) are labeled by ``$*$''. So in
case (a) $L=K_U$, the unknot, while $L\cup K$ is the {\it Hopf
link}. In case (b) again $L=K_U$, while $L\cup K$ is the link $4^2_1$,
according to the Rolfsen table. In case (c) we have $K_U$ versus the
{\it Whitehead link} $L_W$, and in case (d) the Hopf link versus the
link $6^3_1$ (the {\it chain link}). When a link is of the form $L^0=L\cup K$
where $K$ is an unknotted component, the procedure (ii) applied to a (suitable)
diagram of $L$ often produces the most economic triangulations of $S^3$ having
$L^0$ realized as a Hamiltonion subcomplex. 
\begin{figure}[ht]
\begin{center}
 \includegraphics[width=10cm]{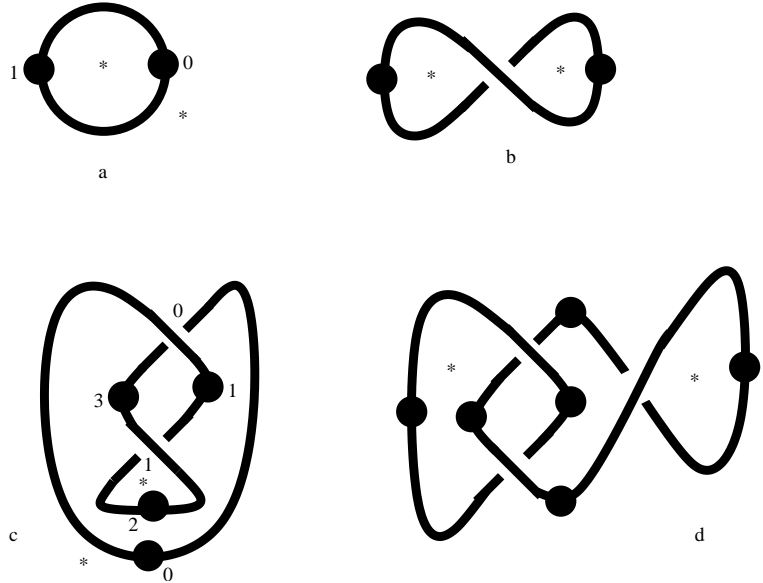}
\caption{\label{CARRIED} Some links carried by diagrams.}
\end{center}
\end{figure}


From now on we always denote by $L^0$ a link carried by a diagram
$\Dd$, hence obtained from $(T,b,H^0)$ as above.

\subsection{From link diagrams to distinguished QH triangulations}\label{DECO}
The next task is to convert $(T,b,H^0)$ into a {\it distinguished} QH
triangulation suited to the computation of the
quantum hyperbolic invariants $\Hh_N(L^0)$.
\medskip

Let $(T,b,d)$, $d=(w,f,c)$, be any QH triangulation supported by
$(T,b)$. We define the {\it total decoration} of an edge $e$ of $T$ as
$$d(e) = (W(e), L(e), C(e))$$ where $W(e)$ is the product over all
tetrahedra glued along $e$ of the $(w_j')^{*_b}$s at $e$ (the $N$th
root cross-ratio moduli or their inverses according to $*_b=\pm 1$), and
$L(e)$ and $C(e)$ are the similar sums of {\it signed} log-branches
$*_bl_j$ and $c_j$s, respectively. The total decoration $d(R)$ of a
spine $2$-region $R$ is defined in a dual way. We are going to impose
to $(T,b,d)$ global constraints in terms of $H^0$ and the $d(e)$s (see
\cite{Top, GT, AGT} for details).
\medskip

{\bf Global conditions on $(w,f)$.} These do not depend on $H^0$.
First we want the cross-ratios to encode the trivial representation
$\rho_{\rm triv}$. It is enough to require that
\begin{equation}\label{gcw}   
W(e)=1, \ \ {\rm for \ \ every\ edge}\ e.
\end{equation} 
Moreover, we also require that
\begin{equation}\label{gcl}  
L(e)=0,  \ \ {\rm for \ \ every\ edge}\ e.
\end{equation} 

{\bf Global charges.} The global conditions on charges encode the
subcomplex $H^0$ of $T$. They are
\begin{equation}\label{gcc}
\left\lbrace \begin{array}{ll}
C(e)=0 & {\rm for\ every} \ e\in H^0\\
C(e)=2 & {\rm for\ every} \ e\in T\setminus H^0.
\end{array}\right.
\end{equation}

We call {\it distinguished} any QH triangulation $\Tt=(T,H^0,b,d)$
satisfying the global constraints (\ref{gcw}), (\ref{gcl}) and
(\ref{gcc}). We denote by $\Gg=\Gg(\Dd,H^0)$ the corresponding QH
o-graph. It is a particular case of the general results of \cite{Top,
GT} that, up to the determined phase ambiguity, the value of the state
sum $\Hh_N(\Tt)$ does not depend on the choice of the distinguished triangulation $\Tt$
so that $\Hh_N(\Tt)$ well defines a link invariant $\Hh_N(L^0)$.
\medskip

Now we use some specific features of $(T,H^0,b)$ in order to specialize
the choice of $\Tt$.
\medskip

{\bf Universal constant system $(w,f)$.} A nice property of the triangulations $(T,b,H^0)$
is the existence of a {\it constant system} $(w,f)$ of cross ratios and flattenings that works for any
diagram $\Dd$ and any choice of wall orientations. A solution is
$$(w_0,f_0, f_1)=(2,0,-1), \quad
(l_0,l_1,l_2)=(\log(2),0,-\log(2)). $$ The conditions (\ref{gcw}) and
(\ref{gcl}) hold because along the boundary of every spine $2$-region
which is also a $\Dd$-region there is an even number of cross-ratios
$w_1=-1$, and at every further spine $2$-region there is a pattern of
pairs $(w_j,l_j)$ with opposite $b$-signs $*_b$.  Note that the same
argument works if instead of $(T,b)$ we take a triangulation obtained
by inserting an {\it odd} number of walls at every edge of $|\Dd|$
(see Remark \ref{more-wall}).
\medskip

\noindent {\bf Convention.} {\it From now on we will use by default the above universal constant system $(w,f)$, so that only the charges are varying parameters.}

\medskip

However, we will find useful in Section \ref{KvH} to use another, more general, way to make $(T,b,H^0)$ a distinguished QH triangulation.
\medskip

{\bf Idealization.} Systems of cross-ratios verifying (\ref{gcw}) can
be obtain as follows. We identify $(\C,+)$ with the subgroup of
$SL(2,\C)$ acting via Mo\"{e}bius transformations as translations on
$\C \subset \C\cup \infty = {\mathbb P}^1(\C)$.  The coboundary $z$ of
any $\C$-valued $0$-cochain on $(T,b)$ can thus be considered as a
$PSL(2,\C)$-valued $1$-cocycle on $(T,b)$ that represents $\rho_{\rm
triv}$. If $z$ is nowhere vanishing (which is the case if the
$0$-cochain is {\it injective}, since $T$ is quasi-regular), we say
that $z$ is {\it idealizable} (with base point $0$). In such a case,
if $x_0,x_1,x_2$ and $x_3$ are the vertices of a branched tetrahedron
$(\Delta,b)$ of $(T,b)$, ordered by using the branching $b$, we can
define four distinct points of $\C$ by
$$u_0=0,\ u_1=z([x_0,x_1])(0),\ u_2=z([x_0,x_2])(0),\
u_3=z([x_0,x_3])(0).$$ The associated cross-ratio is
$$w_0=[u_0:u_1:u_2:u_3] = u_3(u_2-u_1)/u_2(u_3-u_1)\in \mc\setminus
\{0,1\}.$$ 
Cross-ratios obtained in this way have so-called {\it
canonical log-branches}, which satisfy condition (\ref{gcl}):
\begin{equation}\label{can-flat}
\begin{array}{l}
l_0:=\log(u_2-u_1)+\log(u_3)-\log(u_2)-\log(u_3-u_1)\\
l_1:=\log(u_2)+\log(u_3-u_1)-\log(u_1)-\log(u_3-u_2)\\
l_2:=\log(u_3-u_2)+\log(u_1)-\log(u_3)-\log(u_2-u_1)
\end{array}
\end{equation}
The corresponding canonical flattenings are $f_j:=(l_j - \log(w_j))/\sqrt{-1}\pi$. Note that canonical flattenings are {\it even}, in the sense that
both $f_0$ and $f_1$ belong to $2\Z$.
\smallskip

It is easy to see that any constant cross ratio-system on $(T,b,H^0)$
can be obtained by idealization. For simplicity we will show it for
{\it knot} diagrams, the general case being not much harder. Assume at
first that the knot diagram $\Dd$ is {\it alternating}. We orient
every wall in such a way that the dual oriented edge is of the form
$[v^-,v^+]$, where the two endpoints are possibly associated to
different crossings of $\Dd$.  Note that this is possible because
$\Dd$ is alternating. Next we fix a 0-cochain $\gamma$ of the form:
$$\gamma (V_\pm)=\pm 1,\quad \gamma (v^\pm)=\pm a.$$ For a fixed
generic $a$ the idealization procedure gives, at every tetrahedron
$(\Delta,b)$ of $(T,b)$, the four points
$$ (u_0,u_1,u_2,u_3)=(0,a-1,a+1,2a)$$
with constant first cross-ratio
$$w_0= 4a/(a+1)^2 \ .$$  The corresponding canonical flattening is also constant. For example, $w_0=2$ iff
 $a=\pm\sqrt {-1}$;  if $a=\sqrt {-1}$, we
get the constant canonical flattening $(w_0,f_0, f_1)=(2,0,-2)$.
\medskip

If $\Dd$ is not alternating, we can modify the above procedure as follows
in order to obtain anyway any constant cross-ratio
$w_0= 4a/(a+1)^2$. 

\begin{lem}\label{solvecros} 
Given any knot diagram $\Dd$, there is a way to select a crossing segment at every double
point of $|\Dd|$ so that there is
exactly one segment endpoint on each edge of $|\Dd|$.
\end{lem}
\noindent {\it Proof.}  Orient $|\Dd|$. Select one crossing segment at a double point, and move along $|\Dd|$ according to the orientation. Pass across the next
visited double point without selecting any segment, continue and
select on the next visited double point the crossing segment according
to our move along the graph. Continue by alternating in this way:
``select'', ``pass across'', ``select'', ``pass across'', etc. If we
complete the circuit without obstructions we are done. Assume on the
contrary that we reach a first obstruction. This means that, for the
first time, either we visit again a double point with an already
selected crossing segment, and the rule would impose that we should
select now also the other crossing segment, or we visit again a double
point with no selected crossing segment, and the rule would impose
that we should again select no segment. In both situations we have
created a loop in $|\Dd|$ with an {\it odd} number of legs of
crossings pointing into the encircled region. So there is one leg that
is trapped. This is absurd.\cvd

\begin{remark}
{\rm Lemma \ref{solvecros} is obviously true for any alternating knot
diagram, and given a arbitrary knot diagram $\Dd$, if we define $\Dd'$ by stipulating that
every selected segment on $|\Dd|$ is over-crossing, then $\Dd'$ is
alternating. That is, Lemma \ref{solvecros} is equivalent to the fact
that for every knot diagram $\Dd$ there is an alternating knot diagram
$\Dd'$ such that $|\Dd|=|\Dd'|$.} \end{remark}

Now let $\Dd$ and $\Dd'$ be as in the last remark. Let $\gamma'$ be
the 0-cochain defined as above on the alternating diagram $\Dd'$; then, by
moving along $|\Dd|$ as in the proof of Lemma \ref{solvecros}, we
define a 0-cochain $\gamma$ on $\Dd$ such that $\gamma=\gamma'$ at
every crossing where $\Dd$ and $\Dd'$ coincide, and
$$\gamma (V_\pm)=\pm 1,\quad \gamma (v^\pm)= -\gamma'(v^\pm)$$
elsewhere. It turns out that the idealization of $\gamma$ has the constant
$w_0= 4a/(a+1)^2$, as desired. However, note that the canonical flattening is not constant.

\section{Enhanced QH Yang-Baxter operators}\label{YBO}
We deal with a triangulation $(T,b)$ associated to a diagram $\Dd$ of
a link $L$, according to the construction of Section \ref{DIAG}. It is endowed with the universal constant system $(w_0,f_0,f_1)=(2,0,-1)$, so
that it remains to manage with the charge.
\medskip

{\bf Notations for charge variables.} For the aim of future
computations, it is convenient to fix a name for the charge variables
on $\Gg$. In Figure \ref{SQUARE} we show the charge variables
$$(R,S,U,V,A,B,C,E,Y,Z,T,X)$$
at the four crossings of $\Gg$ corresponding to a crossing of
$\Dd$, and the charge variables
$$(P,F,H,M,G,K)$$ at the two crossings of $\Gg$ corresponding to a
wall. On the wall we have also indicated state variables $i,j,k,l$ for
future use. 
\begin{figure}[ht]
\begin{center}
 \includegraphics[width=9cm]{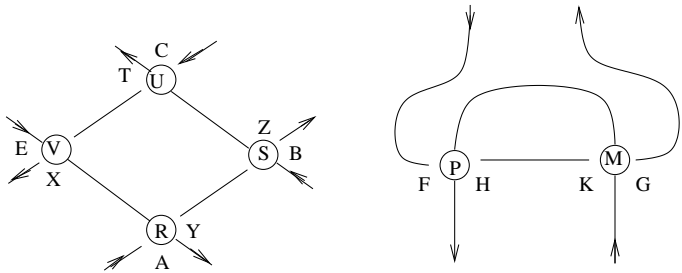}
\caption{\label{SQUARE} Crossing and wall charge variables.}
\end{center}
\end{figure}

Because of the symmetries there is no way for the moment to fix the
position of the variables $A$ rather than $C$, and so on.
\medskip

We refine the previous constructions by assuming that every link
diagram $\Dd$ is endowed with an {\it auxiliary orientation}.  We use
the orientation in order to:
\medskip

(O1) Fix the wall orientations according to the convention of Figure
\ref{WALL-OR}.
\begin{figure}[ht]
\begin{center}
 \includegraphics[width=9cm]{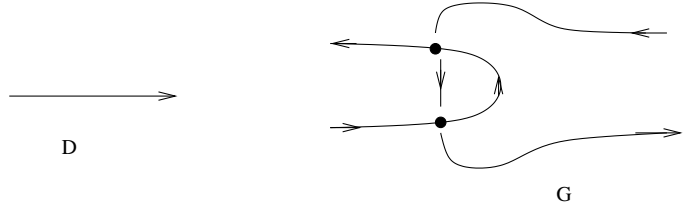}
\caption{\label{WALL-OR} Wall orientation.}
\end{center}
\end{figure}
\smallskip

(O2) Fix the notations of the charge variables at each crossing
according to the convention of Figure \ref{A-Lo-D}. Here we show only
the labelings of the four germs of $\Dd$-regions at a crossing. The
other labelings follow in agreement with Figure \ref{SQUARE}.
\begin{figure}[ht]
\begin{center}
 \includegraphics[width=11cm]{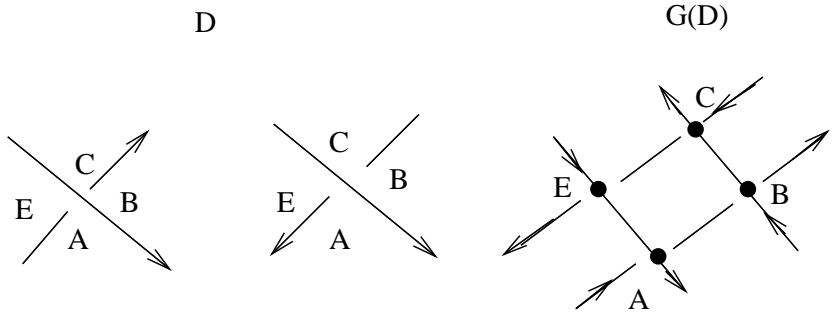}
\caption{\label{A-Lo-D} Charge labeling via oriented diagrams.}
\end{center}
\end{figure}
\smallskip

(O3) Fix a partition by pairs of the walls of the triangulation
$(T,b)$: each pair is made of two walls located at the {\it outgoing}
edges of a crossing of $\Dd$.
\medskip

\subsection{Specialization to closed braids and Yang-Baxter charges}
\label{YBCHARGE}
Suppose now that {\it $\Dd$ is the closure of a braid $\Bb$}. We
stipulate that:
\smallskip

\noindent {\bf Convention.} {\it braids are vertical and directed from bottom
to top, with the closing arcs on the right}. 

\smallskip

Hence $\Dd$ is oriented. Every closing arc has one maximum and one minimum with
respect to the vertical direction. 

Modify the triangulation $(T,b)$
associated to $\Dd$ by inserting a further wall at each such
maximum/minimum point, oriented by the rule of Figure \ref{WALL-OR};
these walls are indicated by {\it black squares} on fat diagrams (see
Figure \ref{curl}). As for $(T,b)$, the edges of the resulting
triangulation $(T',b')$, which are dual to the walls make a {\it non}
Hamiltonian subcomplex $H'$ realizing $L$. The universal constant
system $(w,f)$ works as well on $(T',b')$.
\smallskip

We are going to define a family of charges made up from {\it local
constant} pieces associated to the crossings and the walls on the
(oriented) fat diagrams of the o-graph $\Gg'$ corresponding to
$(T',b')$. By using the above orientation conventions (O1)--(O3) and
those of Figure \ref{SQUARE}, we denote the charge variables as
follows:
\begin{equation}\label{localvar}\left\lbrace\begin{array}{l}
(R1,S1,\dots,A1,B1,C1,E1)\ {\rm at\ every\ {\it positive}\ crossing};\\
(R2,S2,\dots,A2,B2,C2,E2)\ {\rm at\ every\ {\it negative}\ crossing};\\
(P,F,H,M,-F,K)\ {\rm at\ every\ wall\ associated\ to\ a\ crossing};\\
(P1,F1=-1, H1,M1,G1=1,K1)\ {\rm at\ every\ wall\ associated\ to\ a\
maximum};\\
(P2,F2=-1, H2,M2,G2=1,K2)\ {\rm at\ every\ wall\ associated\ to\ a\
minimum}.
\end{array}\right.\end{equation} 
Next we are going to impose invariance of (\ref{localvar}) under the
stabilization moves, a braid Reidemeister move III and a composition
of braid Reidemeister moves II. See the figures \ref{curl},
\ref{R-II} and \ref{R-III}.
\begin{figure}[ht]
\begin{center}
 \includegraphics[width=10cm]{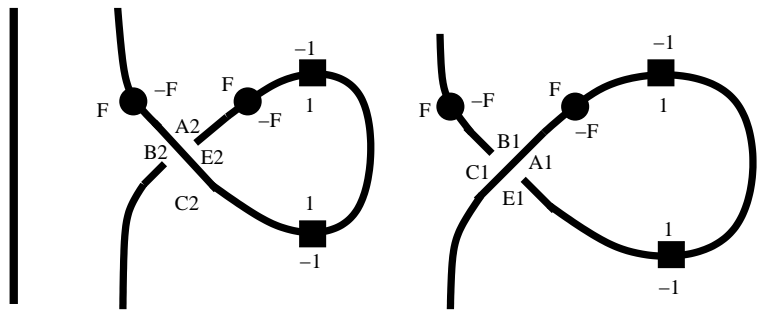}
\caption{\label{curl} Stabilization.}
\end{center}
\end{figure}

\begin{figure}[ht]
\begin{center}
 \includegraphics[width=10cm]{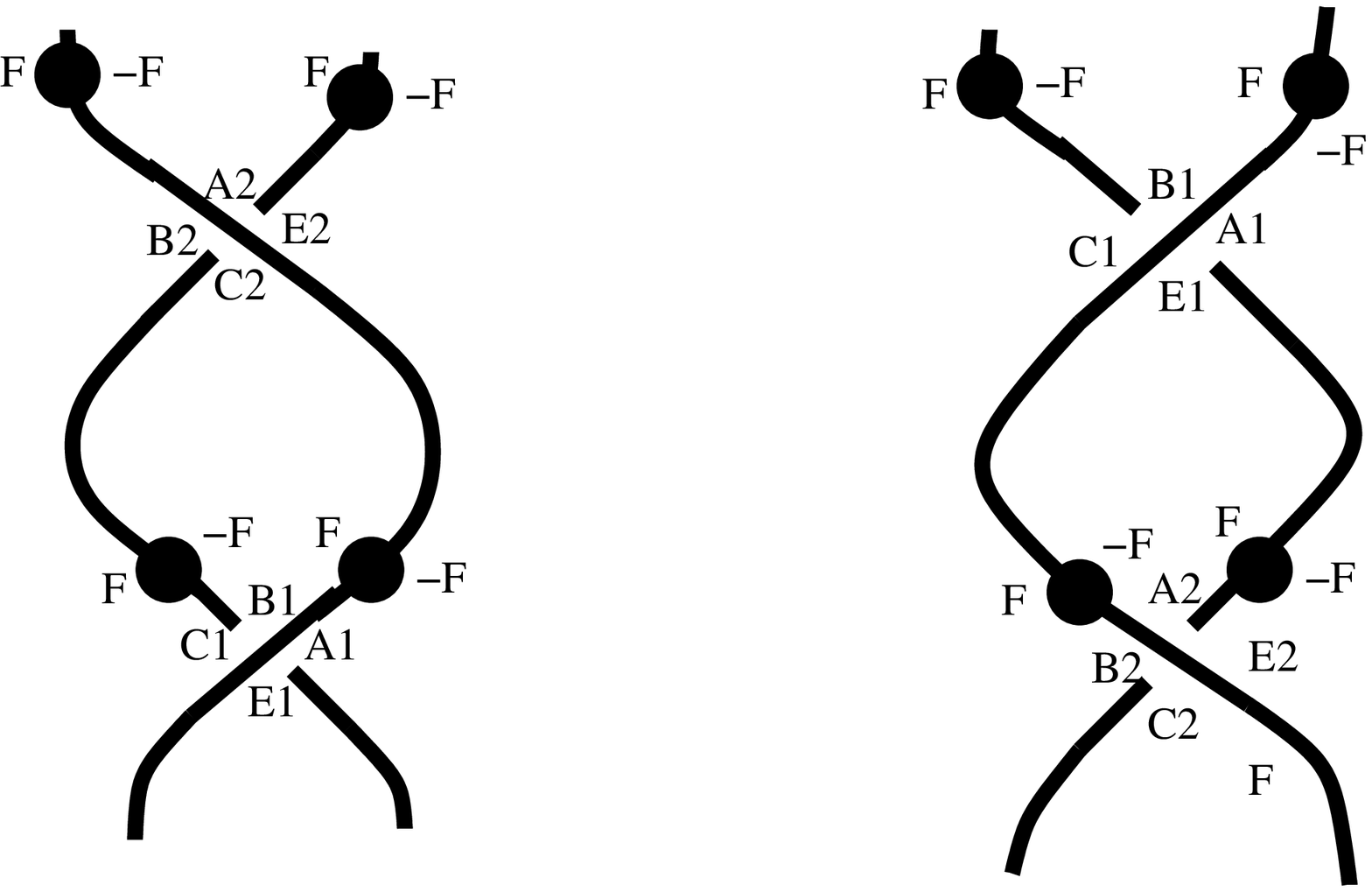}
\caption{\label{R-II} Braid Reidemeister move II.}
\end{center}
\end{figure}

\begin{figure}[ht]
\begin{center}
 \includegraphics[width=11cm]{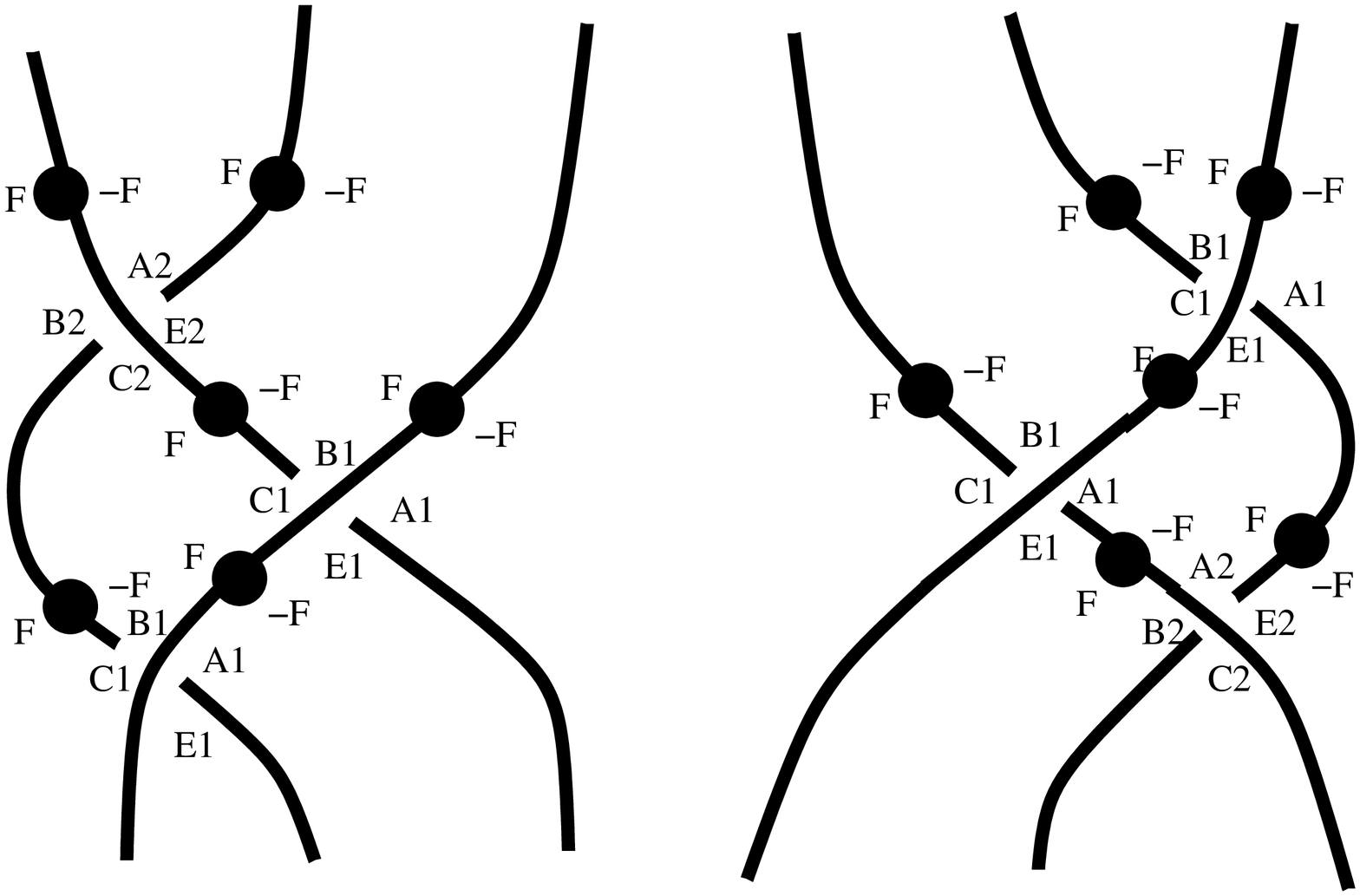}
\caption{\label{R-III} Braid Reidemeister move III.}
\end{center}
\end{figure}

In every case we have two portions of fat diagrams of QH o-graphs,
supporting portions of QH branched spines. The two portions of spines
carry sets of $2$-regions ``with boundary'' corresponding to edges of
$T'$ with incomplete star, which are in natural 1-1 correspondence one
to each other, and some internal $2$-regions, corresponding to edges
having complete star. The invariance of the charge variables (\ref{localvar}) should be a
consequence of {\it QH transits} relating the portions of QH branched
spines, in the sense of \cite{Top, GT, AGT}. Equivalently:
\begin{itemize}
\item Corresponding regions with
boundary have the same total charge; 
\item Each internal region has total charge equal to $2$. 
\end{itemize}

First we study these conditions on the $\Dd$-regions, that is, for the
charge variables that appear in the (planar) figures. Then we will see
how these lift to conditions on the whole charges.
\smallskip

The stabilization moves (i.e. the Reidemeister moves I) lead to the
set of independent conditions
\begin{equation}\label{lsys1}\begin{array}{llll}
B2+2F=0,& A2+B2+C2+E2=2, & C2-2+A2=0, \\
C1+2F=0,& A1+B1+C1+E1=2, & E1-2+B1=0.\end{array}
\end{equation}    
Note that the maxima/minima walls (having $F=-1$, $G=1$) contribute to
get a total charge equal to $2$ on the internal $\Dd$-region created
by the curl.
\smallskip

Consider a composition of two opposite Reidemeister moves II (recall
our convention on link diagrams in Section \ref{DIAG}). It gives a
further independent condition,
\begin{equation}\label{lsys2} B1+C2=2.\end{equation}
The Reidemeister move III yields, in turn, $F=0$. Then, the system
(\ref{lsys1}) \& (\ref{lsys2}) has the {\it one-parameter family of
solutions}
\begin{equation}\label{sol1}
A1=C1=0,\quad B1+E1=2, \quad B2=E2=0, \quad B1=A2,\quad E1=C2, \quad F=0.
\end{equation}
It works also for all other instances of Reidemeister moves III. 
\smallskip

Next we have to prove the existence of charges on $(T',H')$ satisfying
(\ref{sol1}). Recall the condition (C3) in Section \ref{o-gr}, that
will be always assumed. By using the above solution, at crossings we
have:
$$\begin{array}{cccc} X1=-1-V1+B1,& Y1=1-R1, & Z1=1-S1-B1, & T1=1-U1,\\
X2=1-V2, & Y2=1-R2-B1, & Z2= 1-S2, & T2= -1-U2+B1.
\end{array}$$
At walls we have
\begin{equation}\label{cond1}
\begin{array}{lll} H=1-P,& H1=2-P1,& H2=2-P2\\
K=1-M,& K1=-M1,& K2=-M2.\end{array} 
\end{equation}
The composition of Reidemeister moves II and the {\it positive} curl
introduce respectively the independent relations
$$\begin{array}{c}
T1+Z1=T1+X1=2-(P+M), \quad Y1+Z1=X1+Y1=P+M\\
T2+X2=X2+Y2=2-(P+M),\quad T2+Z2=Y2+Z2=P+M\\
U1+V1+S1+R1=0,\quad U2+V2+S2+R2=0,
\end{array}$$
and
$$\begin{array}{c} P+M+H2+K2=P1+M1+H+K = 2\\
P2+M2+H1+K1=P2+M2+X2+T2=2.
\end{array}$$
Together with (\ref{cond1}) the latter yields
$$P+M=H+K=P1+M1=H1+K1=P2+M2=H2+K2.$$ Finally we realize that neither
the {\it negative curl} nor the Reidemeister move III add independent
relations. Summing up the above computations we get:
\begin{prop}\label{constant-R} 
The variables (\ref{localvar}) lift to a global charge on the
triangulation $(T',H')$ which is invariant with respect to the
Reidemeister moves I and III and the composition of Reidemeister moves
of Figures \ref{R-II}, if and only the following relations, depending
on the free parameters $(U1,U2,B1,P,P1,P2)$, hold true:
$$\begin{array}{c} R1=U1, \quad S1=1-B1-U1, \quad V1=B1-U1-1\\
A1=0,\quad C1=0, \quad B1+E1=2\\
R2=2+U2-2B1, \quad S2=-1-U2+B1,\quad V2=S2\\
B2=0, \quad E2=0,\quad A2=B1, \quad C2=E1\\
F=0,\quad M=1-P,\quad G=0\\
F1=-1,\quad M1=1-P1,\quad G1=1\\
F2=-1,\quad M2=1-P2, \quad G2=1.\end{array}$$
\end{prop}
We call {\it Yang-Baxter charge} any of these solutions.

\subsection{From Yang-Baxter charges to QH link invariants}\label{YBCtoLINV}
Denote by $\Tt(\Bb,c)$ the QH triangulation of $S^3$ given by $(T',H',b')$, a fixed Yang-Baxter charge $c$, and the universal constant system $(w,f)$ of cross-ratio
moduli and flattenings. 

The QH triangulation $\Tt(\Bb,c)$ is {\rm not}
distinguished for $(S^3,L)$. In fact, every edge of the subcomplex
$H'$ realizing $L$ has total charge equal to $0$, while the other
edges have total charge equal to $2$, with the exception of the edge
dual to the region $\Omega_0$, which has total charge equal to
$-2$. Hence $\Hh_N(\Tt(\Bb,c))$ is {\it not} a state sum of the
quantum hyperbolic invariant $\Hh_N(L)$. However, by recalling the
constructions (i) and (ii) in Section \ref{A-LCARRIED}, there are two
natural ways to modify $\Tt(\Bb,c)$ in order to get QH link invariants.
The first is contained in the following lemma.
\begin{lem}\label{convert} 
By inserting two new walls after the minimum/maximun ones at the
closing arc adjacent to the region $\Omega_0$, the QH triangulation
$\Tt(\Bb,c)$ can be extended to a {\rm distinguished} one
$\Tt'(\Bb,c')$, such that $\Hh_N(L)=_N \Hh_N(\Tt'(\Bb,c'))$.
\end{lem} 
\Dim The universal constant system $(w,f)$ extends to the two new walls, so that it remains
to fix the charges on them. We do it as follows. The first wall
(according to the arc orientation) carries the same charges as black
disk walls at crossings; the second wall carries charges of the form
$$(P0,F0)=(P0,2), \quad (M0,G0)=(M0,0)$$
satisfying
$$P0+M0+1=0,\quad 1+H0+K0=2.$$ By the basic charge condition (3) of
Section \ref{o-gr} this reduces to $M0=-P0-1$. Any extension $c'$ of
$c$ is then determined by such a choice of $P_0$ and $M_0$. \cvd
\medskip

In Figure \ref{T2W} we represent both $\Tt(\Bb,c)$ (by forgetting the added
 dots near the region $\Omega_0$) and $\Tt'(\Bb,c')$ for a closed
 braid presentation of the Whitehead link. Values of a
 specific Yang-Baxter charge $c$ on the $\Dd$-regions are indicated ($B1=2$ and the $F=0$ values on the black disk walls at crossings are omitted).

\begin{figure}[ht]
\begin{center}
 \includegraphics[width=6cm]{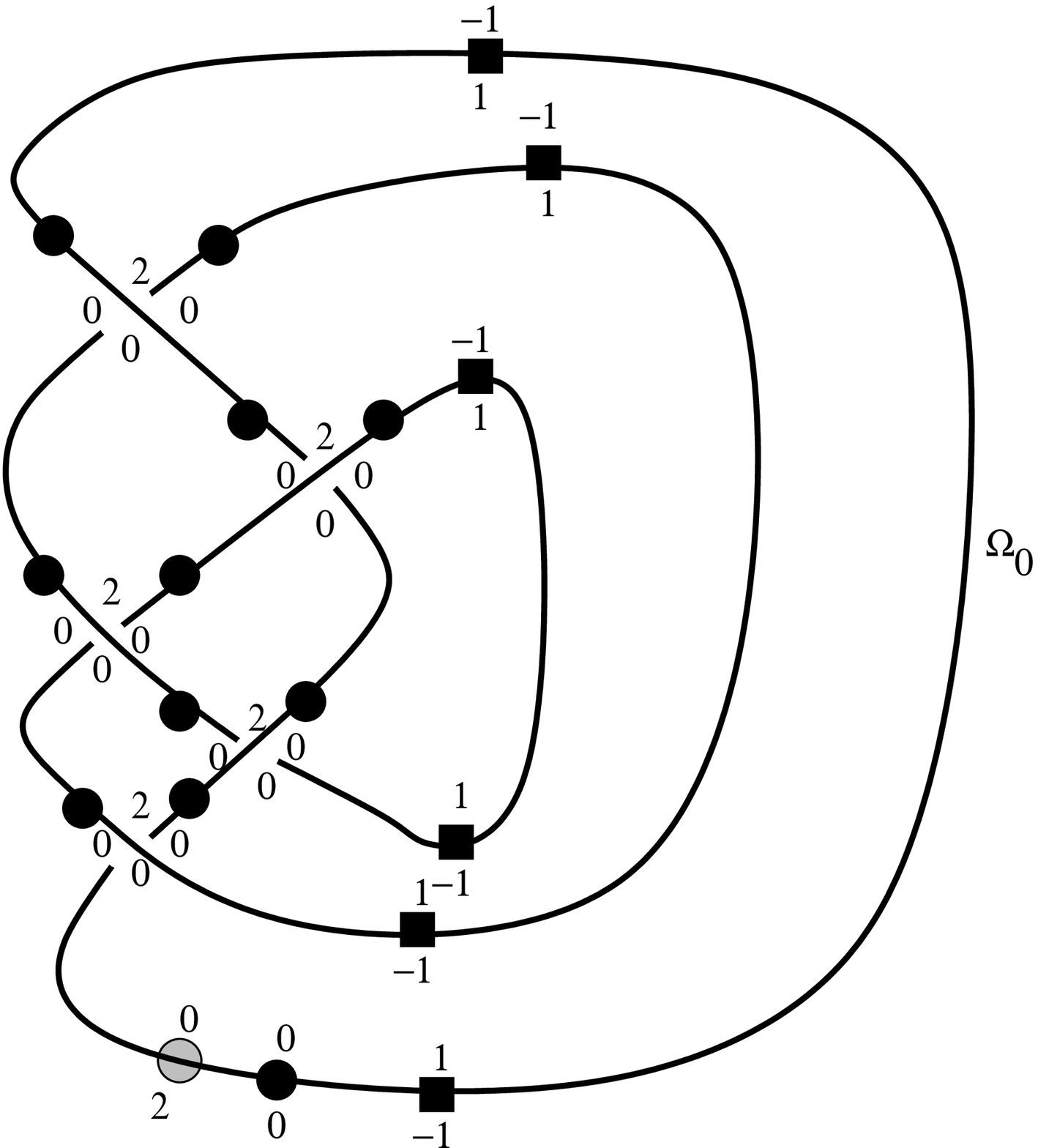}
\caption{\label{T2W} The distinguished QH triangulation $\Tt'(\Bb,c')$
for the Whitehead link.}
\end{center}
\end{figure}
\smallskip

The second way is contained in the following lemma.

\begin{lem}\label{convert2} By removing from $\Tt(\Bb,c)$ 
the maximum/minimum walls on the closing arc adjacent to the region
$\Omega_0$, we get a distinguished QH triangulation $\Tt''(\Bb,c'')$
($c''$ being the restriction of $c$) that carries the link $L^0= L\cup
K_m$, where $K_m$ is the meridian of the component of $L$ that contains
that arc. Hence $\Hh_N(\Tt''(\Bb,c''))=_N\Hh_N(L^0) =_N N\Hh_N(L)$. In
particular this does not depend on the choice of the component of $L$
supporting the meridian $K$.
\end{lem}
All the statements are clear with the exception of the last one, which
is a consequence of Lemma \ref{[*]-H} (2) below.

\begin{remark}\label{also 1-1} 
{\rm The above results hold as well when $\Dd$ is more generally the
closure of any oriented $(1,1)$-tangle diagram
in {\it normal position} with respect to the vertical direction (that is, all
crossings are directed from bottom to top like for braids). Every
$(1,1)$ tangle can be normalized by possibly rotating some crossing
and introducing maxima/minima. Braid presentations correspond to
special $(1,1)$-tangles obtained by re-opening the closing arc adjacent
to the region $\Omega_0$.}
\end{remark} 

Finally we can consider also the state sum $\Hh_N(\Tt(\Bb,c))$ itself.
The general invariance properties of QH state sums imply that its
value does not depend on the choice of the closed braid $\Dd$ and the Yang-Baxter charge $c$, up to the usual phase ambiguity. Hence we formally dispose of
further link invariants, say $$[L]_N =\Hh_N(\Tt(\Bb,c)).$$

\subsection{From Yang-Baxter charges to enhanced Yang-Baxter operators}
\label{ChargetoOp}
In order to finalize the discussion, we fix now the following specific
Yang-Baxter charge $c_0$:
$$R1=U1=0, \ S1=-1, \  V1=1, \ A1=0,\ C1=0, \ E1=0,  B1=2$$
$$R2=-2, \ U2=0,\ \ S2=1,\ V2=1, \ B2=0, \ E2=0,\ A2=2, \ C2=0$$
$$F=0, \ P=0, \ M=1,\ G=0$$
$$Fj=-1,\ Pj=0 \ Mj=1,\ Gj=1, \ \ j=1,2 \ . $$

Similarly, in Lemma \ref{convert} we fix $P0=0$. Note, however, that the following
discussion works as well for any Yang-Baxter charge $c$.  
\medskip

On the QH o-graph $\Gg$ corresponding to $\Tt(\Bb,c_0)$),
we point out a few distinguished {\it local configurations}:
\medskip

{\bf Walls:} There are two types of walls, either near a crossing or
at a maximum/minimum. We call them {\it C-wall} and {\it M-wall}
respectively.
\medskip

{\bf Crossings:} At every positive (negative) crossing we distinguish
two local configurations, called {\it braiding} and {\it complete crossing}
respectively.
\medskip

For every odd $N$, every such a local portion of $\Gg$ supports a
{\it QH tensor} in the following sense:

\begin{defi} The \emph{local QH tensor} of a portion of $\Gg$ is the result of tracing, like in formula \eqref{ssum}, the pattern of matrix dilogarithms associated to $\Gg$, and normalizing by a factor $N^{-1}$ for each wall (hence ``complete crossings'' below are normalized by a factor $N^{-2}$).\end{defi}

The local QH tensors can be read directly from the diagram $\Dd$, and
the normalization distributes the factor $N^{-(V-2)}$ in \eqref{ssum}.
\smallskip

The ``complete crossing'' local portions (in terms of
S-graphs), and the corresponding QH tensors (in {\it graphical}
representation) are shown in Figure \ref{QH-Rmatrix+} and Figure
\ref{QH-Rmatrix-}; all state variables belong to $\Ii_N$.

\begin{figure}[ht]
\begin{center}
 \includegraphics[width=11cm]{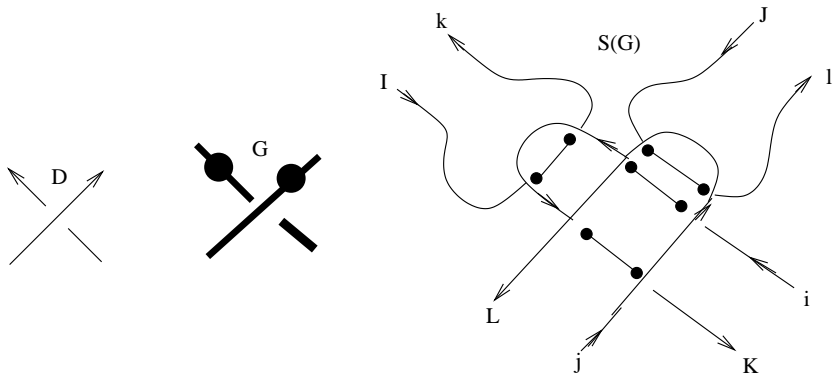}
\caption{\label{QH-Rmatrix+}  Positive complete crossing.} 
\end{center}
\end{figure} 
\begin{figure}[ht]
\begin{center}
 \includegraphics[width=11cm]{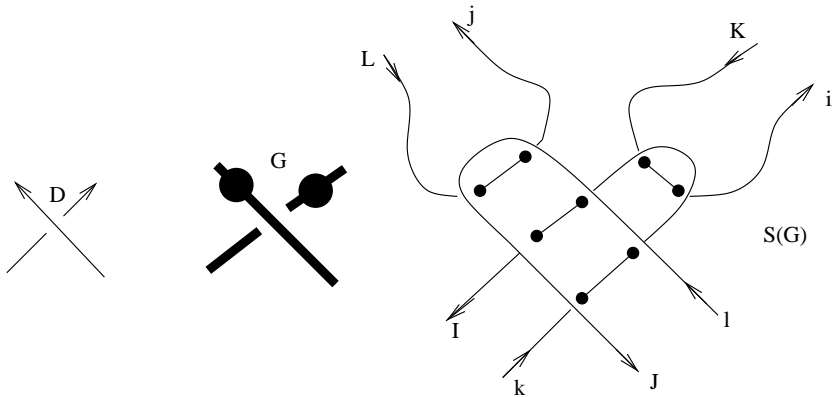}
\caption{\label{QH-Rmatrix-} Negative complete crossing.} 
\end{center}
\end{figure}

The position of the state indices is somehow
reminescent of the one for matrix dilogarithms in Figure
\ref{A-Rtens+} and Figure \ref{A-Rtens-}. The corresponding {\it
litteral} notation for the same tensors will be respectively
\begin{equation}\label{HNtens}
\Hh_N(C,+)^{i,j,I,J}_{k,l,K,L} \ \ ,\ \
\Hh_N(C,-)_{i,j,I,J}^{k,l,K,L} \ . 
\end{equation}
To get the ``braiding'' local portions, we just eliminate the
walls from Figure \ref{QH-Rmatrix+} and Figure \ref{QH-Rmatrix-}
(thus getting pictures like on the left of Figure \ref{SQUARE}), and we
keep the same state index distribution. The corresponding litteral
notations will be respectively
\begin{equation}
\Hh_N(B,+)^{i,j,I,J}_{k,l,K,L} \ \ ,\ \ \Hh_N(B,-)_{i,j,I,J}^{k,l,K,L} 
\ . 
\end{equation}
Concerning the walls, we name the state variables as on the right of
Figure \ref{SQUARE}. The corresponding litteral notation will be
\begin{equation}
\Hh_N(\Ww_X)^{j,k}_{i,l}
\end{equation}
where $X=C,M$ according to wall types. As the Yang-Baxter charge $c_0$
is fixed, these local QH tensors are constant at every
positive (resp. negative) crossing, and constant and equal at the
maxima/minima.
\medskip

{\bf Convention.} {\it We denote by ``$=_N$'' the equality of QH tensors modulo sign and multiplication by $N$th roots of unity that do not depend on states.}
\medskip

Let us anticipate some features of the {\it QH enhanced
Yang-Baxter operators} $({\rm R}_N,M_N,1,1)$ that we are going
to construct. They will include:
\medskip

(1) The  {\it QH R-matrix} ${\rm R}_N = {\rm R}_N(+)$ with entries ${\rm R}_N(+)_{l,k}^{i,j}$, $i, j, k, l \in \Ii_N$, associated to any positive crossing of $\Dd$ according to the graphical encoding on the left of Figure \ref{QH-Rmatrix4}, and similarly
for ${\rm R}_N(-)={\rm R}_N(+)^{-1}$, which is associated to any negative crossing. In Figure \ref{QH-Rmatrix4} we show also some
values of $c_0$.

\begin{figure}[ht]
\begin{center}
 \includegraphics[width=8cm]{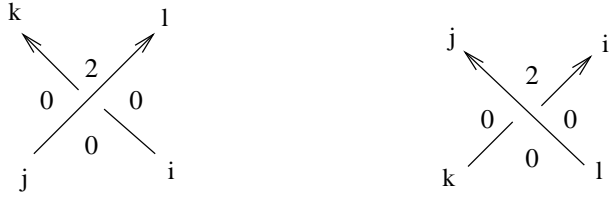}
\caption{\label{QH-Rmatrix4} QH R-matrix.} 
\end{center}
\end{figure} 

\smallskip

(2) An $N\times N$-matrix $M_N$ with entries $(M_N)^{k}_{i}$, $k,i \in \Ii_N$.
\medskip

We will construct $R_N$ and $M_N$ by using the above local QH tensors. However,
note that this is not so immediate as, for instance, the types of these
tensors are different. 
\medskip

{\bf Discrete Fourier transformation.} The first (main) modification
of the local QH tensors consists in replacing from the very beginning the
matrix dilogarithms $\Rr_N(*_b,d)$ by their (discrete)
{\it Fourier transform} $\tilde{\Rr}_N(*_b,d)$, as explained in
Section \ref{discreteFT}. Clearly, the value of $\Hh_N(L)$ is
unaltered by such a transformation. Thus every local QH tensor $\Hh_N(*)$
is replaced by the corresponding $\tilde{\Hh}_N(*)$.
\medskip

{\bf Conversion.} This is a purely formal manipulation of the local QH
tensors, producing tensors of different type. The idea is to
convert the local QH tensors into endomorphisms with source and target
given by the link orientation, in the spirit of quantum hyperbolic
field theory \cite{AGT}. Hence, for instance, the QH tensors of
crossings become endomorphisms of $(\mc^N \otimes \mc^N)^{\otimes 2}$
directed from bottom to top.  Moreover the litteral representations will be coherent with the current conventions adopted for ``planar''
R-matrices. We specify the conversion results by defining the entries.
\medskip

{\bf Complete crossing:}
\begin{equation}\label{Cr+}
{\rm Cr}_N(+)^{i,K,j,L}_{l,J,k,I} :=
\tilde{\Hh}_N(C,+)^{i,j,I,J}_{k,l,K,L}\quad ,\quad 
{\rm Cr}_N(-)^{l,J,k,I}_{i,K,j,L}:=
\tilde{\Hh}_N(C,-)_{i,j,I,J}^{k,l,K,L} \ .
\end{equation}
{\bf Braiding:}
\begin{equation}\label{Br+}
{\rm Br}_N(+)^{i,K,j,L}_{l,J,k,I} :=
\tilde{\Hh}_N(B,+)^{i,j,I,J}_{k,l,K,L}\quad ,\quad 
{\rm Br}_N(-)^{l,J,k,I}_{i,K,j,L}:=
\tilde{\Hh}_N(B,-)_{i,j,I,J}^{k,l,K,L} \ .
\end{equation}
{\bf Wall:}
\begin{equation}\label{XW}
({\rm XW}_N)^{k,l}_{i,j}:= \tilde{\Hh}_N(\Ww_X)^{j,k}_{i,l} \ .
\end{equation}
Note that
\begin{equation}
{\rm Cr}_N(\pm)= ({\rm CW}_N\otimes {\rm CW}_N)\circ {\rm Br}_N(\pm) \ .
\end{equation}
Finally, we can state our main step towards the construction of
the QH enhanced Yang-Baxter operators.
\begin{lem}\label{reducetensor} 
Denote by $V$ the ``diagonal'' subspace of $\mc^N\otimes \mc^N$ with
basis $e_i\otimes e_i$, where $i=0,\ldots,N-1$. Then:

1) The complete crossing tensors ${\rm Cr}_N(\pm)$ (resp. the
    wall tensors {\rm XW}) are supported by
$V\otimes V$ (resp. by $V$) and define automorphisms of it.

2) ${\rm Cr}_N(+)=_N {\rm Cr}_N(-)^{-1}$, providing we restrict the tensors to $V\otimes V$.
\end{lem}
\proof The first claim follows from Lemma \ref{braiding1} and
Corollary \ref{wall-conv} in Section \ref{compute}. Indeed, Lemma \ref{braiding1} shows that the braiding tensors ${\rm
Br}_N(\pm)$ are automorphisms of $(\mc^N \otimes \mc^N)^{\otimes 2}$
mapping $V^{\otimes 2}$ to itself. Also, Corollary \ref{wall-conv}
states that the wall tensors {\rm XW}, X=C, M, are endomorphisms of
$\mc^N \otimes \mc^N$ supported by $V$ and invertible on it. The
conclusion then follows from the fact that ${\rm Cr}_N(\pm)$
is obtained by composing braidings and walls.

For the second claim, consider the endomorphism $A$ of $V\otimes V$
supported by either the left or right member of Figure \ref{R-II}. 
Slide all walls to the top. Then, by applying one
Reidemeister move II at the middle of the figure we see that $A^2=_NA$;
such a move is done by reducing to QH transits as in the proof of
Proposition \ref{constant-R}. Since $A$ is invertible, $A=_N$ Id. \cvd\medskip

As $V$ is equipped with a given basis, it will be canonically
identified with $\C^N$. Define the endomorphisms ${\rm
R}_N(\pm)$ of $\C^N\otimes \C^N$ and ${\rm W}_{X,N}$ of $\C^N$ by
\begin{equation}
{\rm R}_N(+)^{i,j}_{l,k}:={\rm Cr}_N(+)^{i,i,j,j}_{l,l,k,k}\quad ,\quad 
{\rm R}_N(-)^{k,l}_{j,i}:= {\rm Cr}_N(-)^{k,k,l,l}_{j,j,i,i}
\end{equation}
\begin{equation}
({\rm W}_{X,N})^k_i := ({\rm XW}_N)^{k,k}_{i,i}.
\end{equation}
Finally set
\begin{equation}
M_N:= {\rm W}_{M,N}\circ {\rm W}_{M,N} \ .
\end{equation}
We define also the braidings restrictions:
\begin{equation}
{\rm B}_N(+)^{i,j}_{l,k}:={\rm Br}_N(+)^{i,i,j,j}_{l,l,k,k}\quad ,\quad 
{\rm B}_N(-)^{k,l}_{j,i}:= {\rm Br}_N(-)^{k,k,l,l}_{j,j,i,i} \ .
\end{equation}
 Thus ${\rm R}_N(\pm)$ corresponds to the {\it automorphism} of
$V\otimes V$ induced by ${\rm Cr}_N(\pm)$, and we have \begin{equation}\label{invres} {\rm R}_N(+)=_N {\rm R}_N(-)^{-1}\end{equation} according to the previous Lemma. Hence the local QH tensors are invariant under
the Reidemeister move II. Moreover, the QH tensors are invariant under QH
transits up to sign and multiplication by $N$th roots of unity (see
\cite{Top,GT,AGT}), and any two triangulations $\Tt(\Bb,c_0)$
differing by the Reidemeister moves of Proposition \ref{constant-R}
can be connected by using a finite sequence of QH transits. By
restricting to $V=\C^N$ we deduce:
\begin{itemize}
 \item ${\rm R}_N$ is an R-matrix, that is, we have the {\it
quantum Yang-Baxter equation}:
$$({\rm R}_N\otimes {\rm Id})({\rm Id} \otimes {\rm R}_N)({\rm R}_N \otimes {\rm Id})
=_N ({\rm Id} \otimes {\rm R}_N)({\rm R}_N \otimes {\rm Id})({\rm Id} \otimes {\rm
R}_N)$$
\item $M_N$ is an {\it enhancement} of ${\rm R}_N$, that is, we have
the identities:
$$(M_N \otimes M_N){\rm R}_N =_N {\rm R}_N(M_N \otimes M_N)$$
$$ {\rm Tr}_2({\rm R}_N^{\pm 1}(id \otimes M_N)) =_N {\rm Id} \ . $$
\end{itemize}
The commutation of $M_N^{\otimes 2}$ with ${\rm R}_N$ can be seen by
sliding the pairs of walls associated to consecutive maxima and minima
along the two strands of a positive crossing. The last identity
corresponds to the Reidemeister moves I.  Here, the partial
contraction
$${\rm Tr}_j: {\rm End}((\mc^N)^{\otimes k})\longrightarrow {\rm
End}((\mc^N)^{\otimes (k-1)}),\quad k\geq j\geq 1$$ is defined by
$${\rm Tr}_j(f)(v_{i_1}\otimes \ldots \otimes \widehat{v}_{i_j}
\otimes \ldots \otimes v_{i_{k}}) =
\sum_{j_1,\ldots,j,\ldots,j_{k}=1}^N
f^{j_1,\ldots,j,\ldots,j_{k}}_{i_1,\ldots,j,\ldots,i_{k}}\
v_{j_1}\otimes \ldots \otimes \widehat{v}_j \otimes \ldots \otimes
v_{j_{k}}$$ where $\textstyle f(v_{i_1}\otimes \ldots \otimes
v_{i_{k}})=\sum_{j_1,\ldots,j_{k}=1}^N
f^{j_1,\ldots,j_{k}}_{i_1,\ldots,i_{k}}\ v_{j_1}\otimes \ldots \otimes
v_{j_{k}}$ for a basis $\{v_i\}$ of $\mc^N$.
\smallskip

Summing up we have (see \cite{Tu}):
\begin{prop}\label{EYBO} 
The $4$-tuple $({\rm R}_N,M_N,1,1)$ is an enhanced
Yang-Baxter operator up to sign and multiplication by $N$th roots of unity.
\end{prop}
Explicit formulas are given in Section \ref{YB-esplicito}.
\smallskip

Let us come back to the situation of Section \ref{YBCtoLINV}. So $L$
is a link with a diagram $\Dd$ that is the closure of a braid $\Bb$,
say with $p$ strands. By composing elementary tensors of the form ${\rm Id} \otimes \ldots \otimes {\rm Id} \otimes R_N \otimes {\rm Id} \otimes \ldots \otimes {\rm Id}$ along the
braid one
gets a tensor
$$T_N(\Bb): V^{\otimes p} \rightarrow V^{\otimes p} \ .$$
Then
$$[L]_N =_N \Hh_N(\Tt(\Bb,c_0))=_N {\rm Trace}\left(M_N^{\otimes p}
\circ T_N(\Bb): V^{\otimes p} \rightarrow V^{\otimes p}\right).$$
Consider now the $(1,1)$ tangle diagram $\Dd_0$
of $L$ obtained by opening up the strand of $\Dd$ adjacent to the
region $\Omega_0$. The associated tensor is an endomorphism
\begin{equation}\label{QHtangle} T_N(\Dd_0): V\rightarrow
V\end{equation} that satisfies (recall Lemma \ref{convert2}) 
$$\Hh_N(L\cup K)=_N \Hh_N(\Tt''(\Bb,c_0'')) =_N {\rm
Trace}\left(T_N(\Dd_0): V \rightarrow V\right).$$ 



\begin{lem}\label{[*]-H}  For every odd $N>1$, we have:
\smallskip

(1) $[K_U]_N=_N 0$, $\Hh_N(K_U)=_N 1$, $\Hh_N(L_H)=_N N$
where $K_U$ is the unknot and $L_H$ the Hopf link.
\smallskip

(2) Let $L=L_1 \cup L_2$ be any split link (ie. $L_1$ and $L_2$ are
    unlinked). Then
$$ \Hh_N(L) = [L_1]_N\times \Hh_N(L_2) = \Hh_N(L_1)\times [L_2]_N \
. $$
\end{lem}
\Dim Statement (1) will be proved in Corollary \ref{no-crossing}. 
Statement (2) follows from the full
invariance with respect to Reidemeister moves (Proposition
\ref{constant-R}) and the fact that we can freely place the last walls
in the proof of Lemma \ref{convert} either at $L_1$ or $L_2$, without
effecting the value of $\Hh_N(L)$.  \cvd
\begin{cor}\label{split0} For every odd $N>1$, we have:
\smallskip

(1) $[L]_N=0$ for every link $L$.
\smallskip

(2) $\Hh_N(L)=0$ for every split link $L=L_1 \cup L_2$.
\end{cor}
\Dim Take $L'=L\cup K_U$, where $K_U$ is not linked with $L$. By Lemma
\ref{[*]-H} we have
$$ [L]_N=[L]_N\times \Hh_N(K_U)= \Hh_N(L)\times [K_U]_N =0$$
and $\Hh_N(L_1\cup L_2)= [L_1]_N\times \Hh_N(L_2) = 0$.\cvd\medskip

Hence we have proved (see Theorem \ref{ATC} in the Introduction): 

\begin{teo}\label{QHstate} 
The QH enhanced Yang-Baxter operators $({\rm R}_N,M_N,1,1)$ define
planar state sum formulas $\Hh_N(\Tt'(\Bb,c'))$ for the QH link invariants $\Hh_N(L)$, which
identify them as generalized Alexander invariants.
\end{teo}
\subsection{Puzzles}\label{R++}
Let $(T,b)$ be the triangulation associated to an oriented
link diagram $\Dd$, as at the beginning of this section. Grouping the
$C$-walls by pairs at each crossing has been a natural choice in order to have
the same {\it local} configurations. However, there are other
possible distributions of the $C$-walls that lead to the same
final link invariants having, for instance, some computational advantages.
\smallskip

Recall Lemma \ref{solvecros}: we can select at every crossing of $\Dd$ a
traversing segment (either over or under crossing) such that there is
exactly one segment endpoint on each edge of $|\Dd|$. Of course there in not a canonical
way to do it, and each way depends on the implementation of some {\it global} procedure.
\smallskip

{\bf Convention.} {\it Let us fix such a segment selection, and move every $C$-wall to
the corresponding segment end-point.}
\smallskip

This leads us to deal with new crossing tensors ${\rm R}_N(\epsilon_0,\epsilon_1)$ composed of two $C$-walls and one braiding, where $\epsilon_i = \pm$, $\epsilon_0$ is the
crossing sign and $\epsilon_1 =+$ if the selected
arc is over-crossing, and $\epsilon_1 =-$ otherwise. For instance:
$${\rm R}_N(-,+)= (id \otimes {\rm W}_{C,N}) \circ {\rm B}_N(-) \circ ({\rm
W}_{C,N}\otimes id)$$
that is 
$${\rm R}_N(-,+)^{r,i}_{k,s}= \sum_{j,l=0}^{N-1} ({\rm W}_{C,N})^l_s
 {\rm B}_N(-)^{j,i}_{k,l}({\rm W}_{C,N})^r_j \ .$$
By \eqref{invres} and the fact ${\rm W}_{C,N}^2
=_N {\rm Id}$ (Lemma \ref{wall-conv}), we have:
$${\rm R}_N(-)^{-1}=_N 
{\rm B}_N(-)^{-1}\circ ({\rm W}_{C,N}\otimes {\rm W}_{C,N})=_N ({\rm W}_{C,N}\otimes {\rm W}_{C,N})\circ {\rm B}_N(+)=_N {\rm R}_N(+)\ .$$
It follows that
$$ {\rm B}_N(+)=_N ({\rm W}_{C,N}\otimes {\rm W}_{C,N})\circ {\rm B}_N(-)^{-1}
\circ ({\rm W}_{C,N}\otimes {\rm W}_{C,N}).$$
Hence
$$ {\rm R}_N(+,+) =_N ({\rm id} \otimes {\rm W}_{C,N})\circ {\rm B}_N(-)^{-1} \circ
({\rm W}_{C,N}\otimes {\rm id} )$$
$${\rm R}_N(+,-)=_N ({\rm W}_{C,N}\otimes {\rm id} ) \circ {\rm B}_N(-)^{-1}\circ
({\rm id}  \otimes {\rm W}_{C,N})$$
and similarly for the others ${\rm R}_N(\epsilon_0,\epsilon_1)$.
\medskip

Different crossing tensors can be puzzled in order to produce QH link
partition functions in much more flexible way than the one strictly
suggested by the Yang-Baxter operator setup. In Figure \ref{puzzle} we show a few
examples of puzzles, that will be useful later.  The top left diagram computes the QH invariants of the
Whitehead link $L_W$ by:
\begin{equation}\label{LW}
\Hh_N(L_W)=_N \sum_{i,r,k,p=0}^{N-1} {\rm R}(-,+)^{r,i}_{k,i}  
{\rm R}(+,+)^{p,p}_{k,r} \ .
\end{equation}
On the top right we see the Hopf link:
\begin{equation}\label{LH}
\Hh_N(L_H)=_N \sum_{i,j,r,k,p=0}^{N-1} {\rm R}(+,+)^{k,r}_{i,i}  
{\rm R}(+,+)^{p,p}_{k,r} \ .
\end{equation}
Both $L_H$ and the link $4^2_1$ (see Figure \ref{CARRIED}) are carried by
the diagram with one crossing, so that we have also
$$\Hh_N(L_H)=_N \sum_{i,j=0}^{N-1} {\rm R}(-,+)^{j,j}_{i,i}$$
and
\begin{equation}\label{421}
\Hh_N(4^2_1)=_N \sum_{i,j=0}^{N-1} {\rm R}(+,+)^{j,i}_{j,i} \ .
\end{equation}
On the bottom we see the figure-eight knot $4_1$, with similar state
sums computing $\Hh_N(4_1\cup K_m)$ (recall Lemma \ref{convert2}), involving a few crossing tensors ${\rm R}_N(\epsilon_0,\epsilon_1)$.

\begin{figure}[ht]
\begin{center}
 \includegraphics[width=11cm]{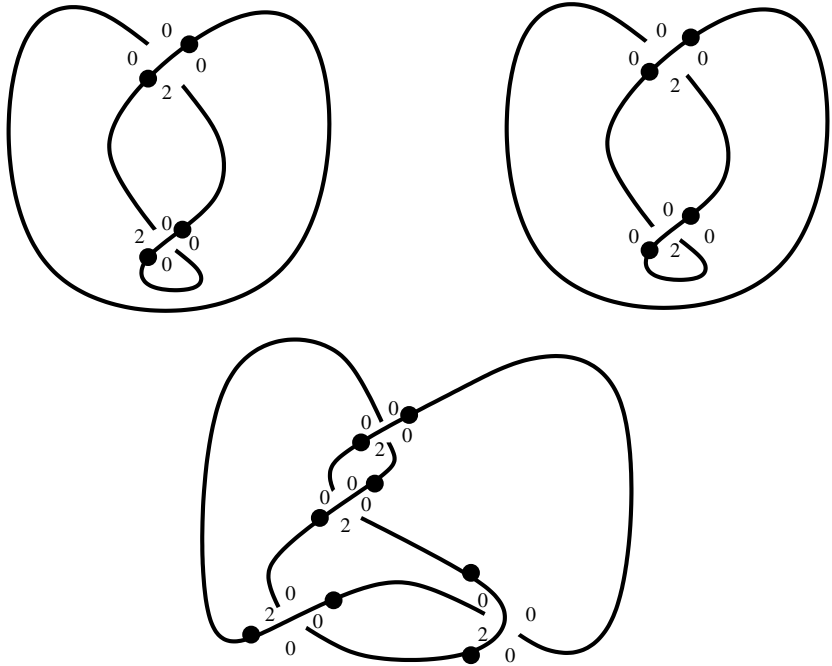}
\caption{\label{puzzle} Crossing tensor puzzle.}
\end{center}
\end{figure}
 
By using Lemma \ref{solvecros} it is not hard to prove the following:
\begin{prop}\label{QHpuzzle}
For every link diagram $\Dd=\Dd(L)$, there exists a branched triangulation
$(T,b)$ supporting a puzzle of QH crossing tensors
${\rm R}_N(\epsilon_0,\epsilon_1)$, whose total contraction gives
$\Hh_N(L\cup K_m)=_N N\Hh_N(L)$.
\end{prop}
Note that $\Dd$ is not oriented; we just use {\it local}
orientations to identify the tensors ${\rm R}_N(\epsilon_0,\epsilon_1)$;
these local orientations can conflict, but it does not matter. They select also
the wall orientations, hence the branching
$b$. We stress that $\Dd$ is not normalized in the sense of Remark \ref{also 1-1}, and has no $M$-walls. In other words, Proposition \ref{QHpuzzle} means that we
can puzzle local contributions of the Yang-Baxter charge $c_0$ in
order to produce a distinguished QH triangulation $\Tt(\Dd,c)=(T,b,c)$ such that
$\Hh_N((L\cup K_m)=_N \Hh_N(\Tt(\Dd,c))$. In practice, it is enough to
puzzle the tensors ${\rm R}_N(\epsilon_0,\epsilon_1)$ on the link diagram
in such a way that the $\Dd$-regions have the right total charge; then
there exists automatically a distinguished global charge.

\subsection{Proof of Theorem \ref{equivalenceHK}}\label{K=H} 
Given link $L$ consider a configuration
like in Lemma \ref{convert} (or Lemma \ref{convert2}) with our
favourite Yang-Baxter charge $c_0$. Apply Lemma \ref{solvecros}
to replace the QH R-matrices by a network of tensors ${\rm
R}_N(\sigma_0,\sigma_1)$. According to Proposition \ref{QHRmatrix} and Corollary \ref{wall-conv} (2), we have ${\rm
R}_N(+,+) = {\rm R}_N(+,-)$, and the following relation between $({\rm R}_N,M_N,1,1)$ and the Kashaev enhanced Yang-Baxter operator $(R_{K,n},\mu_{K,n},-s,1)$:
\begin{displaymath}\label{main-id}
 R_{K,N}(-)^{l,k}_{i,j} = \zeta^{1+(m+1)(l+k-i-j)}\
 {\rm R}_N(+,\pm)^{l,k}_{i,j}
\end{displaymath}
\begin{equation}\label{enhKMM}
 (\mu_{K,N})^j_i = \zeta^{m+1} M_N(c_0)^j_i,
\end{equation}
where $s=\exp(\sqrt{-1}\pi/n)$, and we note that $s=-\zeta^{m+1}$ when $n=N=2m+1$ is odd. A similar relation holds between  $R_{K,N}(+)$ and  ${\rm R}_N(-,\pm)$. When tracing the network of tensors ${\rm R}_N(\sigma_0,\sigma_1)$,
the factors $\zeta^{(m+1)(l+k-i-j)}$ compensate. Hence, by recalling
the notations 
Lemma \ref{convert} and \ref{convert2} we get:
\begin{teo}\label{embedding}  For every distinguished QH triangulation $\Tt'$ associated to a link
$L$ we have $\Hh_N(\Tt')=_N \ <\bar{L}>_N$ and $\Hh_N(\Tt'')=_N \ N<\bar{L}>_N$.
\end{teo}
The occurrence of the mirror image $\bar{L}$ depends on the orientation conventions we have adopted to define the quantum hyperbolic tetrahedra (see Remark \ref{altrosegno} and Remark
\ref{sign-or}), and the equality $\Hh_N(\bar{L}) \ =_N \ \overline{\Hh_N(L)}$ (see Proposition 5.8 of \cite{AGT}).

\section{QH and Kashaev's state sums}\label{KvH}
In this section we describe the
relations between the QH state sums and the 3-dimensional and planar state sums that had been proposed by Kashaev in \cite{K1} and \cite{K2}. 

\subsection{Generalities on the Kashaev's state sums} 
\subsubsection{The $3$-dimensional state sums.}\label{3DK} In \cite{K1}, for every odd $N>1$ a 3-dimensional state sum $K_N(\TG)$ is associated to any quasi-regular triangulation $(T,H)$ of any pair $(W,L)$, where $W$ is a compact
closed oriented $3$-manifold, $L$ is a link in $W$, and $H$ a
Hamiltonian subcomplex of $T$. The triangulation $(T,H)$ is equipped
with a decoration consisting of:
\smallskip

(i) A global charge $c$ on $(T,H)$, like in (\ref{gcc}) above; 
\smallskip

(ii) A total ordering of the vertices of $T$;
\smallskip

(iii) An injective $\C$-valued $0$-cochain $\gamma$ defined on the set
of vertices of $T$.
\smallskip

\noindent We denote by $\TG$ the resulting decorated triangulation.
\medskip

More precisely, it is required in \cite{K1} that $c$ takes half integer values and verifies half the global charge conditions mod($N$). Since $N$ is odd, it is
not restrictive to assume that $c$ lifts to an integral charge, like in the QH setup. 

It is not proved in \cite{K1} that the state sums $K_N(\TG)$ always
exist or define topological invariants of $(W,L)$, but rather that they are invariant under certain decorated versions of usual elementary triangulation moves.

\subsubsection{The planar state sums.} \label{PK} In \cite{K2}, a notion of {\it charged} oriented link diagram
$(\Dd,\hat{c})$ is introduced, together with a state sum $\KG_n(\Dd,\hat{c})$ for every $n>1$ ({\it not
necessarily odd}), involving an R-matrix (in the form of Boltzmann weights). Theorem 1 of \cite{K2} states that:
\begin{itemize}
 \item The planar state sums $\KG_n(\Dd,\hat{c})$ are invariant under charged
versions of the Reidemeister moves, and define link invariants
$<*>_n$. 
\item Decorated triangulations $\TG$ of
$(S^3,L)$ as in Section \ref{3DK} can be associated to certain charged diagrams $(\Dd,\hat{c})$, so that for every odd $N>1$, $K_N(\TG)$ is equal to $\KG_N(\Dd,\hat{c})$ up to
multiplication by $N$th roots of unity. From this it is claimed that the state sums $K_N(\TG)$ define invariants of links in $S^3$, up to the same ambiguity. 
\end{itemize}


\subsection{Relations between the $3$-dimensional and QH state sums} 

One can describe the state sums $K_N(\TG)$ of Section \ref{3DK} as follows. The vertex total
 ordering induces a branching $b$ by taking on every edge the
 orientation from the lowest to the biggest endpoint. Every
 tetrahedron $(\Delta,b)$ of $(T,b)$ is endowed as usual with a $*_b$-sign. The coboundary of  the $0$-cochain with respect to the edge $b$-orientation, $z=\delta
 \gamma$, is a nowhere vanishing $\C$-valued 1-cocycle. We denote by $z(e_0),z(e_1),z(e_2)$ the cocycle values on the
 edges of $\Delta$, named and ordered as in Figure \ref{A-LoG+} and Figure \ref{A-LoG-}. Recall that $e_j'$ denotes the opposite edge. Set
$$q_0=z(e_0)z(e_0'), \ q_1=z(e_1)z(e_1'), \ -q_2=z(e_2)z(e_2') .$$
By taking $z(e_j)^{1/N}z(e_j')^{1/N}$ for each of the $q_j$ we get a set of Nth roots $q_0'$, $q_1'$,
$-q'_2$ (see Section \ref{basic} for
our conventions on $z^{1/N}$). For every $N=2m+1$ and every charge $c$, one associates to $(\Delta,b,z,c)$ a tensor $T_N(\Delta,b,z,c)$ depending on $b$, $(q_j')_j$, and $c$. The 3-dimensional Kashaev's state sums have the form (note that Remark \ref{diff-norm} (2) applies also in this case)
\begin{equation}\label{Kformula}
K_N(\TG)= N^{-(V-2)}\prod_{e\in T\setminus H} (z(e)^{1/N})^{-2m}
\sum_{s}\prod_{(\Delta,b,z,c)} T_N(s,\Delta,b,z,c). 
\end{equation}
Recall the idealization procedure of Section \ref{DECO}. In \cite{Top} it is noted that (see also \cite{AGT}, where the role
of the canonical flattening is stressed): 
\begin{prop}\label{idealKT} The decorated triangulation $\TG$ defines a distinguished QH triangulation $\Tt$ with cross-ratio moduli $w_j=-q_{j+1}/q_{j+2}$ by taking the idealization of $z$ (considered as an $SL(2,\C)$-valued cocycle) and its canonical flattening (\ref{can-flat}). Moreover one has (see Section \ref{dilogs} for the definition of $\Ll_N$)
\begin{equation}\label{formT_N} T_N(\Delta,b,z,c)=_N (-q_2')^{\frac{N-1}{2}}(\Ll_N)^{*_b}(w_0',(w_1')^{-1}). 
\end{equation}
\end{prop} 
The tensors $T_N$ differ from the matrix dilogarithms $\Rr_N$ by the {\it local} normalization factor $(-q_2')^{\frac{N-1}{2}}$, instead of $((w_0')^{-c_1}(w_1')^{c_0})^{\frac{N-1}{2}}$. The {\it global}
normalization factor $$\prod_{e\in T\setminus H} (z(e)^{1/N})^{-2m}$$ occurring in (\ref{Kformula}) compensates the behaviour of the local ones in order to get the invariance with respect to the decorated triangulation moves. 
\smallskip

In \cite{Top} we showed that invariants $\Hh_N(W,L,\rho,\kappa)$ are defined for any $PSL(2,\C)$-valued character $\rho$ of the fundamental group of $W$ and any {\it cohomological weight} $\kappa\in H^1(W;\mz/2\mz)$, by means of the state sums $\Hh_N(\Tt)$ in \eqref{ssum}, based on distinguished QH triangulations. In the situation of Proposition \ref{idealKT}, $\Hh_N(\Tt)$ computes the invariant $\Hh_N(W,L,\rho_{{\rm
triv}},\kappa)$, the weight $\kappa$ being encoded by the flattening and the charge. The role of $\kappa$ is missed in \cite{K1}. However, by taking it into
account, the existence and invariance proof we have developed for
the QH invariants can be straightforwardly adapted in order to show:
\begin{cor} The state sums $K_N(\TG)$ compute invariants
$K_N(W,L,\kappa)$ well-defined up to multiplication by powers of $\zeta_N$ (with no further sign
ambiguity).\end{cor}
\begin{remark}\label{Borel}{\rm The definition of $K_N(W,L,\kappa)$ can be
extended to characters $\rho$ of $\pi_1(W)$ with values in a {\it Borel subgroup} of $PSL(2,\C)$ 
(See \cite[Remark 4.31]{Top}). On the other hand, unlike the tensors $T_N$ which depend on cocycle values (that is, the discretization of the
parallel transport associated to the flat connection corresponding to
$\rho$), the matrix
dilogarithms $\Rr_N$ entering the QH invariants depend on {\it cross-ratio moduli}, and fully display common structural features with the
classical Rogers dilogarithm \cite{GT}, which play a key role in order to develop a theory with nice analytic properties \cite{AGT}, dealing with arbitrary $PSL(2,\C)$-valued characters $\rho$ and arbitrary systems of cross-ratio moduli, possibly not arising from the idealization of $1$-cocycles (eg. for cusped hyperbolic 3-manifolds). The quantum coadjoint action \cite{Ba} explains the underlying relationship between the cyclic $6$j--symbols of a {\it Borel} subalgebra of the quantum group
$U_{\zeta_N}(sl_2)$, and those of the full quantum group itself.}
\end{remark}
Consider any situation where $K_N(W,L,\rho,
\kappa)=K_N(\TG)$ and $\Hh_N(W,L,\rho,\kappa)=\Hh_N(\Tt)$ are both defined. It follows from \eqref{formT_N} that
$K_N(\TG)$ and $\Hh_N(\Tt)$ have a common state dependent part formed by the entries of tensors $(\Ll_N)^{*_b}$, and differ by
non vanishing scalar factors $S_K$ and $S_{QH}$, respectively. Clearly, the ratio $S_K/S_{QH}$ is bounded
when $N\to +\infty$. Hence the invariants are ``asymptotically
equivalent'', that is:
\begin{cor}\label{asympt} $\limsup \{ \log|\Hh_N(W,L,\rho,\kappa)|/N \}\ =   \ 
\limsup \{ \log|K_N(W,L,\rho,\kappa)|/N\}$.
\end{cor}
In the case of the link invariants $K_N(L) := K_N(S^3,L,0)$ and $\Hh_N(L) = \Hh_N(S^3,L,\rho_{{\rm triv}},0)$, we can say even more:
\begin{prop}\label{H=K} 
For every link $L$ and odd integer $N>1$ we have $K_N(L)\ =_N \ \Hh_N(L)$.
\end{prop}
\noindent {\it Proof.} Consider a configuration like in Lemma \ref{convert}, with our favourite Yang-Baxter charge $c_0$. In order to compute both $K_N(\TG)$ and $\Hh_N(\Tt)$, take a 0-cochain having $a= \sqrt{-1}$ to realize the constant cross ratio system with $w_0=2$, and take the corresponding canonical flattening, as in Section \ref{DECO}. It
is enough to show that $S_{QH}/S_K\ =_N \ 1$. Multiply both scalar factors by $N^{(V-2)}$ and then consider rather $S_{QH}^{2/(N-1)}$ and $S_K^{2/(N-1)}$, which we still denote by $S_{QH}$ and $S_K$ for simplicity. As $w_1=-1$, we have $S_{QH}=_N\sqrt[N]{2}^{-2(C+1)}$. Now note that: 
\begin{itemize}
\item $z(e)= 2$ when the edge $e$ is dual to a $\Dd$-region;
\item $z(e) = \pm 2i$ when the edge $e$ is dual to a wall or a square-shaped region at a diagram crossing;
\item $z(e)=\pm \sqrt{2}e^{\pm \sqrt{-1}\pi/4}$ elsewhere;
\item $q_2=2$.
\end{itemize}
Denote by $C$ the number of crossings in the link diagram, $B$ the number of braid strands, $T$ the number of tetrahedra of the supporting triangulation, $R_D$ the number of $\Dd$-regions. Recalling how the hamiltonian subcomplex $H$ is featured in Lemma
\ref{convert}, we realize that
$$S_K=_N \sqrt[N]{2}^{T-2(R_D-1+C +3C +2B+1+1)}\ . $$
Since $T=8C+4B +4$ and $2=\chi(S^2)=C-2C+R_D$, we have also $S_K=_N \sqrt[N]{2}^{-2(C+1)}$, as desired. \cvd
\subsection{The planar state sums in the QH setup} Let $\Dd$ be any oriented diagram of a link $L$. The charged diagrams $(\Dd,\hat{c})$ considered in Section \ref{PK} are defined as follows. First assume
that $n=N=2m+1$ is odd. Let $(T,H,b)$ be a branched triangulation associated to $\Dd$ and carrying $L$, as in (i)
of Section \ref{A-LCARRIED}. Let $c$ be a global charge on
$(T,H)$. {\it Assume} that all the walls with total charge equal
to $0$ have charge values $F=G=0$, while the wall with total charge
equal to $2$ has $F=2$ and $G=0$, like in Lemma \ref{convert}. We define $\hat{c}$ as the labelling of the germs of $\Dd$-regions
at every crossing $v$ of $\Dd$ by {\it half} the residues mod$(N)$ of
the corresponding values of $c$ on the dual edges, that is, by variables $A'=[(m+1)A]_N\in \Ii_N$, and
similarly for $B'$, $C'$ and $E'$, as in Figure \ref{A-Lo-D}. Note that these
variables satisfy half the usual charge conditions mod($N$) about the vertices and faces
of $\vert \Dd\vert$. For arbitrary $n$, the labellings $\hat{c}$ are defined by variables in $\Ii_n$ satisfying the same conditions. 
\smallskip

Define an {\it $n$-state} of $(\Dd,\hat{c})$ as a labeling of every arc $e$
by an index in $\Ii_n=\{0,\dots,n-1\}$, such that the label is $0$ on
the edge adjacent to $\Omega_0$ and carrying the wall $B_0$. For every
$n$-state $s$, associate to every crossing $v$ of $(\Dd,c)$ with
crossing sign $\pm$ a {\it Boltzmann weight}
\begin{equation}\label{Boltz} 
R_n(\pm,v,s|\hat{c}) = R_n(\pm ,i,j,k,l\vert A',B',C',E')\end{equation}
according to Figure \ref{K-R} (see \cite[(2.8)]{K2}). The ``planar'' Kashaev state sums are then given by \cite[(3.6)]{K2}
\begin{equation}\label{Kss} 
\KG_n(\Dd,\hat{c})=\sum_s \prod_v R_n(\pm,v,s|\hat{c}) \prod_e
\zeta^{s(e)}. \end{equation}
\begin{figure}[ht]
\begin{center}
 \includegraphics[width=9cm]{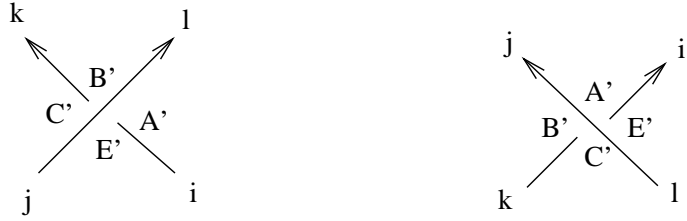}
\caption{\label{K-R} Graphical representation of Kashaev's Boltzmann weights.}
\end{center}
\end{figure}

An enhanced Yang-Baxter operator recovering the state sums \eqref{Kss} should include a {\it constant} R-matrix, whence a specialization of the variables $A'$,
$B'$, $C'$ and $E'$ to some fixed value. In \cite[(2.12) \& (2.15)]{K2} such a specialization is suggested. It is given by
\begin{equation}\label{cs1} A'=2\ ,\ B'=C'=E'=0 \quad \mbox{and}\quad B'=2\ , \ A'=C'=E'=0                                                                              
 \end{equation}                                                                                                             
 at negative and positive crossings, respectively. The corresponding R-matrix is given in terms of the Boltzmann weights \eqref{Boltz} by
\begin{equation}\label{Bw1}
(R_{K,n})^{l,k}_{i,j} = R_n(-,i,j,k,l\vert 1,0,0,0)\zeta^{k+l}
\end{equation}
and we have
\begin{equation}\label{Bw2}
(R_{K,n}^{-1})_{l,k}^{i,j} = R_n(+,i,j,k,l\vert 0,1,0,0)\zeta^{i+j}.
\end{equation}
(See Section \ref{KvH3} for explicit formulas). Now, to make the connection with the planar QH state sums, note that if $\Dd$ is the closure of a braid diagram $\Bb$ and $n=N$ is odd, the specialization \eqref{cs1} is induced by our favourite Yang-Baxter charge $c_0$, extended as in Section \ref{YBCtoLINV} to the triangulation $\Tt'(\Bb,c_0')$ obtained from $(T,H)$ by glueing the maximum/minimum walls. More precisely, by removing these walls we can turn $c_0'$ into labellings $\hat{c}_0$ of $\Bb$ that differ from \eqref{cs1} only at some determined top or bottom crossings, where the Boltzmann weights can be computed by tracing $R_{K,n}^{\pm 1}$ with an enhancing homomorphism $\mu_{K,n}$ (compare eg. with \cite[(2.17)]{K2}). By collecting terms in the state sum $\Hh_N(\Tt')$ of Theorem \ref{embedding}, it is not hard to check that: 
\begin{prop}\label{KpQH} For every odd $N>1$, every link $L$, every braid diagram $\Bb$ of $L$, and every distinguished QH triangulation $\Tt'=\Tt'(\Bb,c_0')$, we have $\KG_n(\Bb,\hat{c}_0) =_N \Hh_N(\Tt')$. In particular, $(R_{K,n},\mu_{K,n},-s,1)$ is an enhanced Yang-Baxter operator for the state sums $\KG_n(\Dd,\hat{c})$, which thus well define the Kashaev's link invariant $<L>_n$ for every charged oriented link diagram $(\Dd,\hat{c})$.
\end{prop}
Hence the Kashaev state sums $\KG_n(\Dd,\hat{c})$ compute $<L>_n$ by using arbitrary charged oriented link diagrams $(\Dd,\hat{c})$ (not necessarily associated to braid closures), in a way similar to the puzzles of Section \ref{R++}. 


\begin{cor}\label{Kteo1} Let $\TG = (T,H,c,\gamma)$ be a decorated triangulation as in Section \ref{3DK}, associated to a diagram $\Dd$ of a link $L$, and $(\Dd,\hat{c})$ a charged diagram associated to $(T,H,c)$. We have
\begin{equation}\label{teo1}
K_N(L) = K_N(\TG)= \KG_N(\Dd,\hat{c})=<\bar{L}>_N. 
\end{equation}
\end{cor}
Note that the equality $K_N(\TG)= <\bar{L}>_N$ does not imply that $K_N(L)$ is well defined by using state sums supported by {\it arbitrary} decorated triangulations $\TG$ of $(S^3,L)$, but only for those associated to link diagrams. However, by Proposition \ref{H=K} we know that $K_N(L)$ is fully well defined by $K_N(\TG)$.

Corollary \ref{Kteo1} is an unfolding in QH terms of the \cite[Theorem 1]{K2}. In that paper one considers decorated triangulations $\TG$
associated to link diagrams where a singular 3-ball dual to the branched spine shown in Figure \ref{Kconf} is associated to each crossing. Note that two walls intersect transversally at the middle; by sliding them we recover the configurations of Figure \ref{QH-Rmatrix+} and Figure \ref{QH-Rmatrix-}. The argument of Theorem 1 of \cite{K2} is purely local, based on a computation of the R-matrix $R_{K,N}$ by means of a combination of a version of the matrix dilogarithms (see the page 1417 and formula (4.25)).


\begin{figure}[ht]
\begin{center}
 \includegraphics[width=9cm]{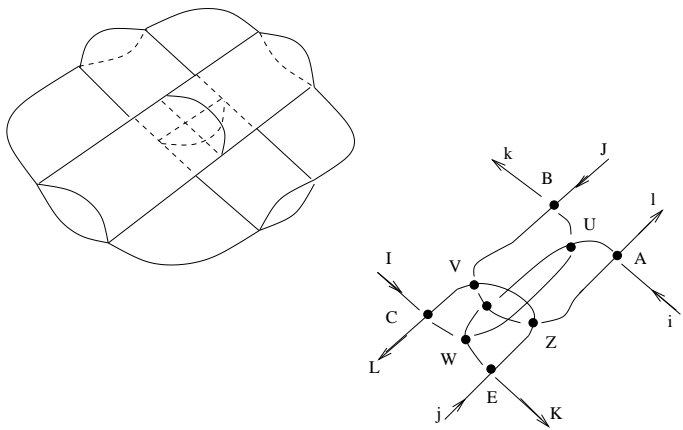}
\caption{\label{Kconf} Kashaev's configuration at crossings.}
\end{center}
\end{figure}

\section{Disproving the AbS conjecture}\label{AbS}

Recall that a link $L$ in $S^3$ is {\it hyperbolic} if $M=S^3\setminus L$
is a hyperbolic cusped 3-manifold, ie. if it admits a complete hyperbolic structure, which is necessarily of finite volume, and unique up to isometry by Mostow rigidity. 
\begin{conj} \label{KVCstatement} {\rm (Volume Conjecture, \cite{K4})} 
For every hyperbolic link $L$ in $S^3$ we have
$$2\pi\lim_{n\to +\infty} \log |<L>_n|/n = {\rm Vol}(M).$$
\end{conj}
Thanks to Theorem \ref{MM} the Kashaev Volume Conjecture can be
equivalently stated in terms of $J'_n(L)$, and in this form it will be
indicated as the {\it Kashaev-Murakami-Murakami Volume Conjecture}
(KMM VC) (see \cite{MM}). The KMM VC 
is known to hold true for a few knots (eg. the figure-eight knot $4_1$, or the knots $5_2$, $6_1$ and $6_2$), and also \cite{VdV} for the infinite family of {\it Whitehead Chain links}, including
the classical Whitehead's link $L_W$.
\smallskip

By Theorem \ref{equivalenceHK} the KMM VC can be recasted into the general framework of QH invariants. A formally similar problem concerns the semiclassical limit
of the QH invariants of hyperbolic cusped 3-manifold $M$. We have already proposed: 
\begin{conj}\label{CVC} {\rm (Cusped QH VC)} 
For every cusped hyperbolic $3$-manifold $M$, there is a weight
$\kappa$ such that $$2\pi \limsup_{N\to +\infty}  \{ \log|\Hh_N(M,\kappa)|/N \} = {\rm Vol}(M) \ . $$
\end{conj}
This conjecture has been checked when $M$ is the complement of the figure-eight knot $K$. In that case we have
\begin{equation}\label{firstsum}H_N(S^3\setminus K,\kappa)=N^2\frac{|g(w_0')|^2}{|g(1)|^2} \ \left| 1+ \sum_{\beta=1}^{N-1}
\zeta^{\beta^2} \prod_{k=1}^{\beta} \frac{w_1'^{-1}}{1-w_0'\zeta^k} \right|^2\end{equation}
where the function $g$ is defined in Section \ref{basic}, and the $N$th root cross-ratio moduli $w_0'$ and $w_1'$ have modulus $1$ and depend on the weight $\kappa$ and the complete hyperbolic structure of $S^3\setminus K$. Consider the ``diagonal'' sub state sum 
\begin{equation}\label{possum}
H_N^0(S^3\setminus K,\kappa)= N^2\frac{|g(w_0')|^2}{|g(1)|^2}
\ \left(1+\sum_{\beta=1}^{N-1} \prod_{k=1}^{\beta}
\frac{1}{\vert 1-w_0'\zeta^k\vert^2}\right)\ .
\end{equation}
By replacing formally $w_0'$ with $1$ we find the Kashaev's formula
\begin{equation}\label{Kformeight} <K>_N = 1+\sum_{\beta=1}^{N-1}
\prod_{i=1}^{\beta}|1-\zeta^i|^2 = N^{2}\left(1+\sum_{\beta=1}^{N-1}
\prod_{k=1}^{\beta} \frac{1}{\vert 1-\zeta^k\vert^2}\right)\ .\end{equation} 
\begin{remark}{\rm There are weights $\kappa$ such that \eqref{firstsum} and \eqref{possum} have the same semiclassical limit as \eqref{Kformeight}. 
Since the limsup in Conjecture \ref{CVC} vanishes for some weights (e.g. for the figure-eight knot and $\kappa=0$), these play a
subtle and rather mysterious r\^ole.}\end{remark}  
Because of coincidences like \eqref{possum}-\eqref{Kformeight}, it is rather natural to compare the asymptotic
behaviour of a general QH state sum $\Hh_N(\Tt)$ with its ``degeneration'' $\Hh_N(\Tt_\infty)$, defined as follows. For every QH tetrahedron $(\Delta,b,d)$ of
$\Tt$, consider the limit system of {\it signatures}
$$\sigma_j :=\lim_{N\to +\infty} w'_j = (-1)^{f_j-*_bc_j} .$$ If
$\sigma$ is {\it tame}, that is $\sigma_0=-1$
whenever $*_b=-1$, no singularities appear by replacing $w'_j$ with $\sigma_j$ in the matrix dilogarithm $\Rr_N(\Delta,b,d)$, so we get a {\it limit state sum} $\Hh_N(\Tt_\infty)$, $\Tt_\infty=(T,b,\sigma)$. When
$\sigma$ is not tame $\Hh_N(\Tt_\infty)$ can be defined anyway by
continuous extension. Then one can expect: 
\begin{equation}\label{AbSconj}\limsup \{ \log|\Hh_N(\Tt)|/N \} = \limsup \{
\log|\Hh_N(\Tt_\infty)|/N \}.\end{equation}
We call \eqref{AbSconj} the {\it Asymptotic by Signature (AbS) Conjecture}. We are going to disprove it.
\begin{lem}\label{signature} For every link $L$ there are QH triangulations
$\Tt_0$ and $\Tt_1$ supported by a same tringulation $(T,b)$ of $S^3$,
having a same tame signature $\sigma$ (hence the same
$\Tt_\infty$) , and such that for every odd $N$ we have
\begin{align*} \Hh_N(\Tt_0) =_N & \Hh_N(L\cup K_m) =_N   N\Hh_N(L)\\
\Hh_N(\Tt_1) =_N & \Hh_N(L+K_U) =  0 \end{align*}
where $L+K_U$ is the split link made by $L$ and the unknot $K_U$.
\end{lem}
A proof is illustrated by the puzzles of Figure
\ref{puzzle2}. The small picture at the bottom indicates that we start
with any o-graph associated to a $(1,1)$-tangle presentation of $L$, with a charge carrying $L\cup K_m$ and a tame signature. In the same spirit, the puzzles of Figure \ref{puzzle}
corresponding to the Whitehead and Hopf links are
supported by the same triangulation $(T,b)$ of $S^3$ and have the
same tame signature.

\begin{cor} The AbS conjecture is false.
\end{cor} 
\Dim Lemma \ref{signature} and \eqref{AbSconj} imply that $\limsup \{ \log|\Hh_N(L)|/N \}=0$ for every link $L$. This contradicts the KMM VC.\cvd

\begin{figure}[ht]
\begin{center}
 \includegraphics[width=8cm]{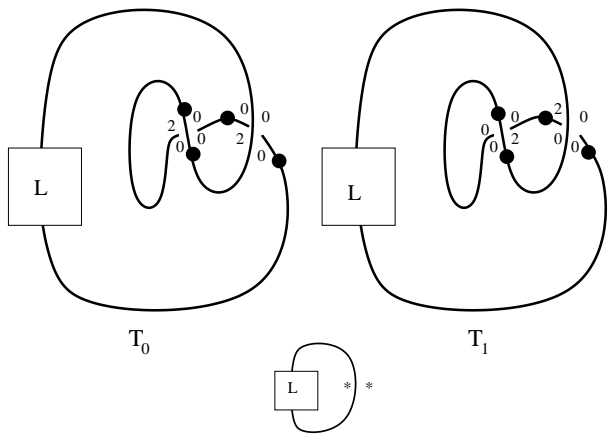}
\caption{\label{puzzle2} Sharing signature puzzles.}
\end{center}
\end{figure}

\section{Tensor computations}\label{compute}
As usual, let $N>1$ be any odd integer. The following functions are the basic ingredients of all our computations.
\smallskip

\subsection{Basic functions}\label{basic}
\begin{itemize}
\item For every $x\in \C^*$, we denote by $\log(x)$ the standard branch
of the Neperian logarithm, equal to the real $\log$ when $x>0$ and such that
$\log(-1)= \sqrt{-1}\pi$. The function $x^{1/N} := \exp(\log(x)/N)$
is extended to $0^{1/N}:=0$ by continuity.
\item For any $n\in \mz$ we denote by $[n]_N \in \Ii_N$ the residue of
$n$ {\rm mod}$(N)$.
\item For generic $u,v \in \C$ and any $n\in \mz$ we define $\omega(u\vert n)$ by the recurrence relation
$$ \omega(u \vert n+1) := \omega(u \vert n)\ (1-u\zeta^{n+1})\quad , \quad \omega(u \vert 0) :=1\ ,$$
and we set 
$$\omega(u,v \vert n) := \frac{v^{n}}{\omega(u\vert n)}\ .$$
In particular, $\omega(u,v\vert 0)=1$ and $\omega(u,v \vert n) = \prod_{j=1}^{n} \frac{v}{1-u\zeta^j}$ for any positive integer $n$.
\item We put $[x]:=N^{-1}\frac{1-x^N}{1-x}$, extended to $[1]:=1$ by continuity, and
$$g(x) :=  \prod_{j=1}^{N-1}(1 - x\zeta^{-j})^{j/N}\quad ,\quad h(x) :=  g(x)/g(1) \ .$$
\end{itemize}
\begin{remark}\label{first-rem}{\rm (1) Assume that $u^N+v^N=1$. Then $w(u,v\vert n)$ is $N$-periodic in the integer argument, so that $\omega(u,v \vert n)= \omega(u,v \vert [n]_N)$, and} $\omega(u,v \vert l)\ \omega(u\zeta^l,v\vert n) = \omega(u,v\vert l+n)$.

{\rm (2) We have the inversion relation $\omega(x\vert -n)\ \omega(x^{-1}\zeta^{-1} \vert n)= (-x)^{-n}\zeta^{\frac{n(n-1)}{2}}$.}
\end{remark}
Consider the rational function $f$ defined on the affine surface
$z^N=\frac{1-x^N}{1-y^N}$ by
\begin{equation}\label{deff}
 f(x,y| z)=\sum_{n=1}^{N} \prod_{j=1}^n \frac{1-y\zeta^j}{1-x\zeta^j}
\ z^n \ .
\end{equation}
The next lemma will be used frequently in the sequel. To simplify its notations, let us put $$(x)_n:= (1-x)(1-x\zeta)\dots (1-x\zeta^{[n]_N})=(1-x)\omega(x\vert  [n]_N)\ .$$ Denote by $^*$ the complex conjugation.
\begin{lem}\label{algres} We have:
\begin{enumerate}
\item [(i)] $x\ f(x,0| z\zeta) = (1-z) \ f(x,0|z)$
\item [(ii)] $g(x)\ g(z/\zeta)\ f(x,0|z) =_N x^{N-1}g(1)$
\item [(iii)] For all $m$, $n\in \mz$ it holds:
\begin{displaymath}
f(x\zeta^n,x\zeta^{-1}\vert \zeta^m) = \left\lbrace\begin{array}{ll}
x^{N-1-[m-1]_N}\ [x]^{-1} & \textrm{if}\ [n]_N =0 \\ 0 & \textrm{if}\
[n]_N \ne0,\ [m-1]_N < [n]_N\\ x^{N-1-[m-1]_N}\ [x]^{-1} & \\ \quad
\quad \ \frac{\zeta^{-nm}(\zeta)_{m-2}(x\zeta)_{n-1}}
{(\zeta^{n+1-m})^*_{m-n-2}(\zeta)_{n-1}^*}&
\textrm{if}\ [n]_N \ne0,\ [n]_N \leq [m-1]_N \end{array}\right.
\end{displaymath}
In particular, for $[n]_N \ne0$ and $[m]_N=0$ we get
\begin{equation}\label{eqnorm2}
f(x\zeta^n,x\zeta^{-1}\vert 1) = [x]^{-1}(x\zeta)_{n-1}.
\end{equation}

\item [(iv)] $g(x\zeta^n) = g(x) \ \omega(x,(1-x^N)^{1/N}\vert n)$
\end{enumerate}
\end{lem}
\noindent {\it Proof.} (i), (ii) and (iv) are proved in \cite{GT},
Lemma 8.2-8.3 (see also the Appendix of \cite{KMS}). We have also (see \cite{GT}, Proposition 8.6, page 569)
\begin{align}
f(x,x\zeta^{-1}\vert \zeta) = & x^{N-1}[x]^{-1}\label{eqnorm}\\
f(x,y \vert z\zeta) = & \frac{1-z}{x-yz\zeta}\ f(x,y \vert z)
\label{auto1}\\ f(x\zeta,y \vert z) = &
\frac{(1-x\zeta)(x-yz)}{z(x-y)}\ f(x,y \vert z).\label{auto2}
\end{align}
The case $[n]_N =0$ in (iii) follows directly from \eqref{eqnorm}-\eqref{auto1}. When $[m-1]_N < [n]_N$ we get
\begin{equation}\label{auto3}
f(x\zeta^n,x\zeta^{-1} \vert \zeta^m) = x^{-[m-1]_N} \
\prod_{k=1}^{[m-1]_N} \frac{1-\zeta^{m-k}}{\zeta^n-\zeta^{m-k}}\
f(x\zeta^n,x\zeta^{-1} \vert \zeta).
\end{equation}
By using \eqref{auto2} we see that $f(x\zeta^n,x\zeta^{-1} \vert \zeta^m)=0$ if
moreover $[n]_N\ne 0$. Finally, when $[n]_N \ne0$ and $[n]_N \leq [m-1]_N$
there is a simple pole in the product, and \begin{multline*}
f(x\zeta^n,x\zeta^{-1} \vert \zeta^m) = x^{-[m-1]_N} \ (\zeta)_{m-2} \
(x\zeta)_{n-1}\ f(x,x\zeta^{-1}\vert \zeta)\ \frac{(\zeta^k-1)_{\vert
k=0}}{(\zeta^n-\zeta^{m-k})_{\vert [m-k]_N=[n]_N}}\\ \prod_{k}
\frac{1}{\zeta^n-\zeta^{m-k}}
\end{multline*}
where $k$ goes from $1$ to $[m]_N$ and
$[m-k]_N\ne[n]_N$ in the product, which is easily seen to be equal to
$\zeta^{-n(m-1)}/((\zeta^{n+1-m})^*_{m-n-2}(\zeta)^*_
{n-1})$. The last case of (iii) then follows from \eqref{eqnorm} and
\begin{displaymath}
 \frac{(\zeta^k-1)_{k=0}}{(\zeta^n-\zeta^{m-k})_{[m-k]_N=[n]_N}} =
 \zeta^{-n}.
\end{displaymath}
The reduced formula for $m=0$ is a consequence of $(\zeta)_{N-2} =
(\zeta^{1-m})^*_{N-n-2}(\zeta)^*_{n-1} = N$. This
concludes the proof.\cvd
\begin{remark}
{\rm The case $[n]_N \ne0$, $[n]_N \leq [m-1]_N$ is not made explicit in the
proof of Proposition 8.6 of \cite{GT}. In fact, at page 569,
line -7 and -2, of that paper, there are two possible cases
corresponding to $[n]_N \leq [m-1]_N$ and $[n]_N > [m-1]_N$, but the
Kronecker symbol $\delta_N(i+j-k-l)$ in line 3 selects the second
one.}
\end{remark}

{\bf Convention.} {\it From now on we set $N=2m+1$, $m\geq 1$, so that
``$m$'' is a reserved character}.

\subsection{The matrix dilogarithms}\label{dilogs}
The $N$-{\it matrix dilogarithm} of a branched tetrahedron $(\Delta,b)$ with
QH decoration $d=(w,f,c)$ (see Section \ref{o-gr}) is given by

\begin{equation}\label{symqmatdil}
\mathcal{R}_N(\Delta,b,d) = \mathcal{R}_N(*_b,d) =
\bigl((w_0')^{-c_1}(w_1')^{c_0}\bigr)^{\frac{N-1}{2}}\
(\Ll_N)^{*_b}(w_0',(w_1')^{-1})\quad \in {\rm Aut}(\mc^N\otimes \mc^N)
\end{equation}
where 
\begin{align*} \Ll_N(u,v)_{k,l}^{i,j}& =  h(u)\
\zeta^{kj+(m+1)k^2}\ \omega(u,v\vert i-k) \ \delta_N(i + j - l) \\
\bigl( \Ll_N(u,v)^{-1}\bigr)^{k,l}_{i,j}& =  \frac{[u]}{h(u)}\
\zeta^{-kj-(m+1)k^2}\ \frac{\delta_N(i+j-l)}{\omega(u/\zeta,v\vert
i-k)}.
\end{align*} 
Note that Remark \ref{first-rem} (1) applies in this case.

\subsection{Discrete Fourier transform}\label{discreteFT} 
We call {\it (discrete) Fourier transformation} the conjugation by
tensor powers of the automorphism $F$ of $\C^N$ with entries
$F^i_j=\zeta^{ij}/\sqrt{N}$. Hence, the Fourier transform of the $N$-matrix
dilogarithm $\Rr_N(*_b,d)$ is
$$\tilde{\Rr}(*_b,d) = F^{\otimes 2}\circ \Rr_N(*_b,d)\circ
(F^{-1})^{\otimes 2}. $$
In general, for every QH triangulated polyhedron $\Yy$ (possibly with
free 2-faces, as in Section \ref{ChargetoOp}) we denote by
$\Hh_N(\Yy)$ the QH tensor obtained by using the
original matrix dilogarithms, and by $\tilde{\Hh}_N(\Yy)$ its Fourier transform. Clearly, for
every QH triangulation $\Tt$ of a {\it closed} pseudomanifold we have
$\Hh_N(\Tt)=\tilde{\Hh}_N(\Tt)$.
\begin{lem}\label{FFourier} We have
\begin{align*}\tilde{\Rr}_N(+,d)^{i,j}_{k,l}= &
\frac{\bigl((w_0')^{-c_1+2}(w_1')^{c_0}\bigr)^{\frac{N-1}{2}}}
{Ng((w_1')^{-1}/\zeta)}\
\frac{\zeta^{(k-i)(j-l)+(m+1)(j^2-l^2)}}{\omega((w_1')^{-1}/\zeta,w_0'\vert
l-i)}\\
\tilde{\Rr}_N(-,d)_{i,j}^{k,l}= &
g((w_1')^{-1}/\zeta)\ [(w_1')^{-1}]\
\bigl((w_0')^{-c_1-2}(w_1')^{c_0}\bigr)^{\frac{N-1}{2}}\
\\ & \hspace*{5cm} \zeta^{(i-k)(j-l)+(m+1)(l^2-j^2)}\ \omega((w_1')^{-1},w_0'\vert
l-i) \ . \end{align*}
\end{lem}

\noindent {\it Proof.} By direct substitution we find
\begin{align*} & \tilde{\Ll}_N(w_0',(w_1')^{-1})^{i,j}_{k,l}  
= N^{-2}\sum_{\alpha,\beta,\gamma,\delta=0}^{N-1}
\zeta^{k\alpha+l\beta-i\gamma-j\delta}
\Ll_N(w_0',(w_1')^{-1})_{\alpha,\beta}^{\gamma,\delta}
\\&
=N^{-2}h(w_0')\sum_{\alpha,\gamma=0}^{N-1}\zeta^{(k-\gamma+(m+1)
\alpha)\alpha+(j-i)\gamma}\ \omega(w_0',(w_1')^{-1}\vert
\gamma-\alpha) \sum_{\beta=0}^{N-1}\zeta^{\beta(\alpha+l-j)}\\
&
 =N^{-1}h(w_0')\sum_{\alpha,\gamma=0}^{N-1}\zeta^{(k-\gamma+(m+1)\alpha)
 \alpha+(j-i)\gamma}\ \omega(w_0',(w_1')^{-1}\vert \gamma-\alpha) \
 \delta_N(\alpha+l-j)\\&
 =N^{-1}h(w_0')\zeta^{k(j-l)+(m+1)(j-l)^2}\sum_{\gamma=0}^{N-1}
 \zeta^{\gamma(l-i)}\ \omega(w_0',(w_1')^{-1}\vert \gamma+l-j)\\&
 =N^{-1}h(w_0')\zeta^{(k+l-i)(j-l)+(m+1)(j-l)^2}\sum_{\gamma=0}^{N-1}
 \zeta^{(\gamma+l-j)(l-i)}\ \omega(w_0',(w_1')^{-1}\vert\gamma+l-j)\\&
 =_N N^{-1}\left(g((w_1')^{-1}/\zeta)\right)^{-1}(w_0')^{N-1}
 \frac{\zeta^{(k-i)(j-l)+(m+1)(j^2-l^2)}}{\omega((w_1')^{-1}/\zeta,w_0'\vert
 l-i)}.\end{align*} In the last equality we use Lemma
 \ref{algres}(i)-(ii). For negative branching orientation, since the
 discrete Fourier transform of the inverse is the inverse of the
 discrete Fourier transform it is enough to check the formula
 \begin{displaymath}\tilde{\Ll}_N^{-1}(w_0',(w_1')^{-1})^{i,j}_{k,l} =
 g((w_1')^{-1}/\zeta)\frac{[(w_1')^{-1}]}{(w_0')^{N-1}}\
 \zeta^{(k-i)(l-j)+(m+1)(j^2-l^2)}\ \omega((w_1')^{-1},w_0'\vert
 j-k).\end{displaymath} 
 
We have
\begin{align*}& \left(\tilde{\Ll}_N(w_0',(w_1')^{-1})\circ 
\tilde{\Ll}_N^{-1}(w_0',(w_1')^{-1})\right)^{i,j}_{m,n} = \\& =N^{-1}\
[(w_1')^{-1}]\
\zeta^{(m+1)(j^2-n^2)}\sum_{k,l=0}^{N-1}\zeta^{(m-k)(l-n)-(k-i)(j-l)}
\frac{\omega((w_1')^{-1},w_0'\vert
j-k)}{\omega((w_1')^{-1}/\zeta,w_0'\vert n-k)}\\& = \delta_N(m-i)\
[(w_1')^{-1}]\
\zeta^{(m+1)(j^2-n^2)+i(j-n)}\sum_{k=0}^{N-1}\zeta^{k(n-j)}
\frac{\omega((w_1')^{-1},w_0'\vert
j-k)}{\omega((w_1')^{-1}/\zeta,w_0'\vert n-k)}\\&= \delta_N(m-i)
[(w_1')^{-1}] \zeta^{(m+1)(j^2-n^2)+(i-n)(j-n)}\\ &
\hspace{2cm}\sum_{k=0}^{N-1}\zeta^{(n-k)(j-n)} \
\frac{\omega((w_1')^{-1},w_0'\vert
j-n)\omega((w_1')^{-1}\zeta^{j-n},w_0'\vert
n-k)}{\omega((w_1')^{-1}/\zeta,w_0'\vert n-k)}\end{align*}

\begin{align*}& = \delta_N(m-i)\ [(w_1')^{-1}]\ 
\zeta^{(m+1)(j^2-n^2)+(i-n)(j-n)} \ \omega((w_1')^{-1},w_0'\vert j-n)\
\frac{\delta_N(j-n)}{[(w_1')^{-1}]}\\& = \delta_N(m-i)\
\delta_N(j-n).\end{align*} The sum in the fourth equality is computed
by using Lemma \ref{algres} (iii).\cvd

\subsection{Braidings} Figure \ref{braidcomplex} 
(bottom) shows a {\it tunnel crossing}, that is, the portion of
branched spine corresponding to the portion of o-graph on the left of
Figure \ref{SQUARE}. Note that for the moment no wall has been
inserted within the tunnels. The dual {\it singular} octahedron {\bf
O}, which has two pairs of identified edges, is shown on the top right. The indices $i1,i2, \dots \ \in \Ii_N$ refer to state variables.
\begin{figure}[ht]
\begin{center}
\includegraphics[width=13cm]{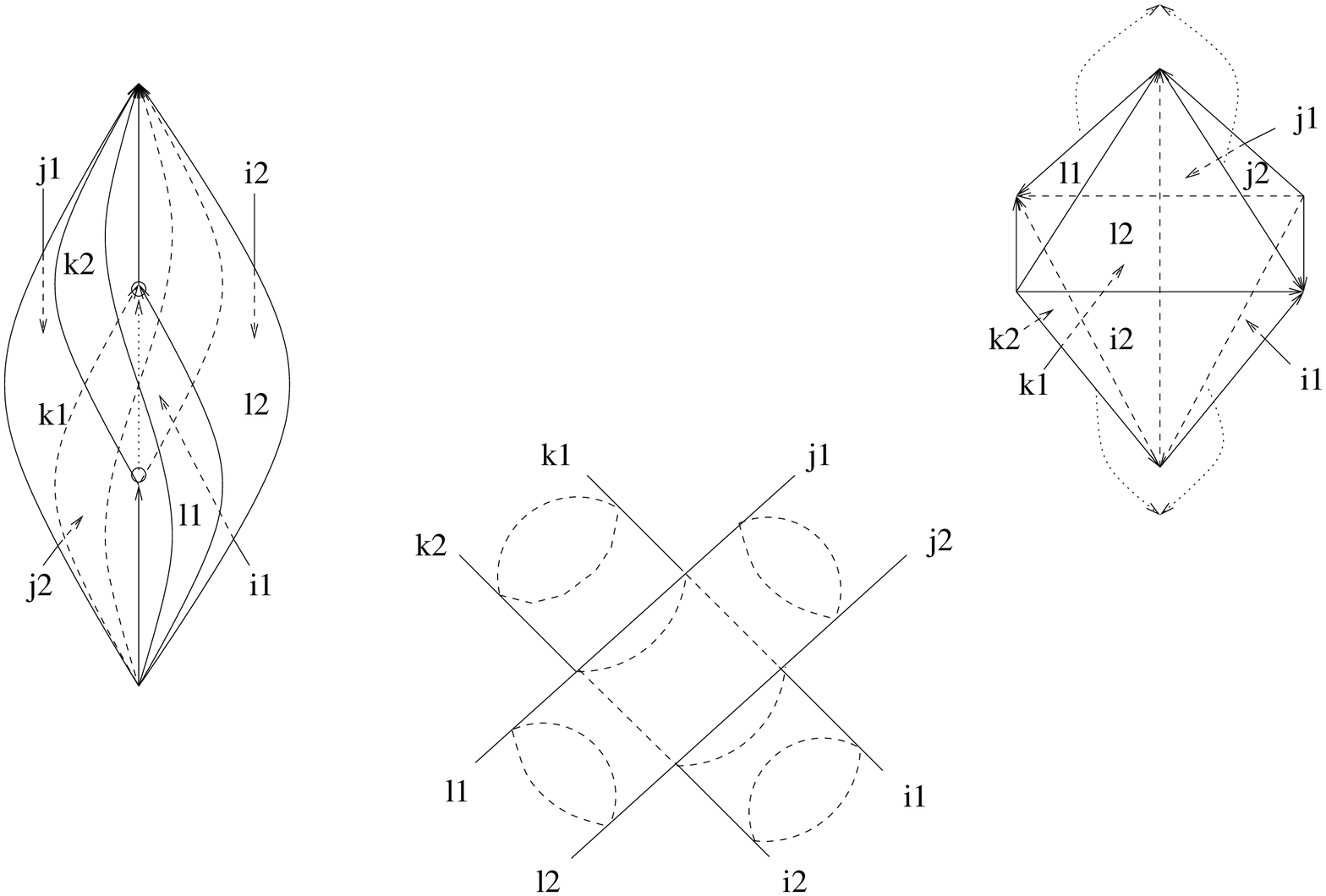}
\caption{\label{braidcomplex} A tunnel crossing and the dual singular
octahedron ${\bf O}$.}
\end{center}
\end{figure} 
\smallskip  

Denote by $\mathcal{O}$ any {\it distinguished} QH polyhedron
supported by {\bf O}, with tetrahedra
$\underline{\Delta}^i=(\Delta^i,b^i,d^i)$, $i=1,\ldots,4$ ordered in
the counterclockwise way about the central axis, starting from the
front $3$-simplex. So $\underline{\Delta}^1$ contains the edge dual to
the planar regions at l2 and i2, $\underline{\Delta}^2$ corresponds to
the regions at i1 and j2, and so on. $\underline{\Delta}^1$ and
$\underline{\Delta}^3$ (resp. $\underline{\Delta}^2$ and
$\underline{\Delta}^4$) have negative (resp. positive) branching
orientation. Here ``distinguished QH polyhedron'' means that we are
using the usual universal constant $(w,f)$, and that the charges $c_i^j$ at the internal
edge satisfy
\begin{equation}
\label{eqfc1}
c_1^1+c_1^2+c_1^3+c_1^4=2.
\end{equation}
At the top left of Figure \ref{braidcomplex} we consider $\Oo$ as a singular QH cobordism
between twice punctured $2$-disks with identified punctures. By looking {\it from right to left} at the bottom picture we consider $\Oo$ as associated to a {\it negative}
crossing. With the notations of Section
\ref{ChargetoOp}, it corresponds to the braiding
tensor (recall that it is converted and based on the discrete Fourier transform):
$${\rm Br}_N(-,c)^{j1,j2,i1,i2}_{k1,k2,l1,l2}$$
that belongs to ${\rm End}((\C^N)^{\otimes 4})$. Here the charge $c$
is indicated as a varying parameter. Similarly, by looking {\it from bottom to top} we consider $\Oo$ as a {\it positive}
crossing, and it corresponds to the braiding tensor
$${\rm Br}_N(+,c)^{i1,i2,l2,l1}_{j2,j1,k1,k2}\ . $$
It can be checked (see the proof of Lemma \ref{braiding1}
below) that ${\rm Br}_N(-,c)$ is equal to
\begin{equation}\label{operatorform}
(\tilde{\Rr}_N(\underline{\Delta}^4)_{23}^{t_2}\circ P_{23})\circ(\tilde{\Rr}_N(\underline{\Delta}^1)_{34}
^{t_3t_4}\circ P_{34})\circ
(\tilde{\Rr}_N(\underline{\Delta}^3)_{12} \circ P_{12})\circ(\tilde{\Rr}_N(\underline{\Delta}^2)_{23}^
{t_3} \circ P_{23}),\end{equation} where eg. $\Rr_N(\underline{\Delta}^3)_{12}$
means $\Rr_N(\underline{\Delta}^3)$ acting on the first and second
tensor factor of $(\mc^N)^{\otimes 4}$, $t_i$ is the transposition
on the $i$th factor, and $P_{ij}$ the flip map. In particular, by invertibility of the matrix
dilogarithm and their partial transpose we see that
the braiding tensors ${\rm Br}_N(-,c)$ are automorphisms of $(\C^N)^{\otimes 4}$.

Put
$$\begin{array}{l} K_\mathcal{O} = N\ [((w_1^1)')^{-1}][((w_1^3)')^{-1}]\
\frac{g(((w_1^1)')^{-1}/\zeta)g(((w_1^3)')^{-1}/\zeta)}{g(((w_1^2)')^{-1}/
\zeta)g(((w_1^4)')^{-1}/\zeta)}\times \\
\bigl(((w_0^1)')^{-c_1^1-2}((w_1^1)')^
{c_0^1}((w_0^3)')^{-c_1^3-2}((w_1^3)')^{c_0^3}((w_0^2)')^
{-c_1^2+2}((w_1^2)')^{c_0^2}((w_0^4)')^{-c_1^4+2}((w_1^4)')^
{c_0^4}\bigr)^{\frac{N-1}{2}}\end{array}$$
and $\textstyle \bar{K}_\mathcal{O} = N\ K_\mathcal{O}^{-1}\prod_{i=1}^3
[((w_1^i)')^{-1}]((w_0^2)')^{1-N}$.
\begin{lem}\label{braiding1} {\rm (Braiding tensor)} We have
\begin{multline*}{\rm Br}_N(-,c)^{j1,j2,i1,i2}_{k1,k2,l1,l2}=
{\rm Br}_N(+,c)^{i1,i2,l2,l1}_{j2,j1,k1,k2}= \\  
   K_\mathcal{O}\ \delta_N(i_{12}-k_{12})\ \delta_N(l_{12}-j_{12})\
   \zeta^{(l_1-j_1)(k_1-i_1-l_{12})} \\
\times \ \frac{\omega(((w_1^1)')^{-1},(w_0^1)'\vert
l_2-i_2) \omega(((w_1^3)')^{-1},(w_0^3)'\vert
j_1-k_1)}{\omega(((w_1^2)')^{-1}/\zeta,(w_0^2)'\vert
j_2-i_1)\omega(((w_1^4)')^{-1}/\zeta,(w_0^4)'\vert l_1-k_2)}\end{multline*}
\begin{multline*}{\rm Br}^{-1}_N(-,c)^{l2,l1,k2,k1}_{i2,i1,j2,j1}=
\bar{K}_\mathcal{O}\ \delta_N(i_{12}-k_{12})\
\delta_N(l_{12}-j_{12})\
\zeta^{(l_1-j_1)(i_1-k_1+l_{12})}\\ \times\
\frac{\omega(((w_1^2)')^{-1},(w_0^2)'\vert
j_2-i_1-1)\omega(((w_1^4)')^{-1}/\zeta,(w_0^4)'\vert
l_1-k_2)}{\omega(((w_1^1)')^{-1}/\zeta,(w_0^1)'\vert l_2-i_2
)\omega(((w_1^3)')^{-1}/\zeta,(w_0^3)'\vert j_1-k_1)}
\end{multline*}
where $i_{12}=i_1-i_2$, and similarly for $j_{12}$, $k_{12}$ and $l_{12}$.
\end{lem}
\noindent {\it Proof.} The equality 
$${\rm Br}_N(-,c)^{j1,j2,i1,i2}_{k1,k2,l1,l2}= {\rm
Br}_N(+,c)^{i1,i2,l2,l1}_{j2,j1,k1,k2}$$
depends on the fact that both tensors are formal conversions (see Section \ref{ChargetoOp}) of a same QH tensor
$\tilde{\Hh}_N(\mathcal{O})$.  We compute this last:
$$\begin{array}{l} \sum_{\alpha,\beta,\gamma,\delta=1}^N
\tilde{\Rr}_N(\underline{\Delta}^1)^{\delta,l_2}_{i_2,\alpha}\
\tilde{\Rr}_N(\underline{\Delta}^2)^{i_1,\alpha}_{\beta,j_2}\
\tilde{\Rr}_N(\underline{\Delta}^3)^{\beta,j_1}_{k_1,\gamma}\
\tilde{\Rr}_N(\underline{\Delta}^4)^{k_2,\gamma}_{\delta,l_1}\\
\hspace{2.5cm}= N^{-3}K_\mathcal{O}\
\frac{\omega(((w_1^1)')^{-1},(w_0^1)'\vert
l_2-i_2)\omega(((w_1^3)')^{-1},(w_0^3)'\vert
j_1-k_1)}{\omega(((w_1^2)')^{-1}/\zeta,(w_0^2)'\vert
j_2-i_1)\omega(((w_1^4)')^{-1}/\zeta,(w_0^4)'\vert l_1-k_2)}\ \times
\\
\hspace{1.5cm}\zeta^{(m+1)(j_1^2+l_2^2-j_2^2-l_1^2)}
\sum_{\alpha,\beta,\gamma,\delta=1}^N\zeta^{(\beta-i_1)(\alpha-j_2)+
(\delta-k_2)(\gamma-l_1)+(k_1-\beta)(\gamma-j_1)+(i_2-\delta)(\alpha-l_2)}.
\end{array}$$ 
The last sum equals 
$$\begin{array}{l}
\sum_{\alpha,\beta,\gamma=1}^N\zeta^{(\beta-i_1)(\alpha-j_2)-k_2(\gamma-l_1)
+(k_1-\beta)(\gamma-j_1)+i_2(\alpha-l_2)}\sum_{\delta=1}^N
\zeta^{\delta(\gamma-l_1-\alpha+l_2)}\hspace*{3cm}\\
 \hspace*{2cm}= N\
 \zeta^{i_1j_2+k_2l_2+k_1(l_1-l_2-j_1)-i_2l_2}
\sum_{\alpha,\beta=1}^N\zeta^{\alpha(i_2-i_1+k_1-k_2)+\beta(j_1-j_2+l_2-l_1)}\\
 \hspace*{2cm}= N^3\ \zeta^{(k_1-i_1)(l_1-j_1)}\ 
\delta_N(i_{12}-k_{12})\ \delta_N(l_{12}-j_{12}).
\end{array}$$
This is non vanishing if and only if $j_2=j_1+l_2-l_1$, which implies
$\zeta^{(m+1)(j_1^2+l_2^2-j_2^2-l_1^2)}=
\zeta^{j_1l_1+l_1l_2-l_1^2-j_1l_2}=\zeta^{-(l_1-j_1)l_{12}}$. This proves
the first two equalities.

To compute the inverse, recall (\ref{operatorform}) and that overall
transposition commutes with taking the inverse. By definition,
$$\tilde{\Rr}_N^{-1}(-,d)^{i,j}_{k,l} = N^{-1}\left(g((w_1')^{-1}/\zeta)\right)^{-1}\bigl((w_0')^{c_1+2}(w_1')^{-c_0}
\bigr)^{\frac{N-1}{2}}\
\frac{\zeta^{(k-i)(j-l)+(m+1)(j^2-l^2)}}{\omega((w_1')^{-1}/\zeta,w_0'\vert
l-i)}$$
If $*_b=1$ we find as in Lemma \ref{FFourier} that
$$ \begin{array}{lll}
\left((\tilde{\Rr}_N(+,d)_{12}^{t_2})^{-1}\right)^{i,j}_{k,l}
&=&
\bigl((w_0')^{c_1-4}(w_1')^{-c_0}\bigr)^{\frac{N-1}{2}}g((w_1')^{-1}/
\zeta)\ [(w_1')^{-1}]\hspace{4cm}\\ & & \hspace{1cm}
\zeta^{(k-i)(j-l)+(m+1)(l^2-j^2)}\omega((w_1')^{-1},w_0'\vert l-k-1)\\
\\
\left((\tilde{\Rr}_N(+,d)_{12}^{t_1})^{-1}\right)^{i,j}_{k,l}&=
&N^{-1}\bigl((w_0')^{c_1-2}(w_1')^{-c_0}\bigr)^{\frac{N-1}{2}}g((w_1')^{-1}/
\zeta)\hspace{4cm}\\ & & \hspace{1cm}
\zeta^{(k-i)(j-l)+(m+1)(j^2-l^2)}\omega((w_1')^{-1}/\zeta,w_0'\vert
j-i).\end{array}$$ Hence the entries of ${\rm Br}^{-1}_N(-,c)$ are computed by
$$\begin{array}{l} 
\sum_{\alpha,\beta,\gamma,\delta=1}^N
\left((\tilde{\Rr}_N(\underline{\Delta}^4)_{23}^{t_2})^{-1}\right)_
{\delta,\gamma}^{k_2,l_1}
\
\left(\tilde{\Rr}_N^{-1}(\underline{\Delta}^1)_{34}^{t_3t_4}\right)^
{\delta,l_2}_{i_2,\alpha}\hspace{3cm}\\
\hspace*{8cm}
\left(\tilde{\Rr}_N^{-1}(\underline{\Delta}^3)_{12}\right)_{\beta,j_1}^
{k_1,\gamma}\left((\tilde{\Rr}_N(\underline{\Delta}^2)_{23}^{t_3})^{-1}
\right)_{i_1,j_2}^{\beta,\alpha}\\ = N^{-3}\bar{K}_\mathcal{O}\
\frac{\omega(((w_1^2)')^{-1},(w_0^2)'\vert
j_2-i_1-1)\omega(((w_1^4)')^{-1}/\zeta,(w_0^4)'\vert
l_1-k_2)}{\omega(((w_1^1)')^{-1}/\zeta,(w_0^1)'\vert
l_2-i_2)\omega(((w_1^3)')^{-1}/\zeta,(w_0^3)'\vert j_1-k_1)}\ \times\\
\hspace{1cm}\zeta^{-(m+1)(j_1^2+l_2^2-j_2^2-l_1^2)}
\sum_{\alpha,\beta,\gamma,\delta=1}^N\zeta^{(i_1-\beta)(\alpha-j_2)+
(k_2-\delta)(\gamma-l_1)+(\beta-k_1)(\gamma-j_1)+(\delta-i_2)(\alpha-l_2)}. 
\end{array}$$
At this point we can conclude as above, by computing the
exponents of $\zeta$.\cvd

\subsection{Walls}\label{WALLS}
Our next task is to compute the QH tensors of walls, which are encoded
by the o-graph portion on the right of Figure \ref{SQUARE}.  As usual
we adopt the universal constant system $(w_0,f_0,f_1) = (2,0,-1)$. The
two tetrahedra $\Delta^\pm$, with branching signs $*_b = \pm 1$,
occurring in any wall have decorations $d^\pm = (w,f,c^\pm)$ differing
only for the charges, that is $c^+=(P,F,H)$ and $c^-=(M,G,K)$. As $H$
and $K$ are immaterial in matrix dilogarithm formulas, a generic wall will be
denoted by $\Ww=\Ww(P,F,M,G)$, so that $\Ww_C=\Ww(0,0,1,0)$,
$\Ww_M=\Ww(0,-1,1,1)$ and the wall type introduced in Lemma
\ref{convert} is $\Ww(0,2,-1,0)$. Adopting the
notations of Section \ref{ChargetoOp}, we have to compute the QH
tensors
\begin{align*}\Hh_N(\Ww)^{j,k}_{i,l}= &
N^{-1}\sum_{\alpha,\beta=0}^{N-1}\Rr_N(+,d^+)
^{j,\alpha}_{\beta,l}\ \Rr_N(-,d^-)^{\beta,k}_{i,\alpha}\\
\tilde{\Hh}_N(\Ww)^{j,k}_{i,l} = &
N^{-1}\sum_{\alpha,\beta=0}^{N-1}\tilde{\Rr}_N(+,d^+)
^{j,\alpha}_{\beta,l}\ \tilde{\Rr}_N(-,d^-)^{\beta,k}_{i,\alpha}\end{align*}\ .
Denote by $w_j'$ and $z_j'$ the $N$th root moduli of $\Delta^-$ and
$\Delta^+$. 
We have:
\begin{itemize}
\item For $\Ww_C$: $w_0'=\sqrt[N]{2}\zeta^{(m+1)}$, $w_1'=-1$,
 $z_0' = w_0'\zeta^{-(m+1)}$ and $z_1' = w_1'$;
\item For $\Ww_M$: $w_0'=\sqrt[N]{2}\zeta^{(m+1)}$,
$w_1'=\exp(\pi\sqrt{-1}/N)$, $z_0' = w_0'\zeta^{-(m+1)}$ and $z_1' =
w_1'$.
\end{itemize}
\begin{lem}\label{defwalls} {\rm (Wall QH tensors)} Let $\Ww=\Ww_C$ or $\Ww= \Ww_M$. Then: 
\begin{align*}
\tilde{\Hh}_N(\Ww)^{j,k}_{i,l} & =_N N^{-1}
\delta_N(i-j)\delta_N(k-l)\zeta^{(m+1)(k-i)}(w_1')
^{\frac{N-1}{2}}\frac{1-w_1^{-1}}{1-(w_1')^{-1}\zeta^{k-i}}\\
\Hh_N(\Ww)^{j,k}_{i,l} & =_N N^{-1}
(w_1')^{\frac{N-1}{2}-[j-i-(m+1)]_N}\delta_N(l-k-(j-i))\ .
\end{align*}
 \end{lem}
\noindent {\it Proof.} In both cases we have the same relations
between variables $w_j'$ and $z_j'$, so the respective QH tensors 
have the same form (also by taking care of the scalar factors, which depend on the charges). From Lemma \ref{FFourier} we obtain easily:
$$\begin{array}{l}\tilde{\Hh}_N(\Ww)^{j,k}_{i,l} =
N^{-1}\sum_{\alpha,\beta=0}^{N-1}\tilde{\Rr}_N(+,d^+)^{j,\alpha}_{\beta,l}\
\tilde{\Rr}_N(-,d^-)^{\beta,k}_{i,\alpha} \hspace{6cm}\\ \hspace{1cm}
= N^{-2}[(w_1')^{-1}](w_1')^{\frac{N-1}{2}}\zeta^{(m+1)}\
\frac{\omega((w_1')^{-1},w_0'\vert
k-i)}{\omega((w_1')^{-1}/\zeta,w_0'\zeta^{-(m+1)}\vert l-j)} \\
\hspace{4cm} \times\ \sum_{\alpha,\beta=0}^{N-1}\
\zeta^{(\beta-j)(\alpha-l) +
(m+1)(\alpha^2-l^2)+(i-\beta)(\alpha-k)+(m+1)(k^2-\alpha^2)}\\\hspace{1cm}
= N^{-2}[(w_1')^{-1}](w_1')^{\frac{N-1}{2}}\zeta^{(m+1)}\
\frac{\omega((w_1')^{-1},w_0'\vert
k-i)}{\omega((w_1')^{-1}/\zeta,w_0'\zeta^{-(m+1)}\vert l-j)}\\
\hspace{7.5cm} \times \ \zeta^{jl-ik-(m+1)l^2+(m+1)k^2}\
\sum_{\alpha,\beta=0}^{N-1}\ \zeta^{\beta(k-l)+\alpha(i-j)}\\
\hspace{1cm} = [(w_1')^{-1}] (w_1')^{\frac{N-1}{2}}\zeta^{(m+1)(1+k-j)}\delta_N(i-j)
\delta_N(k-l)\frac{1-(w_1')^{-1}}{1-(w_1')^{-1}\zeta^{k-i}}.\end{array}$$
We will now compute directly the QH tensors rather than apply the
discrete Fourier transform to the result we have just obtained. From
the formulas of the $N$-matrix dilogarithms we get
$$\begin{array}{lll}\Hh_N(\Ww)^{j,k}_{i,l} & =&
N^{-1}\sum_{\alpha,\beta=0}^{N-1}\Rr_N(+,d^+)^{j,\alpha}_{\beta,l}\
\Rr_N(-,d^-)^{\beta,k}_{i,\alpha} \hspace{6cm}\\ & = &
N^{-1}(w_1')^{\frac{N-1}{2}}\ h(w_0'\zeta^{-(m+1)})\ h(w_0')^{-1}\
[w_0']\ \times \\ & & \hspace{1cm}
\sum_{\alpha,\beta=0}^{N-1}\frac{\omega(w_0'\zeta^{-(m+1)},(w_1')^{-1}\vert
j-\beta)}{\omega(w_0'/\zeta,(w_1')^{-1}\vert i-\beta)}\
\delta_N(j+\alpha-l)\delta_N(i+\alpha-k)\\ & = &
N^{-1}(w_1')^{\frac{N-1}{2}}\ g(w_0'\zeta^{-(m+1)})\ g(w_0')^{-1}\
[w_0']\ \delta_N(l-k-(j-i))\ \times \\ & &
\hspace{1cm}\sum_{\beta=0}^{N-1}\frac{\omega(w_0'\zeta^{-(m+1)},(w_1')^{-1}
\vert j-\beta)}{\omega(w_0'/\zeta,(w_1')^{-1}\vert
i-\beta)}.\end{array}$$ Factorizing as in Remark \ref{first-rem} (1),
\begin{multline}\label{factorization}\omega(w_0'\zeta^{-(m+1)},
(w_1')^{-1}\vert j-\beta) = \omega(w_0'\zeta^{-(m+1)},(w_1')^{-1}\vert
j-i)\hspace{2cm}\\ \omega(w_0'\zeta^{j-i-(m+1)},(w_1')^{-1}\vert
i-\beta)\end{multline} and using Lemma \ref{algres} (iv) to compute
the ratio of $g$ functions we find
\begin{multline*} 
\Hh_N(\Ww)^{j,k}_{i,l}=_N N^{-1}(w_1')^{\frac{N-1}{2}}\ [w_0']\
\delta_N(l-k-(j-i))\ \times \hspace{6.6cm}\\ \hspace{3.5cm}
\frac{\omega(w_0'\zeta^{-(m+1)},(w_1')^{-1}\vert
j-i)}{\omega(w_0'\zeta^{-(m+1)},(w_1')^{-1}\vert m+1)}\
f(w_0'\zeta^{j-i-(m+1)},w_0'\zeta^{-1} \vert 1)\end{multline*} where
the function $f$ is defined in (\ref{deff}). From Lemma \ref{algres}
(iii) we get
\begin{equation}\label{1red}
f(w_0'\zeta^{j-i-(m+1)},w_0'\zeta^{-1} \vert 1) =
\left\lbrace\begin{array}{ll} [w_0']^{-1} & \textrm{if}\ [j-i-(m+1)]_N
=0 \\ (w_0'\zeta)_{j-i-(m+1)-1}\ [w_0']^{-1} & \textrm{if}\
[j-i-(m+1)]_N \ne 0 \ .\end{array}\right.
\end{equation}
The formula for $\Hh_N(\Ww)^{j,k}_{i,l}$ follows immediately from this
and \eqref{factorization}. We can check it is coherent with the
formula of $\tilde{\Hh}_N(\Ww)^{j,k}_{i,l}$ we had previously
obtained, as follows:
\begin{align*}
\tilde{\Hh}_N(\Ww)^{j,k}_{i,l} & =_N N^{-2} \sum_{I,J,K,L=0}^{N-1}
\zeta^{Ii+Ll-Jj-Kk}\ \Hh_N(\Ww)^{J,K}_{I,L}\\ & =_N N^{-3}
(w_1')^{\frac{N-1}{2}} \sum_{I,J,K,L=0}^{N-1} \zeta^{Ii+Ll-Jj-Kk}\
(w_1')^{-[J-I-(m+1)]_N}\delta_N(L-K-(J-I))\\ & =_N
N^{-3}(w_1')^{\frac{N-1}{2}}\zeta^{(m+1)(l-j)}\sum_{K,I=0}^{N-1}
\zeta^{I(i-j)+K(l-k)}\sum_{[J-I-(m+1)]_N=0}^{N-1}(w_1'\zeta^{j-l})
^{-[J-I-(m+1)]_N}\\ & =_N N^{-1}
\delta_N(i-j)\delta_N(k-l)\zeta^{(m+1)(k-j)}(w_1')^{\frac{N-1}{2}}
\frac{1-w_1^{-1}}{1-(w_1')^{-1}\zeta^{k-i}}.\end{align*} 
\cvd\medskip

Recall the conversion of QH tensors and the related notations of Section
\ref{ChargetoOp}.
\begin{cor}\label{wall-conv} (1) The converted tensor ${\rm XW}_N$ of $\tilde{\Hh}_N(\Ww)$ (either equal to ${\rm CW}_N$ or ${\rm MW}_N$ according to $\Ww =\Ww_C$ or $\Ww=\Ww_M$) is an endomorphism supported by
and invertible on the diagonal subspace $V$ of $\mc^N\otimes \mc^N$
with basis $e_i\otimes e_i$.
\smallskip

(2) Denote by ${\rm W}_{X,N}$ the restriction of ${\rm XW}_N$
to $V$. Then ${\rm W}_{C,N}^2=_N {\rm Id}$ and
\begin{align*} (M_N)^k_i :=  ({\rm W}_{M,N}^2)^k_i =_N
\delta_N(1+i-k)\ .
\end{align*}
\end{cor} 
\noindent {\it Proof.} The first claim in (1) is a direct consequence
of Lemma \ref{defwalls}. We compute ${\rm W}_{X,N}\circ {\rm W}_{X,N} $ by applying
in the order the Fourier transform and the conversion procedure
to:
\begin{align*}
\sum_{\alpha,\beta=0}^{N-1} \Hh_N(\Ww_m)^{\beta,k}_{\alpha,l}
\Hh_N(\Ww_m)^{j,\alpha}_{i,\beta}= N^{-2}
(w_1')^{N-1}\sum_{\alpha,\beta=0}^{N-1}
(w_1')^{-[l-k-(m+1)]_N-[j-i-(m+1)]_N} & \\
\times \delta_N(l-k-(\beta-\alpha))\ \delta_N(j-i-(\beta-\alpha)) \  .
\end{align*}
Since $w_1' = \exp(\pi\sqrt{-1}/N)$ for $\Ww_M$, in this case the
last term above is equal (mod $=_N$) to $N^{-1}\zeta^{k-l} \delta_N(l-k-(j-i))$. Now the Fourier transform is:
\begin{align*}
(\tilde{\Hh}_N(\Ww_m)\circ \tilde{\Hh}_N(\Ww_m))^{l,k}_{j,i} & =
N^{-3}\sum_{I,J,K,L=0}^{N-1}\zeta^{K-L+iI+lL-jJ-kK}\
\delta_N(L-K-(J-I))\\ 
& =
\delta_N(l-k)\ \delta_N(1+i-l)\ \delta_N(1+j-l).
\end{align*}
By restricting to $V$ this proves (2) for the $M$-wall. For the $C$-wall
the discussion is similar, using $w_1'=-1$.
\cvd

\begin{cor}\label{no-crossing} Let $K_U$ be the unknot and $L_H$ the Hopf link. We have:
$$\Hh_N(K_U)=_N 1, \quad [K_U]_N=0, \quad \Hh_N(L_H)=_N N\ .$$
\end{cor}
\noindent {\it Proof.}  The simplest way to compute $\Hh_N(K_U)$ is by
 tracing just one wall. This corresponds to a distinguished
 triangulation of $(S^3,K_U)$ with $4$ vertices (see Figure \ref{UWALL}). We get
\begin{align*} \Hh_N(K_U)= & N^{-(4-2)}\sum_{j,\alpha=0}^{N-1}N\Hh_N(W(0,-1,0,1))^{j,j+\alpha}_{j,j+\alpha}
\\ = & N^{-1}\sum_{j=0}^{N-1}N\Hh_N(W(0,-1,0,1))^{j,0}_{j,0} \\ = &
N^{-1}\sum_{j=0}^{N-1} [w_0']\omega(w_0',w_1'|0)\omega(w_0',1/w_1'|0)= N^{-1}N =1\ .\end{align*} 
We compute $[K_U]_N$ and $\Hh_N(L_H)$ by using the diagram
without crossings, equipped in the first case with two $M$-walls, and two
$C$-walls in the second. From Corollary \ref{wall-conv} (2) we deduce 
\begin{align*} [K_U]=_N & {\rm Trace}({\rm W}_{M,N}\circ {\rm W}_{M,N}) =_N 0\\
\Hh_N(L_H)=_N  & {\rm Trace}({\rm W}_{C,N}\circ {\rm W}_{C,N}) =_N N.\end{align*}
\cvd
\begin{figure}[ht]
\begin{center}
\includegraphics[width=4cm]{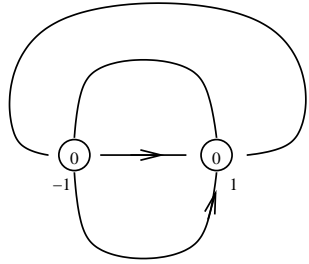}
\caption{\label{UWALL} Computation of the unknot.}
\end{center}
\end{figure} 

\subsection{QH enhanced Y-B operators: formulas}\label{YB-esplicito} Here we give explicit formulas of the Yang-Baxter operators of Theorem \ref{EYBO}. By using our favorite Yang-Baxter charge $c_0$ of Section \ref{ChargetoOp}, the braiding formulas of Lemma \ref{braiding1} depend respectively on: 
\smallskip

In ${\rm B}_N(-)$:
\begin{equation}\label{standardfcw-} 
\begin{array}{ccc} (w_0^2)' = \sqrt[N]{2},\ \ (w_0^1)'=  
(w_0^3)'= \sqrt[N]{2}\ \zeta^{(m+1)}, \ \ (w_0^4)'= \sqrt[N]{2}\
\zeta\\ (w_1^1)'= (w_1^2)' = (w_1^3)' = -1 , \ \  (w_1^4)'=-\zeta^{-1}.
\end{array}
\end{equation}
In ${\rm B}_N(+)$:
\begin{equation}\label{standardfcw+} 
\begin{array}{ccc} (w_0^2)'= (w_0^4)'= \sqrt[N]{2}, \ \
 (w_0^1)'=\sqrt[N]{2}\ \zeta^{(m+1)}, \ \ (w_0^3)'= \sqrt[N]{2}
\zeta^{-(m+1)}, \\ (w_1^1)'= (w_1^2)'= (w_1^4)' = -1, \ \
(w_1^3)' = -\zeta.
\end{array}
\end{equation}
In both cases $K_\Oo=_N N^{-1}$. By restricting ${\rm Br}_N$ to $\C^N\otimes \C^N=V\otimes V$, we get
$${\rm B}_N(-)^{j,i}_{k,l}=
N^{-1}\zeta^{(l-j)(k-i)+(m+1)(j-i-l+k)}
\frac{\omega(-1/\zeta \vert j-i)\omega(-1 \vert l-k)}{\omega(-1 \vert l-i)\omega(-1 \vert j-k)}\ ,$$
$${\rm B}_N(+)^{i,l}_{j,k}=
N^{-1}\zeta^{(l-j)(k-i)+(m+1)(l-i-j+k)}
\frac{\omega(-1/\zeta \vert j-i)\omega(-1/\zeta \vert l-k)}
{\omega(-1 \vert l-i)\omega(-1/\zeta \vert j-k)} \ .$$
The endomorphism $M_N$ has been computed in Lemma \ref{wall-conv} (2). 
Recall that 
$$({\rm W}_{C,N})^k_i=_N N^{-1}\zeta^{(m+1)(k-i)}
\frac{2}{1+\zeta^{k-i}} \ .$$
The QH R-matrices are then given by
$${\rm R}_N(-)^{j,i}_{r,s} = \sum_{k,l=0}^{N-1}({\rm W}_{C,N})^k_r
({\rm W}_{C,N})^l_s {\rm B}_N(-)^{j,i}_{k,l}\ ,$$
$${\rm R}_N(+)^{i,l}_{r,s} = \sum_{j,k=0}^{N-1}({\rm W}_{C,N})^j_r
 ({\rm W}_{C,N})^k_s {\rm B}_N(+)^{i,l}_{j,k} \ .$$
We will not need more explicit formulas.








\subsection{QH tensors and  Kashaev's R-matrices}\label{KvH3} Formulas of the Kashaev R-matrix $R_{K,N}$ can be found in \cite[(2.12) \& (2.15)]{K2} and \cite{MM}. They involve the function $\omega(1\vert [n]_N)$ introduced in Section \ref{basic}, with the residue $[n]_N\in \Ii_N$ taken as the argument, and its complex conjugate $\omega^*(1\vert [n]_N)$ . 
Referring to 
Figure \ref{QH-Rmatrix4}, the entries of $R_{K,N}^{\pm 1} = R_{K,N}(\mp)$ are given by
\begin{align*} R_{K,N}(-)_{i,j}^{l,k} = & N\zeta^{1+(l-j)(1+i-k)}\frac{\theta_N([j-i-1]_N+[l-k]_N)\theta_N([i-l]_N+[k-j]_N)} {\omega(1\vert [j-i-1]_N)\omega(1\vert [l-k]_N)\omega^*(1\vert
[k-j]_N)\omega^*(1\vert [i-l]_N)}\ ,\\
R_{K,N}(+)^{i,j}_{l,k} = & N\zeta^{(j-l)(1+i-k)}\frac{\theta_n([l-i]_n+[j-k]_n)\ 
\theta_n([i-j]_n+[k-l-1]_n)}{\omega(1\vert [j-k]_N)(\omega(1\vert [l-i]_N)
\omega^*(1\vert [k-l-1]_N)\omega^*(1\vert [i-j]_N)} \ .
\end{align*} 
Here the function $\theta:\mz\rightarrow \{0,1\}$ is defined by
\begin{equation}\label{thetadef1}
 \theta_N(n) = \left\lbrace \begin{array}{ll} 1 & \mathrm{if}\ N>n\geq
 0\\ 0 & \textrm{otherwise.}\end{array}\right.
\end{equation}
The same formulas hold for every $n>1$, not necessarily when $n=N$ is odd. 
\begin{remark}\label{sign-or}{\rm We have exchanged the roles of $R_{K,N}(+)$ and $R_{K,N}(-)$ with respect to \cite{MM}, so that we deal with the mirror image invariant $<\bar{L}>_N$ rather than $<L>_N$.}
\end{remark}
The next result is elementary. It provides various caracterizations of the non zero entries of the Kashaev R-matrix $R_{K,N}$. 
\begin{lem}\label{convstates}
{\rm Let $i$, $j$, $k$, $l \in \Ii_N$. The following properties are
equivalent:}
\begin{enumerate}
\item[(i)] $[j-i-1]_N+[l-k]_N+[i-l]_N+[k-j]_N=N-1$;
\item[(ii)] $[j-i-1]_N+[l-k]_N<N$ and $[i-l]_N+[k-j]_N<N$;
\item[(iii)] $l\leq i< j\leq k$, or $i<j\leq k\leq l$, or $j\leq k\leq
l\leq i$, or $k\leq l\leq i<j$;
\item[(iv)] The roots of unity $\zeta^i$, $\zeta^j$, $\zeta^k$ and
$\zeta^l$ are positively cyclically ordered on $S^1$, and $\zeta^i \ne
\zeta^j$ when $i < \max(j,k,l)$.
\end{enumerate}
\end{lem}
In order to relate $R_{K,N}$ to QH tensors it is convenient to deal with the
tensors ${\rm R}_N(\sigma_0,\sigma_1)$ rather than the QH R-matrices. As in Section \ref{R++} consider
$$ {\rm R}_N(+,-) =_N ({\rm W}_{C,N}\otimes {\rm id} )\circ {\rm B}_N(-)^{-1} \circ
({\rm id}  \otimes {\rm W}_{C,N}) \ .$$
First we compute the inverse braiding ${\rm B}_N(-)^{-1}$ by
using the second formula of Lemma \ref{braiding1}, specialized to our
favourite Yang-Baxter charge $c_0$, as in \eqref{standardfcw-}. In such a case
we have $\bar{K}_\mathcal{O} =_N N^{-1}\ 2^{\frac{1-N}{N}}$ and
\begin{equation}\label{Ostandard}
({\rm B}_N(-)^{-1})^{l,k}_{i,j}=
N^{-1}\zeta^{(m+1)(i+l-k-j)+(l-j)(i-k)}\frac{\omega(-1/\zeta\vert
j-k)\omega(-1/\zeta\vert l-i)}{\omega(-1 \vert
j-i-1)\omega(-1 \vert l-k) }\ .
\end{equation}
Put
$$r(x)_{i,j}^{l,k}:= 
N[x]^2\zeta^{(l-j)(i-k)}\frac{\omega(x/\zeta \vert j-k)
\omega(x/\zeta \vert l-i)}{\omega(x \vert
j-i-1)\omega(x \vert l-k)}$$
and
$$({\rm B}_N(-,x)^{-1})^{l,k}_{i,j}:= \zeta^{(m+1)(i+l-k-j)}r(x)_{i,j}^{l,k}\ .$$
Clearly $({\rm B}_N(-)^{-1})^{l,k}_{i,j}=  
({\rm B}_N(-,-1)^{-1})^{l,k}_{i,j}$. 
Moreover, it is easy to check that
\begin{equation}\label{invunit} \omega(x\zeta^{-1}\vert n) = \left\lbrace \begin{array}{lll} \frac{1-x^N}{\omega^*(x\vert N-n)} & \mbox{if}& n\in \Ii_N
\ ,\\ \\ \frac{1}{\omega^*(x\vert -n)} & \mbox{if}& n\in -\Ii_N
\ ,\end{array}\right.\end{equation}
where, abusing of notations, we set $\omega^*(u\vert n):=(\omega(u^*\vert n))^*$, taking the complex conjugate of the argument and value. That is, $$\omega^*(u\vert n) = \prod_{j=1}^n (1-u\zeta^{-j})\ , \quad n \in \Ii_N.$$ Note that $\omega(1\vert n)^{-1}=0$ for all $n<0$, and $j-k, j-i, l-i, l-k \in -\Ii_N \cup \Ii_N$. Hence 
$$r(x)_{i,j}^{l,k} = N[x]^2\zeta^{(l-j)(i-k)}\frac{1-x^N}{\omega(x\vert [j-i-1]_N)\omega(x\vert [l-k]_N)\omega^*(x\vert
[k-j]_N)\omega^*(x\vert [i-l]_N)}$$ if $\theta_N([j-i-1]_N+[l-k]_N)\theta_N([i-l]_N+[k-j]_N)=1,$ and the same formula times some positive power of $1-x^N$ otherwise. In particular, at $x=1$ the entry $r(x)_{i,j}^{l,k}/(1-x^N)$ is well-defined for all state indices, and non vanishing exactly under the conditions of Lemma \ref{convstates}. By comparing with the Kashaev R-matrix we find: 
\begin{equation}\label{KmatHO} R_{K,N}(-)_{i,j}^{l,k} = \zeta^{1+l-j}\ \left(\frac{r(x)_{i,j}^{l,k}}{1-x^N}\right)_{x=1}\ . 
\end{equation}
Let us complete now the computation of $ {\rm R}_N(+,\pm)$.  Set
\begin{equation}\label{hmat}
h(x,\alpha)_i^k = \zeta^{\alpha(k-i)} [x\zeta^{k-i}],\quad x\in \mc,
\alpha\in \mz\ .
\end{equation}
By the wall computation we know that
\begin{align*}
({\rm W}_{C,N})^k_i=_N  \left(x^{\frac{N-1}{2}}h(x,m+1)_i^k\right)_{x=-1} \ .
\end{align*}
A direct substitution gives
\begin{equation}\label{red1fin}
{\rm R}_N(+,-)_{I,j}^{l,K} = \zeta^{(m+1)(I-j+l-K)} \sum_{i,k=0}^{N-1} h(-1,1)^{i}_{I}\
r(-1)^{l,k}_{i,j} \ h(-1,1)^{K}_{k}\ ,
\end{equation}
\begin{equation}\label{red2fin}
{\rm R}_N(+,+)_{i,J}^{L,k} = \zeta^{(m+1)(i-J+L-k)} \sum_{j,l=0}^{N-1}
h(-1,0)^{j}_{J}\
 r(-1)^{l,k}_{i,j}\ h(-1,0)^{L}_{l}\ .
\end{equation}
\begin{lem}\label{conjmat} We have
\begin{align}
h(x,\alpha)h(y,\alpha) = & h(xy,\alpha)\label{hmatrix}\\
(h(y,1)\otimes {\rm Id})r(x)({\rm Id} \otimes h(1/y,1))= & \frac{r(xy)}{1-(xy)^N}\ \label{haction1}\\
({\rm Id} \otimes h(y,0))r(x)(h(1/y,0)\otimes {\rm Id})= &
\frac{r(x/y)}{1-(x/y)^N}.\label{haction2}
\end{align}
 \end{lem}
Note that \eqref{hmatrix} shows that the map $x\mapsto h(x,\alpha)$ defines a linear representation of the multiplicative group $\mc^*$  (compare with \cite[(6.12)--(6.14)]{K3}).  By taking $x=y=-1$ in the lemma and combining \eqref{KmatHO}, \eqref{red1fin} and \eqref{red2fin} we get: 
\begin{prop}\label{QHRmatrix} We have 
\begin{align*}
{\rm R}_N(+,\pm)^{l,k}_{i,j} = & 
\zeta^{(m+1)(i-j+l-k)}\ \left(\frac{r(x)_{i,j}^{l,k}}{1-x^N}\right)_{x=1}\\ 
R_{K,N}(-)^{l,k}_{i,j} = & \zeta^{1+(m+1)(l+k-i-j)}\
 {\rm R}_N(+,\pm)^{l,k}_{i,j}\ .
\end{align*}
In particular ${\rm R}_N(+,+) = {\rm R}_N(+,-)$.
\end{prop}
\noindent {\it Proof of Lemma \ref{conjmat}.} Recall that given $N$-periodic functions $g_1$, $g_2:\mz
\rightarrow \mc$ we have the Poisson formula
\begin{equation}\label{Poissonf}
 \sum_{n=0}^{N-1} g_1(n)g_2(n) = N^{-1} \sum_{n=0}^{N-1} \tilde
 g_1(n)\tilde g_2(-n)
\end{equation}
where $\textstyle \tilde g_i(n) = \sum_{\sigma=0}^{N-1}
\zeta^{n\sigma}g_i(\sigma)$ is the (unnormalized) Fourier transform of $g_i$. As in Lemma \ref{defwalls}, we compute that for fixed $x$, $\alpha$
and $j$ the functions $g_1(i) = h(x,\alpha)_i^j$ and $g_2(i) =
h(x,\alpha)^i_j$ satisfy
\begin{displaymath}
 \tilde g_1(i) = \zeta^{ij}x^{N-1-[\alpha-i-1]_N}\quad ,\quad \tilde
 g_2(i) = \zeta^{ij}x^{N-1-[\alpha+i-1]_N}.
\end{displaymath}
By \eqref{Poissonf} we deduce 
\begin{align*}
\left(h(x,\alpha)h(y,\alpha)\right)_k^j = & N^{-1}
\sum_{i=0}^{N-1}\zeta^{i(j-k)} (xy)^{N-1-[\alpha-i-1]_N}\\ = & N^{-1}
(xy)^{N-1}\zeta^{(j-k)(\alpha-1)}\
\frac{1-(xy)^{-N}}{1-(xy\zeta^{j-k})^{-1}}.
\end{align*}
This proves \eqref{hmatrix}. As for \eqref{haction1}-\eqref{haction2}, we consider the function 
\begin{displaymath}
 F\left (\begin{array}{cc} x & u \\ y & v \end{array} \bigg | z
 \right) = \sum_{\sigma=0}^{N-1}\frac{\omega( y \vert \sigma)\omega(v\vert
 \sigma)}{\omega(x\vert
 \sigma)\omega(u\vert \sigma)} z^\sigma\ ,
\end{displaymath}
where, to ensure that the summand is $N$-periodic with respect to
$\sigma$, we assume that
\begin{displaymath}
 z^N = \frac{(1-x^N)(1-u^N)}{(1-y^N)(1-v^N)}\ .
\end{displaymath}
We are going to use a symmetry relation satisfied by $F$ (see \cite{KMS}, Appendix). Let
$\xi$ be such that
\begin{displaymath}
 \xi^N  = \frac{1-x^N}{1-y^N}
\end{displaymath}
and put
\begin{displaymath}
 g_1(\sigma) = \frac{\omega(y\vert \sigma)}{\omega(x\vert \sigma)}
 \xi^\sigma \quad ,\quad g_2(\sigma) = \frac{\omega(v\vert \sigma)}{\omega(u\vert
 \sigma)} (z/\xi)^\sigma.
\end{displaymath}
By using equation \eqref{auto1} we find
\begin{align*}
 \tilde g_1(\sigma) = & f(x,y\vert \xi\zeta^\sigma) = f(x,y\vert
\xi)\ x^{-\sigma}\ \frac{\omega(\xi\zeta^{-1}\vert \sigma)}{\omega(y\xi x^{-1}\vert
\sigma)}.\end{align*} Similarly, with Remark \ref{first-rem} (2) and
$\frac{1-(z/\xi)^N}{1-(zv/u\xi)^N} = u^N$ we get
\begin{align*}
\tilde g_2(-\sigma) 
= & f(u,v\vert
z\xi^{-1})\ (v\zeta)^\sigma\ \frac{\omega(u\xi(zv\zeta)^{-1}\vert \sigma)}{\omega(\xi z^{-1}\vert
\sigma)}\ .
\end{align*}
Hence, from \eqref{Poissonf} we deduce  
\begin{align}
F\left (\begin{array}{cc} x & u \\ y & v \end{array}\bigg\vert z
\right) = & N^{-1}f(x,y\vert \xi)f(u,v\vert z\xi^{-1})\notag \\ & \hspace{1.5cm}\times \ 
\sum_{\sigma=0}^{N-1} \frac{\omega(\xi\zeta^{-1}\vert \sigma)}{\omega(y\xi x^{-1}\vert
\sigma)}\frac{\omega(u\xi(zv\zeta)^{-1}\vert \sigma)}{\omega(\xi
z^{-1}\vert \sigma)}(v\zeta/x)^\sigma \notag\\ = & N^{-1} \ f(x,y\vert \xi)\ f(u,v\vert
z\xi^{-1}) \ F\left (\begin{array}{cc} y\xi/x & \xi/z \\ \xi/\zeta &
u\xi/vz\zeta \end{array}\bigg\vert v\zeta/x \right)\ .\label{transF}
\end{align}
Now, consider the left hand side of \eqref{haction1}. By Remark \ref{first-rem} (2) we have
\begin{displaymath}
 h(y;\alpha)^i_j = [y]\zeta^{\alpha(i-j)}\frac{\omega(y\zeta^{-1}\vert i-j)}{\omega(y\vert
 i-j)} = [y]\zeta^{(\alpha-1)(i-j)}\frac{\omega((y\zeta)^{-1}\vert j-i)}{\omega(y^{-1}\vert j-i)}\ .
\end{displaymath}
Then
\begin{equation}\label{firstform}
\left((h(y,1)\otimes {\rm Id})r(x)\right)_{I,j}^{l,k} = N\ S_1\
[x]\zeta^{-k(l-j)}\frac{\omega(x/\zeta\vert j-k)}{\omega(x\vert
l-k)}
\end{equation}
where
\begin{align}
 S_1 := & [x][y]\ \sum_{i=0}^{N-1} \zeta^{i(l-j)}\
\frac{\omega(x\zeta^{-1}\vert 
l-i)}{\omega(x\vert j-i-1)}\
 \frac{\omega((y\zeta)^{-1}\vert I-i)}{\omega(y^{-1}\vert I-i)} \notag \label{secform} \\
     = & [x][y]\ \frac{\omega(x\zeta^{-1}\vert [l-I]_N)}{\omega(x\vert
[j-I-1]_N)}\ \zeta^{I(l-j)}\ F\left
(\begin{array}{cc} x\zeta^{j-I-1} & y^{-1} \\ x\zeta^{l-I-1} &
(y\zeta)^{-1} \end{array}\bigg\vert \zeta^{j-l} \right)\ .\notag
\end{align}
(We use Remark \ref{first-rem} (1) in the last equality). From \eqref{transF} with $\xi=1$ we deduce that 
\begin{multline}
S_1 = N^{-1}[x][y]\ \frac{\omega(x\zeta^{-1}\vert [l-I]_N)}{\omega(x\vert
[j-I-1]_N)}\ \zeta^{I(l-j)}\
f(x\zeta^{j-I-1}, x\zeta^{l-I-1}\vert 1) \\ \times \ f(y^{-1}, (y\zeta)^{-1}
\vert \zeta^{j-l}) \ F\left
(\begin{array}{cc} \zeta^{l-j} & \zeta^{l-j} \\
\zeta^{-1} & \zeta^{l-j} \end{array}\bigg\vert
(xy\zeta^{j-I-1})^{-1} \right).\notag
\end{multline}
From Lemma \ref{algres}, the identity \eqref{auto2}, and
\begin{displaymath}
 f(x,y\zeta\vert z) = \frac{x-y\zeta}{(1-y\zeta)(x-yz\zeta)}\ f(x,y\vert z),
\end{displaymath}
we deduce
\begin{align*}
 f(x\zeta^{j-I-1}, x\zeta^{l-I-1}\vert 1) = & [x]^{-1}\
 \frac{\omega(x\vert [j-I-1]_N)}{\omega(x\zeta^{-1}\vert [l-I]_N)}\\
 f(y^{-1}, (y\zeta)^{-1} \vert \zeta^{j-l}) = &
 y^{1-N+[j-l-1]_N}[y^{-1}]^{-1}\\ F\left (\begin{array}{cc}
 \zeta^{l-j} & \zeta^{l-j} \\ \zeta^{-1} &
 \zeta^{l-j} \end{array}\bigg\vert
 (xy\zeta^{j-I-1})^{-1} \right) = &
 f(\zeta^{l-j},\zeta^{-1}\vert (xy\zeta^{j-I-1})^{-1})\\
 = & 
\omega(xy\zeta^{j-I-1}\vert
 [l-j]_N)\\ = & \frac{N[xy]}{1-(xy)^{N}}\ \frac{\omega(xy\zeta^{-1} \vert [l-I]_N)}{\omega(xy\vert
 [j-I-1]_N)}.
\end{align*}
Hence 
\begin{displaymath}
 S_1 = \frac{[xy]y^{[j-l-1]_N}}{1-(xy)^{N}}\ \zeta^{I(l-j)}\
 \frac{\omega(xy\zeta^{-1}\vert [l-I]_N)}{\omega(xy\vert [j-I-1]_N)}.
\end{displaymath}
One computes in a similar way that
\begin{align*}
S_2 := & [x] \sum_{k=0}^{N-1} \zeta^{-k(l-j)+K-k}\
\frac{\omega(x\zeta^{-1}\vert
j-k)}{\omega(x \vert l-k)}\ \frac{1-y^{-N}}{N(1-y^{-1}\zeta^{K-k})} \\ = &
\frac{[xy]y^{1-N+[l-j]_N}}{1-(xy)^{N}}\ \zeta^{-K(l-j)}\
\frac{\omega(xy\zeta^{-1}\vert 
[j-K]_N)}{\omega(xy \vert [l-K]_N)}\ .
\end{align*}
By using \eqref{invunit} and gathering terms we eventually obtain equation \eqref{haction1}:
\begin{align*}
\left((h(y,1)\otimes {\rm Id})\right. & \left. r(x) ({\rm Id} \otimes h(1/y,1))\right)_{I,j}^{l,K} =   N\ S_1 \
S_2\ \zeta^{\beta(l-j)}\\ \hspace{2cm} =& \frac{N[xy]^2}{(1-(xy)^N)^2}\ \zeta^{(I-K)(l-j)}\ \frac{\omega(xy\zeta^{-1} \vert [l-I]_N)\omega(xy\zeta^{-1}\vert
[j-K]_N)}{\omega(xy\vert
 [j-I-1]_N)\omega(xy \vert [l-K]_N)}\\ \hspace{2cm}= & \frac{r(xy)_{I,j}^{l,K}}{1-(xy)^N}\ .
\end{align*}
Equation \eqref{haction2} is proved in a similar way.\cvd


\end{document}